\documentclass[12pt]{amsart}
\usepackage{amssymb}

\newtheorem{Thm}{Theorem}[section]
\newtheorem{notation}[Thm]{Notation}
\newtheorem{problem}[Thm]{Problem}

\newtheorem{corollary}[Thm]{Corollary}
\newtheorem{lemma}[Thm]{Lemma}
\newtheorem{proposition}[Thm]{Proposition}
\newtheorem{definition}[Thm]{Definition}
\newtheorem{remark}[Thm]{Remark}
\newtheorem{conjecture}[Thm]{Conjecture}

\newtheorem{example}[Thm]{Example}
\newtheorem{theorem}[Thm]{Theorem}

\newtheorem*{kadison-singer}{Kadison-Singer Problem}
\newtheorem*{paving conjecture}{Paving Conjecture}
\newtheorem*{bourgain-tzafriri}{Bourgain-Tzafriri Conjecture}
\newtheorem*{feichtinger conjecture}{Feichtinger Conjecture}
\newtheorem*{cocktail}{Cocktail Party Problem}

\def\ldots{\mathinner{\ldotp\ldotp\ldotp}}
\def\ldots{\mathinner{\cdotp\cdotp\cdotp}}

\def \H{{\mathbb H}}
\def \Z{{\mathbb Z}}
\def \N{{\mathbb N}}
\def \D{{\mathbb D}}
\def \DD{{\mathbb D}}
\def \B(l_2){B(\ell_2)}
\def \HN{{\mathbb H}_{n}}

\def \K{{\mathbb K}}

\def \cal{\mathcal}

\def \beq{\begin{eqnarray*}}
\def \eeq{\end{eqnarray*}}

\def \< {\langle}
\def \> {\rangle}

\def \R {\mathbb{R}}

\def \p{\phi}
\def \Aj{\{A_{j}\}_{j=1}^{r}}
\def \<{\langle}
\def \>{\rangle}

\begin{document}

\title[The Kadison-Singer Problem in Mathematics
and Engineering]{The Kadison-Singer Problem in
mathematics and engineering:  A detailed account}
\author[P.G. Casazza, M. Fickus, J.C. Tremain, E. Weber
 ]{Peter G. Casazza, Matthew Fickus, Janet C. Tremain, Eric Weber}

\address{Department of Mathematics \\
University of Missouri-Columbia \\
Columbia, MO 65211}
\email{pete@math.missouri.edu}

\address{Department of Mathematics and Statistics, Air Force Institute
of Technology, Wright-Patterson Air Force Base, Ohio 45433}

\email{Matthew.Fickus@afit.edu}

\address{Department of Mathematics \\
University of Missouri-Columbia \\
Columbia, MO 65211}
\email{janet@math.missouri.edu}

\address{Department of mathematics, Iowa State University,
396 Carver Hall, Ames, IA 50011}
\email{esweber@iastate.edu}

\thanks{The first and second authors were supported by NSF DMS 0405376,
 the last author
was supported by NSF DMS 0355573.}


\maketitle

\begin{abstract}
We will show that the famous, intractible 1959 
Kadison-Singer problem in $C^{*}$-algebras
is equivalent to fundamental unsolved problems in a dozen areas
of research in pure mathematics, applied mathematics and Engineering. This 
gives all these areas common ground on which to interact
as well as explaining why
each of these areas has volumes of literature on their respective
problems without a satisfactory resolution. 
In each of these areas we will
reduce the problem to the minimum which needs to be proved to
solve their version of Kadison-Singer.  In some areas we will
prove what we believe will be the strongest results ever
available in the case that Kadison-Singer fails.
Finally, we will give some
directions for constructing a counter-example to
Kadison-Singer.
\end{abstract}

\section{Introduction}\label{Intro}
\setcounter{equation}{0}

The famous 1959 Kadison-Singer Problem \cite{KS} has defied the best efforts
of some of the most talented mathematicians of our time.

\begin{kadison-singer}[KS]
Does every pure state on the (abelian) von Neumann algebra $\D$ of
bounded diagonal operators on ${\ell}_2$ have a unique extension to
a (pure) state on $B({\ell}_2)$, the von Neumann algebra of all
bounded linear operators on the Hilbert space ${\ell}_2$?
\end{kadison-singer}

A {\bf state} of a von Neumann algebra ${\cal R}$ is a linear
functional $f$ on ${\cal R}$ for which $f(I) = 1$ and $f(T)\ge 0$
whenever $T\ge 0$ (whenever $T$ is a positive operator).
The set of states of ${\cal R}$ is a convex subset of the dual space
of ${\cal R}$ which is compact in the ${\omega}^{*}$-topology.  By the
Krein-Milman theorem, this convex set is the closed convex hull of its
extreme points.  The extremal elements in the space of states are
called the {\bf pure states} (of ${\cal R}$).
This problem arose from the very productive
collaboration of Kadison and Singer in
the 1950's which culminated in their seminal work on triangular
operator algebras.  During this collaboration, they often discussed
the fundamental work of Dirac \cite{Di} on Quantum Mechanics.  
In particular, they kept returning to one part of Dirac's work because it
seemed to be problematic.   Dirac wanted to find a
``representation'' (an orthonormal basis) for a compatible
family of observables (a commutative family of self-adjoint
operators).  On pages 74--75 of \cite{Di} Dirac states:

\vskip8pt
\hskip.5truein\vbox{\hsize4truein
``To introduce a representation in practice
\vskip8pt
(i)  \hskip4pt We look for observables which we would like to have diagonal
either because we are interested in their probabilities or for
reasons of mathematical simplicity;

(ii)  \hskip4pt We must see that they all commute --- a necessary condition
since diagonal matrices always commute;

(iii)  \hskip4pt We then see that they form a complete commuting set, and if
not we add some more commuting observables to make them into a complete
commuting set;

(iv)  \hskip4pt We set up an orthogonal representation with this commuting set
diagonal.
\vskip8pt

{\bf The representation is then completely determined} ...
 {\bf by the observables
that are diagonal} ...''}
\vskip6pt
The emphasis above was added.  Dirac then talks about finding a basis that
diagonalizes a self-adjoint operator, which is troublesome since there
are perfectly respectable self-adjoint operators which do not have
a single eigenvector.  Still, there is a {\it spectral resolution}
of such operators.  Dirac addresses this problem on pages 57-58
of \cite{Di}:
\vskip8pt
\hskip.5truein\vbox{\hsize4truein
``We have not yet considered the lengths of the basic vectors.  With an
orthonormal representation, the natural thing to do is to normalize
the basic vectors, rather than leave their lengths arbitrary, and
so introduce a further stage of simplification into the representation.
However, it is possible to normalize them only if the parameters are
continuous variables that can take on all values in a range, the basic
vectors are eigenvectors of some observable belonging to eigenvalues
in a range and are of infinite length...''}
\vskip8pt

In the case of $\DD$, the representation is
$\{e_i\}_{i\in I}$, the orthonormal basis of $l_2$.  But what happens if our
observables have ``ranges" (intervals) in their spectra?  This led Dirac to
introduce his famous $\delta$-function --- vectors of ``infinite length."
From a mathematical point of view, this is problematic.  What we need is to
replace the vectors $e_i$ by some mathematical object that is essentially the
same as the vector, when there is one, but gives us something precise and
usable when there is only a $\delta$-function.  This leads to the ``pure
states'' of $\B(l_2)$ and, in particular, the (vector) pure states $\omega_x$,
given by $\omega_x(T)=\<Tx,x\>$, where $x$ is a unit vector in $\H$.  Then,
$\omega_x(T)$ is the expectation value of $T$ in the state corresponding to
$x$.  This expectation is the average of values measured in the laboratory
for the ``observable" $T$ with the system in the state corresponding to $x$.
The pure state $\omega_{e_i}$ can be shown to be completely determined by its
values on $\DD$; that is, each $\omega_{e_i}$ has a {\it unique\/} extension
to $\B(l_2)$.  But there are many other  pure states of $\DD$.  (The family of
all pure states of $\DD$ with the $w^*$-topology is $\beta(\Z)$, the
$\beta$-compactification of the integers.)  Do these other pure states have
unique extensions?  This is the Kadison-Singer problem (KS).

By a ``complete" commuting set, Dirac means what is now called a ``maximal
abelian self-adjoint" subalgebra of $\B(l_2)$; $\DD$ is one such.  There are
others.  For example, another is generated by an observable whose``simple"
spectrum is a closed interval.  Dirac's claim, in mathematical form, is that
each pure state of a ``complete commuting set" has a unique state extension
to $\B(l_2)$.  Kadison and Singer show [37] that that is {\it not so\/} for
each complete commuting set other than $\DD$.  They also show that each pure
state of $\DD$ has a unique extension to the uniform closure of the algebra
of linear combinations of operators $T_\pi$ defined by $T_\pi e_i=e_{\pi(i)}$,
where $\pi$ is a permutation of $\Z$.

Kadison and Singer believed that KS had a negative answer.  In particular,
on page 397 of \cite{KS} they state:  ``We incline to the view
that such extension is non-unique''.

This paper is based on two fundamental principles.
\vskip8pt
\noindent {\bf Fundamental Principle I}[Weaver, Conjecture \ref{C2}]:
The Kadison-Singer Problem is a statement about partitioning
projections on finite dimensional Hilbert spaces with small diagonal into
submatrices of norms $\le 1-\epsilon$.
\vskip8pt
\noindent {\bf Fundamental Principle II}[Theorem \ref{FPII}]:
Every bounded operator on a finite dimensional Hilbert space is
a constant times a ``piece'' of a projection operator
from a larger Hilbert space.
\vskip8pt
Armed with these two basic principles, we will make a tour of
many different areas of research.  In each area we will use
Fundamental Principle II (often in disguised form)
to reduce their problem to a statement about (pieces of)
projections.  Then we will apply Fundamental Principle I
to see that their problem is equivalent to the Kadison-Singer
Problem.
\vskip8pt

This paper is a greatly expanded version of \cite{CT}.
Let us now discuss the organization of this paper.
  In Sections
2-8 we will successively
look at equivalents of the Kadison-Singer Problem
 in operator theory, frame theory, Hilbert space theory,
Banach space
theory, harmonic analysis,
time-frequency analysis and finally in engineering.  In section
9 we will address some approaches to producing a counter-example
to KS.  In Section 2 we will establish our first fundamental
principle for showing that very general problems are equivalent
to KS.  In Section 3 we introduce our ``universal language''
of frame theory and introduce our second fundamental principle
for reducing problems to KS.
In Section 4, we will show that
KS is equivalent to a fundamental result concerning inner products.
This formulation of the problem has the advantage that it can
be understood by a student one week into their first course in
Hilbert spaces.  In Section 5 we show that KS is equivalent to
the Bourgain-Tzafriri Conjecture (and in fact, a significantly
weaker form of the conjecture is equivalent to KS).  This also
shows that the Feichtinger Conjecture is equivalent to KS.
 In Section 6, we show that
a fundamental problem in harmonic analysis is equivalent to KS.
We also classify the uniform paving conjecture and the uniform
Feichtinger Conjecture.  As a consequence we will discover a
surprising new identity in the area.  In Section 7, we show that
the Feichtinger Conjecture for frames of translates is equivalent
to one of the fundamental unsolved problems in harmonic analysis.
In Section 8, we look at how KS arises naturally in various
problems in signal-processing, internet coding, coding theory and
more.
\vskip8pt
\noindent {\bf Notation for statements of problems}:  Problem A
(or Conjecture A) {\bf implies} Problem B (or Conjecture B) means
that a positive solution to the former implies a positive solution
to the latter.  They are {\bf equivalent}
if they imply each other.
\vskip8pt
\noindent {\bf Notation for Hilbert spaces}:
Throughout, ${\ell}_2 (I)$ will denote a finite or
infinite dimensional complex Hilbert space with
a fixed orthonormal basis $\{e_i\}_{i\in I}$.
If $I$ is infinite we let ${\ell}_2 = {\ell}_2 (I)$,
and if $|I|=n$ write ${\ell}_2(I) = {\ell}_2^n$ with fixed
orthonormal basis $\{e_i\}_{i=1}^{n}$.  For any Hilbert
space $\H$ we let $B(\H)$ denote the family of bounded
linear operators on $\H$.
 An $n$-dimensional subspace of ${\ell}_2(I)$ will
be denoted $\HN$.  
For an operator $T$ on any one of our Hilbert spaces,
its matrix representation with respect to our fixed orthonormal 
basis is the collection $(\langle Te_i,e_j\rangle )_{i,j\in I}$.
If $J\subset I$,
the {\bf diagonal projection}
$Q_J$ is the matrix whose entries are all zero except
for the $(i,i)$ entries for $i\in J$ which are all one.
For a matrix $A=(a_{ij} )_{i,j\in I}$ let
${\delta}(A) =\max_{i\in I} |a_{ii}|$.
\vskip8pt
\noindent{\bf A universal language}:  We are going to show that
the Kadison-Singer problem is equivalent to fundamental unsolved
problems in a dozen different areas of research in both mathematics
and engineering.  But each of these areas is overrun with technical
jargon which makes it difficult or even impossible for those outside
the field to understand results inside the field.  What we need is a
{\it universal language} for interactions between a broad spectrum
of research.  For our universal language, we have chosen the language
of {\it Hilbert space frame theory} (See Section \ref{FT})
because it is simply
stated and easily understood while being fundamental enough to quickly
pass quite technical results between very diverse areas of research.
Making it possible for researchers from a broad spectrum of research areas
to understand how their problems relate to areas they may know little
about will require certain redundancies.  That is, we will have to reprove
some results in the literature in the format of our universal language.
Also, since frame theory is our universal language, we will prove
some of the fundamental results in this area so that researchers
will have a solid foundation for reading the rest of the paper.
\vskip8pt

\noindent {\bf Acknowledgement}:  We are indebted to Richard
Kadison for numerous suggestions and helpful discussions as
well as making available to us various talks he has given on
the history of KS.

\section{Kadison-Singer in Operator Theory}\label{OT}
\setcounter{equation}{0}

A significant advance on KS was made by Anderson \cite{A} in
1979 when he reformulated KS into what is now known as the
{\bf Paving Conjecture} (See also \cite{A3,A2}).
Lemma 5 of \cite{KS} shows a connection
between KS and Paving.

\begin{paving conjecture}[PC]
For $\epsilon >0$, there is a natural number $r$ so
that for every
 natural number $n$ and
every linear operator
  $T$ on $l_2^n$
whose matrix has zero diagonal,
  we can find a partition (i.e. a {\it paving})
$\{{A}_j\}_{j=1}^r$
  of $\{1, \ldots, n\}$, such that
  $$
  \|Q_{{A}_j} T Q_{{A}_j}\|  \le  \epsilon \|T\|
  \ \ \ \text{for all $j=1,2,\ldots ,r$.}
  $$
\end{paving conjecture}

It is important that $r$ not depend on $n$ in PC.
We will say that an arbitrary operator $T$ satisfies PC if
$T-D(T)$ satisfies PC where $D(T)$ is the diagonal of $T$.

\begin{remark}\label{R1}
There is a standard technique for turning finite dimensional
results into infinite dimensional ones and vice-versa.
We will illustrate this technique here by showing that
PC is equivalent to PC for operators on $\ell_2$ (which is
a known result).  After
this, we will move freely between these cases for our later
conjectures without proving that they are equivalent.
\end{remark}

We can use an abstract compactness argument for proving
this result, but we feel that the following argument is more
illuminating.  We start with a limiting method for
increasing sequences of partitions given in \cite{CCLV}.
Since the proof is short we include it for completeness.

\begin{proposition}\label{prop1}
Fix a natural number $r$ and assume for every natural number
$n$ there is a partition $\{A_j^n\}_{j=1}^{r}$ of $\{1,2,\ldots ,n\}$.
There exist natural numbers $\{k_1<k_2<\cdots \}$ so that if
$m\in A_j^{k_m}$ for some $1\le j\le r$ then $m\in A_j^{k_{\ell}}$,
for all $\ell \ge m$.  Hence, if $A_j = \{m\ |\ m\in A_j^{k_m}\}$ then
\begin{enumerate}
\item $\{A_j\}_{j=1}^{r}$ is a partition of $\N$.
\vspace*{.04cm}
\item If $A_j = \{m_1 < m_2 < \cdots \}$, then for all natural numbers
$\ell$ we have $\{m_1 ,m_2 ,\ldots ,m_{\ell}\}\subset A_j^{k_{m_{\ell}}}$.
\end{enumerate}
\end{proposition}

{\it Proof}:
For each natural number n, 1 is in one of the sets $\{A_j^n\}_{j=1}^{r}$.
Hence, there are natural numbers $n_1^1<n_2^1< \cdots$ and a $1\le j\le r$
so that $1\in A_j^{n_i^1}$ for all $i\in \N$.  Now, for every natural number
$n_i^1$, 2 is in one of the sets $\{A_j^{n_i^1}\}_{j=1}^{r}$.
Hence, there is a subsequence $\{n_i^2 \}$ of $\{n_i^1\}$ and a $1\le j\le
r$ so that $2\in A_j^{n_i^2}$, for all $i\in \N$.  Continuing by induction,
for all $\ell \in \N$ we get a subsequence $\{n_i^{\ell +1}\}$ of
$\{n_i^{\ell}\}$ and a $1\le j\le r$ so that
${\ell} +1\in A_{j}^{n_j^{\ell +1}}$, for all $i\in \N$.  Letting
$k_i = n_i^i$ for all $i\in \N$ gives the
conclusion of the proposition.
\qed

\begin{theorem}\label{IPC}
The Paving Conjecture is equivalent to the Paving Conjecture for
operators on ${\ell}_2$.
\end{theorem}

{\it Proof}:
Assume PC holds for operators on ${\ell}_2^n$.
 Let $T = (a_{ij})_{i,j=1}^{\infty}$ be a
bounded linear operator on $\ell_2$.  Fix $\epsilon >0$.
By our assumption, for every natural
number $n$ there is a partition $\{A_j^n\}_{j=1}^{r}$ of
$\{1,2,\ldots ,n\}$ so that if $T_n = (a_{ij})_{i,j=1}^{n}$ then
for all $j=1,2,\ldots ,r$
$$
\|Q_{A_j^n}T_nQ_{A_j^n}\|\le \frac{\epsilon}{2}\|T_n\| \le \frac{\epsilon}{2}
\|T\|.
$$
Let $\{A_j\}_{j=1}^{r}$ be the partition of $\N$ given in Proposition
\ref{prop1}.  Fix $1\le j\le r$, let $A_j = \{m_1<m_2< \cdots \}$,
and for all $\ell \in \N$ let $Q_{\ell} = Q_{I_{\ell}}$ where
$I_{\ell} =\{m_1,m_2 ,\ldots ,m_{\ell}\}$.  Fix $f\in {\ell}_2(\N)$.
For all large $\ell \in \N$ we have:
\begin{eqnarray*}
\|Q_{A_j}TQ_{A_j}(f)\| &\le&
2\|Q_{\ell}Q_{A_{j}}TQ_{A_j}Q_{\ell}(f)\| \\
&=& 2 \|Q_{\ell}Q_{A_j}^{k_{m_{\ell}}}T_{k_{m_{\ell}}}Q_{A_j}^{k_{m_{\ell}}}
Q_{\ell}(f)\| \\
&\le& 2 \|Q_{A_j}^{k_{m_{\ell}}}T_{k_{m_{\ell}}}Q_{A_j}^{k_{m_{\ell}}}\|
\|Q_{\ell}(f)\|\\
&\le& 2 \frac{\epsilon}{2}\|T\|\|f\|= \epsilon \|T\|\|f\|.
\end{eqnarray*}
Hence, $\|Q_{A_j}TQ_{A_j}\|\le \epsilon \|T\|.$

Conversely, assume PC holds for operators on
${\ell}_2$.  We assume that PC fails for operators on ${\ell}_2^n$
and get a contradiction.  If $(1)$ fails, a little thought will
yield that there must be an $\epsilon>0$, a partition $\{I_n\}_{n=1}^{\infty}$
of $\N$ into finite subsets, 
operators $T_n:{\ell}_2(I_n)\rightarrow {\ell}_2(I_n)$ with
$\|T_n\|=1$ and for every partition $\{A_j^n\}_{j=1}^{n}$ of $I_n$
there exists a $1\le j\le n$ so that
$$
\|Q_{A_j^n}T_nQ_{A_j^n}\|\ge \epsilon.
$$
Let
$$
T = \bigoplus_{n=1}^{\infty}T_n : \left(\bigoplus_{n=1}^{\infty}
{\ell}_2(I_n) \right)_{{\ell}_2} \rightarrow
\left(\bigoplus_{n=1}^{\infty}{\ell}_2(I_n) \right)_{{\ell}_2}.
$$
Then, $\|T\|=\sup_n\ \|T_n\|=1$.  By (2), there is a partition
$\{A_j\}_{j=1}^{r}$ of $\N$ so that for all $j=1,2,\ldots ,r$
$$
\|Q_{A_j}TQ_{A_j}\|\le \epsilon.
$$
For every $n\in \N$ and every $j=1,2,\ldots ,r$ let
$A_j^n = A_j \cap I_n$.  Then, $\{A_j^n\}_{j=1}^{r}$ is a partition
of $I_n$.  Hence, for every $j=1,2,\ldots ,r$ we have
$$
\|Q_{A_j^n}T_n Q_{A_j^n}\| = \|Q_{A_j^n}TQ_{A_j^n}\| \le
\|Q_{A_j}TQ_{A_j}\| \le \epsilon.
$$
If $n\ge r$, this contradicts our assumption about $T_n$.
\qed

It is known \cite{BHKW} that the class of
operators satisfying PC (the {\bf pavable operators})
is a closed subspace of $B({\ell}_2)$.
The only large
classes of operators which have been shown to be pavable are
``diagonally dominant'' matrices \cite{BCHL,BCHL2,BHKW,G},
matrices with all entries real and positive \cite{HKW} 
and Toeplitz operators
over Riemann integrable functions (See also \cite{HKW2} and Section \ref{HA}).
Also, in \cite{BHKW2} there is an analysis of the paving problem for
certain Schatten $C_p$-norms. 
We strongly recommend that everyone read the argument of
Berman, Halpern, Kaftal and Weiss \cite{BHKW} showing that matrices with
positive entries satisfy PC.  This argument is a fundamental principle
concerning decompositions of matrices which has applications across the
board --- here, you will see it used in the proof of
Theorem \ref{TP1}, and it was {\it vaguely} used in producing a generalization of
the Rado-Horn Theorem \cite{CKS} (See Theorem \ref{radohorn}).
We next note that in order to verify PC, 
it suffices to show that PC holds for any one of 
your favorite classes of operators.

\begin{theorem}\label{TOT1}
The Paving Conjecture has a positive solution if any one of the following
classes satisfies the Paving Conjecture:
\begin{enumerate}
\item Unitary operators.
\vspace*{.04cm}
\item Orthogonal projections.
\vspace*{.04cm}
\item Positive operators.
\vspace*{.04cm}
\item Self-adjoint operators.
\vspace*{.04cm}
\item Gram matrices $(\langle f_i,f_j\rangle )_{i,j\in I}$ where
$T:{\ell}_2(I) \rightarrow {\ell}_2(I)$ is a bounded
linear operator, and $Te_i = f_i$,  $\|Te_i \| = 1$
for all $i\in I$.
\vspace*{.04cm}
\item Invertible operators (or invertible operators with zero diagonal).
\end{enumerate}
\end{theorem}

{\it Proof}:
$(1)$:  This is immediate from the fact that every bounded operator
is a multiple of a sum of three unitary operators \cite{Cas1}.

$(2)$:  This follows from the Spectral Theorem (or see Fundamental
Principle II:  Theorem \ref{FPII}).

$(3),(4)$:  Since $(3)$ or $(4)$ immediately implies $(2)$.

$(5)$:  We will show that $(5)$ implies a positive solution
to the Bourgain-Tzafriri Conjecture (See Section \ref{BST})
and hence to PC by Theorem \ref{TBST1}.  Given $T:{\ell}_2(I)
\rightarrow {\ell}_2(I)$ with $\|Te_i\|=1$ for all $i\in I$,
let $G = (\langle Te_i,Te_j \rangle )_{i,j\in I}$.  By $(5)$,
there is a partition $\{A_j\}_{j=1}^{r}$ of $I$ which paves
the Gram operator.  Hence, for
all $j=1,2,\ldots ,r$ and all  $f=\sum_{i\in A_j}a_ie_i$ we have
\begin{eqnarray*}
\|\sum_{i\in A_j}a_i Te_i\|^2 &=& \langle \sum_{i\in A_j}a_i Te_i ,
\sum_{k\in A_j}a_k Te_k \rangle \\
&=&
\sum_{i\in A_j}|a_i|^2\|Te_i\|^2
+ \sum_{i\not= k \in A_j}a_i \overline{a_k}\langle Te_i,Te_j \rangle\\
&=& \sum_{i\in A_j}|a_i|^2 + \langle Q_{A_j}(G-D(G))Q_{A_j}f,f\rangle \\
&\ge& \sum_{i\in A_j}|a_i|^2 - \|Q_{A_j}(G-D(G))Q_{A_j}\|\|Q_{A_j}f\|^2\\
&\ge& \sum_{i\in A_j}|a_i|^2 - \epsilon \sum_{i\in A_j}|a_i|^2\\
&=& (1- \epsilon )\sum_{i\in A_j}|a_i|^2.
\end{eqnarray*}

Hence, the Bourgain-Tzafriri Conjecture holds (See section
\ref{BST}). Now we need to jump ahead to Theorem \ref{BSTCT} to see
that the proof of BT implies KS is done from the definition and
does not need any theorems developed between here and there.

$(6)$:  Given an operator $T$, $(\|T\|+1)I + T$ is invertible and
if it is pavable then so is $T$.  For the second statement, given
an operator $T$, let $S = T+(\|T\|^2+2)U$ where $U=(b_{ij})_{i,j\in I}$
 is the unitary
matrix given by the bilateral shift on $\N$
(the wrap-around shift on ${\ell}_2^n$ if $|I|=n$).
Then, $S-D(S)$ is invertible and has zero diagonal.  By $(6)$, for
$0<\epsilon <1$ there is a partition $\{A_j\}_{j=1}^{r}$ of $I$
so that for all $j=1,2,\ldots ,r$ we have
$$
\|Q_{A_j}(S-D(S))Q_{A_j}\| \le \epsilon.
$$
Note that for any $i\in I$, if $i\in A_j$ then $i+1 \notin A_j$,
since otherwise:
$$
|(Q_{A_j}(S-D(S))Q_{A_j}e_{i+1})(i)| = |\langle Te_i,Te_{i+1}\rangle
+ (\|T\|^2+2)| \ge 1.
$$
Hence, $\|Q_{A_j}(S-D(S))Q_{A_j}\|\ge 1$, which contradicts our
paving of $S-D(S)$.  It follows that
$$
Q_{A_j}(S-D(S))Q_{A_j} = Q_{A_j}(T-D(T))Q_{A_j},\ \ \mbox{for all
$j=1,2,\ldots ,r$}.
$$
So, our paving of $S$ also paves $T$.
\qed

Akemann and Anderson \cite{AA} showed that the following conjecture implies
KS.

\begin{conjecture}\label{C1}
There exists $0<\epsilon ,\delta <1$ with the following property:
for any orthogonal projection $P$ on ${\ell}_2^n$ with
${\delta}(P)\le \delta$, there is a diagonal projection $Q$
such that $\|QPQ\|\le 1-\epsilon$ and $\|(I-Q)P(I-Q)\|\le
1- \epsilon$.
\end{conjecture}

It is important that $\epsilon, \delta$ are independent of $n$
in Conjecture \ref{C1}.
It is unknown if KS implies Conjecture \ref{C1}.
Weaver \cite{W2} showed that a conjectured strengthening of
Conjecture \ref{C1} fails.

Recently, Weaver \cite{W} provided important insight
into KS by showing that a slight weakening of Conjecture \ref{C1}
will produce a conjecture equivalent to KS.  This is our
first Fundamental Principle.

\begin{conjecture}[Fundamental Principle I:  Weaver]\label{C2}
There exist universal constants $0<\delta, \epsilon <1$ and
$r\in \N$ so that for all $n$ and all orthogonal projections $P$
on ${\ell}_2^n$
with ${\delta}(P)\le \delta$, there is a paving $\{A_j\}_{j=1}^{r}$
of $\{1,2,\ldots ,n\}$ so that $\|Q_{A_j}PQ_{A_j}\|\le
1-\epsilon$, for all $j=1,2,\ldots ,r$.
\end{conjecture}

This needs some explanation since there is nothing in \cite{W}
which looks anything like Conjecture \ref{C2}. In \cite{W},
Weaver introduces what he calls ``Conjecture $\mathrm{KS_r}$''
(See Section \ref{FT}).
A careful
examination of the proof of Theorem 1 of \cite{W} reveals that
Weaver shows Conjecture $\mathrm{KS_r}$ implies Conjecture \ref{C2}
which in turn implies KS which (after the theorem is proved)
is equivalent to $\mathrm{KS_r}$.  We will see in Section \ref{FT}
(Conjecture \ref{W4}, Theorem \ref{W5}) that we may assume
$\|Pe_i\| = \|Pe_j\|$ for all $i,j = 1,2,\ldots ,n$ in Conjecture
\ref{C2} with a small restriction on the $\epsilon >0$.

\section{Frame Theory:  The Universal Language}\label{FT}
\setcounter{equation}{0}

A family of vectors $\{f_i\}_{i\in I}$ in a Hilbert space $\H$
is a {\bf Riesz basic sequence} if there are constants $A,B>0$ so that
for all scalars $\{a_i\}_{i\in I}$ we have:
\begin{displaymath}
A\sum_{i\in I}|a_i |^2 \le \|\sum_{i\in I}a_i f_i \|^2 \le
B\sum_{i\in I}|a_i|^2.
\end{displaymath}
We call $\sqrt{A},\sqrt{B}$ the {\bf lower and upper Riesz basis
bounds} for $\{f_i\}_{i\in I}$.  If the Riesz
basic sequence $\{f_i\}_{i\in I}$ spans $\H$
we call it a {\bf Riesz basis} for $\H$.  So $\{f_i\}_{i\in I}$
is a Riesz basis for $\H$ means there is an orthonormal basis
$\{e_i\}_{i\in I}$ so that the operator $T(e_i )=f_i$ is invertible.
In particular, each Riesz basis is {\bf bounded}.  That is,
$0 < \inf_{i\in I}\|f_i\| \le \sup_{i\in I}\|f_i\| < \infty$. 

Hilbert space frames were introduced by Duffin and Schaeffer
\cite{DS} to address some very deep problems in nonharmonic
Fourier series (see \cite{Y}).
A family $\{f_i \}_{i\in I}$ of elements of a (finite
or infinite dimensional) Hilbert space
$\H$ is called a {\bf frame} for $\H$ if
there are constants $0<A\le B < \infty$ (called the
{\bf lower and upper frame bounds}, respectively)
so that for all $f\in \H$
\begin{equation}\label{E5}
A\|f\|^2 \le \sum_{i\in I}|\langle f,f_i \rangle |^2
\le
B\|f\|^2.
\end{equation}
If we only have the right hand inequality in Equation \ref{E5}
 we call $\{f_i\}_{i\in I}$
a {\bf Bessel sequence with Bessel bound B}.
  If $A=B$, we call this an
$A$-{\bf tight frame} and if $A=B=1$, it is called a
{\bf Parseval frame}.  If all the frame elements have the same
norm, this is an {\bf equal norm} frame and if the
frame elements are of unit norm, it is a {\bf unit norm
frame}.  It is immediate that $\|f_i\|^2\le B$.  If also
inf $\|f_i\|>0$, $\{f_i\}_{i\in I}$ is a {\bf bounded frame}.
The numbers $\{\langle f,f_i\rangle\}_{i\in I}$ are
 the {\bf frame coefficients} of the vector $f\in \H$.
If $\{f_i \}_{i\in I}$ is a Bessel sequence,
the {\bf synthesis operator} for $\{f_i\}_{i\in I}$
is the bounded linear operator $T:{\ell}_{2}(I) \rightarrow
\H$ given by $T(e_i )= f_i$ for all $i\in I$.
The
{\bf analysis operator} for $\{f_i\}_{i\in I}$
is $T^{*}$ and satisfies: $T^{*}(f) =
\sum_{i\in I} \langle f,f_i \rangle e_i$.  In particular,
$$
\|T^{*}f\|^2 = \sum_{i\in I}|\langle f,f_i\rangle |^2,\ \
\mbox{for all $f\in \H$},
$$
and hence the smallest Bessel bound for $\{f_i\}_{i\in I}$
equals $\|T^{*}\|^2$.
Comparing this to Equation \ref{E5} we have:

\begin{theorem}\label{FTTT}
Let $\H$ be a Hilbert space and $T:{\ell}_2(I)\rightarrow
\H$, $Te_i = f_i$ be a bounded linear operator.  The following
are equivalent:
\begin{enumerate}
\item $\{f_i\}_{i\in I}$ is a frame for $\H$.
\vspace*{.04cm}
\item The operator $T$ is bounded, linear, and onto.
\vspace*{.04cm}
\item The operator $T^{*}$ is an (possibly into) isomorphism.
\end{enumerate}
Moreover, if $\{f_i\}_{i\in I}$ is a Riesz basis, then the Riesz
basis bounds are $\sqrt{A},\sqrt{B}$ where $A,B$ are the frame
bounds for $\{f_i\}_{i\in I}$.
\end{theorem}

It follows that a Bessel sequence
is a Riesz basic sequence if and only if $T^{*}$ is onto.
   The {\bf frame
operator} for the frame is the positive, self-adjoint invertible
operator $S=TT^{*}:\H \rightarrow \H$.  That is,
$$
Sf = TT^{*}f = T\left ( \sum_{i\in I}\langle f,f_i\rangle e_i
\right ) = \sum_{i\in I}\langle f,f_i\rangle Te_i =
\sum_{i\in I}\langle f,f_i\rangle f_i.
$$
In particular,
$$
\langle Sf,f\rangle = \sum_{i\in I}|\langle f,f_i \rangle|^2.
$$
It follows that $\{f_i\}_{i\in I}$ is a frame with frame
bounds $A,B$ if and only if $A \cdot I \le S \le B \cdot I$.
So $\{f_i\}_{i\in I}$ is a Parseval frame if and only if $S=I$.
{\bf Reconstruction} of vectors in $\H$ is achieved via
the formula:
\begin{eqnarray*}
f &=& SS^{-1}f = \sum_{i\in I}\langle S^{-1}f,f_i \rangle f_i \\
&=& \sum_{i\in I}\langle f,S^{-1}f_i \rangle f_i \\
&=& \sum_{i\in I}\langle f,f_i \rangle S^{-1}f_i \\
&=& \sum_{i\in I}\langle f,S^{-1/2}f_i \rangle S^{-1/2}f_i.
\end{eqnarray*}
It follows that $\{S^{-1/2}f_i\}_{i\in I}$ is a Parseval frame
{\it equivalent} to $\{f_i\}_{i\in I}$.  Two sequences
$\{f_i\}_{i\in I}$ and $\{g_i\}_{i\in I}$ in a
Hilbert space are {\it equivalent} if there is an invertible
operator $T$ between their spans with $Tf_i = g_i$
for all $i\in I$.  
We now show that there is a simple
way to tell when two frame sequences are equivalent.

\begin{proposition}\label{FTP10}
Let $\{f_i\}_{i\in I}$, $\{g_i\}_{i\in I}$ be frames for a Hilbert
space $\H$ with analysis operators $T_1$ and $T_2$, respectively.
The following are equivalent:

(1)  The frames $\{f_i\}_{i\in I}$ and $\{g_i\}_{i\in I}$ are equivalent.

(2)  ker $T_1$ = ker $T_2$.
\end{proposition}

{\it Proof}:
$(1)\Rightarrow (2)$:  If $Lf_i = g_i$ is an isomorphism, then
$Lf_i = LT_1 e_i = g_i = T_2 e_i$ quickly implies our statement
about kernels.

$(2)\Rightarrow (1)$:  Since $T_i|_{{(ker\ T_i)}^{\perp}}$ is an
isomorphism for $i=1,2$, if the kernels are equal, then
$$
T_2 \left ( T_1 |_{(ker\ T_2)^{\perp}}\right )^{-1} f_i = g_i
$$
is an isomorphism.
\qed

 In the finite
dimensional case, if
$\{g_j\}_{j=1}^{n}$ is an orthonormal basis of
${\ell}_2^n$ consisting of eigenvectors for $S$ with respective eigenvalues
$\{{\lambda}_j\}_{j=1}^{n}$, then for every $1\le j\le n$,
$\sum_{i\in I}|\langle f_i ,g_j \rangle|^2 = {\lambda}_j$.  In particular,
$\sum_{i\in I}\|f_i\|^2 =$ trace S ($=n$ if $\{f_i\}_{i\in I}$
is a Parseval frame).  An important result is

\begin{theorem}\label{FT20}
If $\{f_i\}_{i\in I}$ is a frame for $\H$ with frame bounds
$A,B$ and $P$ is any orthogonal projection on $\H$, then
$\{Pf_i\}_{i\in I}$ is a frame for $P\H$ with frame bounds
$A,B$.
\end{theorem}

{\it Proof}:
For any $f\in P\H$,
$$
\sum_{i\in I}|\langle f,Pf_i \rangle |^2 =
\sum_{i\in I}|\langle Pf,f_i \rangle |^2 = \sum_{i\in I}|\langle
f,f_i \rangle |^2.
$$
\qed

A fundamental result in frame theory was proved independently
by Naimark and Han/Larson \cite{C,HL}.  For completeness we
include its simple proof.

\begin{theorem}\label{T3}
A family $\{f_i\}_{i\in I}$ is a Parseval frame for a Hilbert
space $\H$ if
and only if there is a containing Hilbert space $\H \subset {\ell}_2 (I)$
with an orthonormal basis $\{e_i\}_{i\in I}$ so that
the orthogonal projection $P$ of ${\ell}_2 (I)$ onto $\H$ satisfies
$P(e_i) = f_i$ for all $i\in I$.
\end{theorem}

{\it Proof}:
The ``only if'' part is Theorem \ref{FT20}.  For the ``if'' part,
if $\{f_i\}_{i\in I}$ is a Parseval frame, then the synthesis operator
 $T:{\ell}_2(I) \rightarrow \H$ is a partial isometry.  So $T^{*}$
is an isometry and we can associate $\H$ with $T^{*}\H$.  Now, for
all $i\in I$ and all $g=T^{*}f \in T^{*}\H$ we have
$$
\langle T^{*}f,Pe_i \rangle = \langle T^{*}f,e_i \rangle
= \langle f,Te_i \rangle = \langle f,f_i \rangle = \langle T^{*}f,
T^{*}f_i \rangle.
$$
It follows that $Pe_i = T^{*}f_i$ for all $i\in I$.
\qed

Now we can establish our Fundamental Principle II which basically
states that any bounded operator on a finite dimensional Hilbert
space is really just a multiple of a ``piece'' of a projection from a larger
space.

\begin{theorem}[Fundamental Principle II] \label{FPII}
Let $\HN$ be an $n$-dimensional Hilbert space with orthonormal
basis $\{g_i\}_{i=1}^{n}$.  If $T:\HN \rightarrow \HN$
is any bounded linear operator with $\|T\|=1$,
then there is a containing Hilbert
space $\HN \subset {\ell}_2^M$ (M=2n-1)
with an orthonormal basis $\{e_i\}_{i=1}^{M}$ so that the orthogonal
projection $P$ from ${\ell}_2^M$ onto $\HN$ satisfies:
$$
Pe_i = Tg_i, \ \ \mbox{for all $i=1,2,\ldots , n$}.
$$
\end{theorem}

{\it Proof}:
Let $S$ be the frame operator for the Bessel sequence
$\{f_i\}_{i=1}^{n} = \{Tg_i\}_{i=1}^{n}$
having eigenvectors $\{x_i\}_{i=1}^{n}$ with respective
eigenvalues $\{{\lambda}_i\}_{i=1}^{n}$ where $1={\lambda}_1 \ge
{\lambda}_2 \ge \cdots \ge {\lambda}_n$.  For $i=2,3,\ldots ,n$ let
$h_i = \sqrt{(1-{\lambda}_i )}x_i$.  Then, $\{f_i\}_{i=1}^n \cup
\{h_i\}_{i=2}^{n}$ is a Parseval frame for $\H$ since for every
$f\in \H$ we have
\begin{eqnarray*}
\sum_{i=1}^{n}|\langle f,f_i \rangle |^2 + \sum_{i=2}^{n}
|\langle f,h_i \rangle |^2 &=& \langle Sf,f\rangle +
\sum_{i=2}^{n}(1- {\lambda}_i )|\langle f,x_i \rangle|^2 \\
&=&\sum_{i=1}^{n}
{\lambda}_i |\langle f,x_i \rangle|^2 +
\sum_{i=2}^{n}(1- {\lambda}_i )|\langle f,x_i \rangle|^2 \\
&=& \sum_{i=1}^{n} |\langle f,x_i \rangle |^2 = \|f\|^2.
\end{eqnarray*}
Now, by Theorem \ref{T3}, there is a containing Hilbert space
${\ell}_2^{2n-1}$ with an orthonormal basis $\{e_i\}_{i=1}^{2n-1}$
so that the orthogonal projection $P$ satisfies: $Pe_i = Tg_i$
for $i=1,2,\ldots , n$ and $Pe_i = h_i$ for $i=n+1,\ldots , 2n-1$.
\qed

  For an introduction
to frame theory we refer the reader to Christensen \cite{C}.

Weaver \cite{W} established an important relationship between
frames and KS by showing that the following conjecture is
equivalent to KS.

\begin{conjecture}\label{C5}
There are universal constants $B\ge 4$ and $\epsilon > \sqrt{B}$
and an $r\in \N$ so that the following holds:  Whenever
$\{f_i\}_{i=1}^{M}$ is a unit norm $B$-tight frame for ${\ell}_2^n$,
there exists a partition $\{A_j\}_{j=1}^{r}$
of $\{1,2,\ldots ,M\}$ so that for all $j=1,2,\ldots ,r$ and
all $f\in {\ell}_2^n$ we have
\begin{equation}\label{E6}
\sum_{i\in A_j}|\langle f,f_i \rangle |^2 \le (B-\epsilon )\|f\|^2.
\end{equation}
\end{conjecture}

In his work on time-frequency analysis, Feichtinger \cite{CV}
noted that all of the Gabor frames he was using (see Section \ref{TFA})
had the property that they could be divided into a finite number
of subsets which were Riesz basic sequences.  This led to the
conjecture:

\begin{feichtinger conjecture}[FC]
Every bounded frame (or equivalently, every unit norm frame)
is a finite union of Riesz basic sequences.
\end{feichtinger conjecture}

There is a significant body of work on this conjecture
\cite{BCHL,BCHL2,CV,G}.  Yet, it remains open even for Gabor frames.
In \cite{CCLV} it was shown that FC is
equivalent to the weak BT, and hence is implied by KS
(See Section \ref{BST}).
In \cite{CT} it was shown
that FC
is equivalent to KS (See Theorem \ref{BSTCT}).  In fact,
we now know that KS is equivalent
to the {\it weak} Feichtinger Conjecture:  Every unit norm Bessel
sequence is a finite union of Riesz basic sequences (See Section
\ref{BST}).  In \cite{CKS2} it was shown that FC is equivalent to the
following conjecture.

\begin{conjecture}\label{Conj101}
Every bounded Bessel sequence is a finite union of frame sequences.
\end{conjecture}

Let us mention two more useful
equivalent formulations of KS due to Weaver
\cite{W}.

\begin{conjecture}[$\mathrm{KS_r}$]\label{W2}
There is a natural number $r$
and universal constants $B$ and
$\epsilon >0$ so that the following holds.  Let $\{f_i\}_{i=1}^{M}$
be elements of ${\ell}_2^n$ with $\|f_i\| \le 1$ for $i = 1,2,\ldots ,M$
and suppose for every $f\in {\ell}_2^n$,
\begin{equation}\label{EFT1}
\sum_{i=1}^{M}|\langle f,f_i \rangle |^2 \le B\|f\|^2.
\end{equation}
Then, there is a partition $\{A_j\}_{j=1}^{r}$ of $\{1,2,\ldots ,n\}$
so that for all $f\in {\ell}_2^n$ and all $j=1,2,\ldots ,r$,
$$
\sum_{i\in A_j}|\langle f,f_i \rangle |^2 \le (B-\epsilon )\|f\|^2.
$$
\end{conjecture}
Weaver \cite{W} also shows that Conjecture $\mathrm{KS_r}$ is equivalent to
PC if we assume equality in Equation \ref{EFT1} for all
$f\in {\ell}_2^n$.
Weaver further shows that Conjecture \ref{W2} is equivalent to KS even if we
strengthen its assumptions so as to require that the vectors $\{f_i\}_{i=1}^{M}$ 
are of equal norm and that equality holds in \ref{EFT1},
but at great cost to our $\epsilon >0$.

\begin{conjecture}[$\mathrm{KS_r}'$]\label{W3}
There exists universal constants $B\ge 4$ and $\epsilon > \sqrt{B}$
so that the following holds.  Let $\{f_i\}_{i=1}^{M}$
be elements of ${\ell}_2^n$ with $\|f_i\| \le 1$ for $i = 1,2,\ldots ,M$
and suppose for every $f\in {\ell}_2^n$,
\begin{equation}\label{EFT1}
\sum_{i=1}^{M}|\langle f,f_i \rangle |^2 = B\|f\|^2.
\end{equation}
Then, there is a partition $\{A_j\}_{j=1}^{r}$ of $\{1,2,\ldots ,n\}$
so that for all $f\in {\ell}_2^n$ and all $j=1,2,\ldots ,r$,
$$
\sum_{i\in A_j}|\langle f,f_i \rangle |^2 \le (B-\epsilon )\|f\|^2.
$$
\end{conjecture}

 We now strengthen the assumptions in
Fundamental Principle I, Conjecture \ref{C2}.

\begin{conjecture}\label{W4}
There exist universal constants $0<\delta, \sqrt{\delta}\le \epsilon <1$ and
$r\in \N$ so that for all $n$ and all orthogonal projections $P$
on ${\ell}_2^n$
with ${\delta}(P)\le \delta$ and $\|Pe_i\|= \|Pe_j\|$ for all
$i,j = 1,2,\ldots ,n$,
 there is a paving $\{A_j\}_{j=1}^{r}$
of $\{1,2,\ldots ,n\}$ so that $\|Q_{A_j}PQ_{A_j}\|\le
1-\epsilon$, for all $j=1,2,\ldots ,r$.
\end{conjecture}

Using Conjecture \ref{W3} we can see that KS is equivalent to
Conjecture \ref{W4}.

\begin{theorem}\label{W5}
KS is equivalent to Conjecture \ref{W4}.
\end{theorem}

{\it Proof}:
It is clear that Conjecture \ref{C2} (which is equivalent to KS)
implies Conjecture \ref{W4}.  So we assume that Conjecture \ref{W4} holds
and we will show that Conjecture \ref{W3} holds.
Let $\{f_i\}_{i=1}^{M}$
be elements of $\HN$ with $\|f_i\| = 1$ for $i = 1,2,\ldots ,M$
and suppose for every $f\in \HN$,
\begin{equation}\label{EFT4}
\sum_{i=1}^{M}|\langle f,f_i \rangle |^2 = B\|f\|^2,
\end{equation}
where $\frac{1}{B}\le \delta$.
It follows from Equation \ref{EFT4} that $\{\frac{1}{\sqrt{B}}f_i\}_{i=1}^{M}$
is an equal norm Parseval frame and so there is a larger Hilbert space
${\ell}_2^M$ and a projection $P:{\ell}_2^M \rightarrow \HN$ so that
$Pe_i = f_i$ for all $i=1,2,\ldots ,M$.  Now,
$\|Pe_i\|^2 = \langle Pe_i, e_i \rangle = \frac{1}{B}\le \delta$
for all $i=1,2,\ldots ,M$.  So by Conjecture \ref{W4}, there is a paving
$\{A_j\}_{j=1}^{r}$ of $\{1,2,\ldots ,M\}$ so that
$\|Q_{A_j}PQ_{A_j}\|\le
1-\epsilon$, for all $j=1,2,\ldots ,r$.  Now, for all $1\le j\le r$ and
all $f\in {\ell}_2^n$ we have:
\begin{eqnarray*}
\|Q_{A_j}Pf\|^2 &=& \sum_{i=1}^{M}|\langle Q_{A_j}Pf,e_i \rangle |^2
= \sum_{i=1}^{M}|\langle f, PQ_{A_j}e_i \rangle |^2\\
&=& \frac{1}{B}\sum_{i\in A_j}|\langle f,f_i \rangle |^2\\
&\le& \|Q_{A_j}P\|^2 \|f\|^2\\
&=& \|Q_{A_j}PQ_{A_j}\| \|f\|^2  \le (1-
{\epsilon})\|f\|^2.
\end{eqnarray*}
It follows that for all $f\in \HN$ we have
$$
\sum_{i\in A_j}|\langle f,f_i \rangle |^2 \le (B - {\epsilon}B)\|f\|^2.
$$
Since ${\epsilon}B> \sqrt{B}$, we have verified Conjecture \ref{W3}.
\qed

We give one final formulation of KS in Hilbert space frame theory.

\begin{theorem}
The following are equivalent:

(1)  The Paving Conjecture.

(2)  For every unit norm $B$-Bessel sequence $\{f_i\}_{i=1}^{M}$ in $\HN$
and every $\epsilon >0$, there exists $r=f(B, \epsilon)$
and a partition $\{A_j\}_{j=1}^{r}$ of $\{1,2,\ldots ,M\}$
so that for every $j=1,2,\ldots ,r$
and all scalars $\{a_i\}_{i\in A_j}$ we have
$$
\sum_{n\in A_j}|\langle f_n,\sum_{n\not= m\in A_j}a_m f_m \rangle |^2
\le \epsilon \|\sum_{m\in A_j}a_m f_m \|^2.
$$
\end{theorem}

{\it Proof}:
$(1)\Rightarrow (2)$:  Let $G$ be the Gram operator for $\{f_i\}_{i=1}^{M}$.
By PC, we can partition $\{1,2,\ldots ,M\}$ into $\{A_j\}_{j=1}^{r}$
so that for all $j=1,2,\ldots ,r$ we have
$$
\|P_{A_j}(G-D(G))Q_{A_j}\| \le \epsilon.
$$
Now, for any $j=1,2,\ldots ,r$ and any scalars $\{a_m\}_{m\in A_j}$ we
have
\begin{eqnarray*}
\sum_{n\in A_j}|\langle f_n,\sum_{n\not= m \in A_j}a_m f_m \rangle |^2
&=& \|Q_{A_j}(G-D(G))Q_{A_j}(\sum_{m\in A_j}a_m f_m)\|^2\\
&\le& \epsilon \sum_{n\in A_j}|a_n|^2\\
&\le& \frac{\epsilon}{1-\epsilon}\|\sum_{n\in A_j}a_nf_n\|^2,
\end{eqnarray*}
where the last inequality follows from the $R_{\epsilon}$-Conjecture
(actually, its proof using PC, see section \ref{HST}).

$(2)\Rightarrow (1)$:  Given (2), we have
\begin{eqnarray*}
\|\sum_{n\in A_j}a_n f_n\|^2 &=& \sum_{n\in A_j}|a_n|^2 +
\sum_{n\not= m\in A_j}a_n \overline{a_m}\langle f_n,f_m \rangle\\
&=& \sum_{n\in A_j}|a_n|^2 + \sum_{n\in A_j}a_n \langle f_n,
\sum_{n\not= m\in A_j}a_mf_m\rangle.
\end{eqnarray*}
Using (2) we now compute:
\begin{eqnarray*}
\left | \sum_{n\in A_j}a_n \langle f_n,\sum_{n\not= m\in A_j}a_mf_m
\rangle \right |^2 &\le& \left ( \sum_{n\in A_j}|a_n|^2\right )
 \sum_{n\in A_j}|\langle f_n,
\sum_{n\not= m\in A_j}a_m f_m \rangle |^2\\
&\le& \left ( \sum_{n\in A_j}|a_n|^2 \right ) \cdot
\epsilon \|\sum_{m\in A_j}a_m f_m \|^2\\
&\le& \left ( \sum_{n\in A_j}|a_n|^2 \right ) \cdot
\epsilon \cdot B \sum_{n\in A_j}|a_n|^2.
\end{eqnarray*}
This is enough to verify the $R_{\epsilon}$-Conjecture
(See Section \ref{HST}).
\qed

An  important open problem in frame theory is:
\begin{problem}
Classify the equal norm Parseval frames with {\it special properties}.
\end{problem}
The {\it special properties} here could be {\it translation invariance}
(Section \ref{TFA}), {\it reconstruction after erasures}
(Section \ref{Eng}), {\it frames which decompose into good
frame sequences} (see Section \ref{Eng}), etc.  The idea here
is to build up a ``bookshelf'' of equal norm frames with special
properties which can be used for applications such as we have for
wavelets.  As we will see, this problem shows up in many formulations
of KS.

\section{Kadison-Singer in Hilbert space theory}\label{HST}
\setcounter{equation}{0}

In this section we will see that KS is actually a fundamental
result concerning inner products.  Recall that a family of
vectors $\{f_i\}_{i\in I}$ is a {\bf Riesz basic sequence} in
a Hilbert space $\H$ if there are constants $A,B>0$ so that
for all scalars $\{a_i\}_{i\in I}$ we have:
\begin{displaymath}
A\sum_{i\in I}|a_i |^2 \le \|\sum_{i\in I}a_i f_i \|^2 \le
B\sum_{i\in I}|a_i|^2.
\end{displaymath}
We call $\sqrt{A},\sqrt{B}$ the {\bf lower and upper Riesz basis
bounds} for $\{f_i\}_{i\in I}$.  If $\epsilon >0$ and
$A = 1-\epsilon, B=1+\epsilon$ we call $\{f_i\}_{i\in I}$ an
$\epsilon$-{\bf Riesz basic sequence}.  If $\|f_i\|=1$ for all
$i\in I$ this is a {\bf unit norm} Riesz basic sequence.
A natural question is whether we can improve the Riesz basis
bounds for a unit norm Riesz basic sequence by partitioning the sequence
into subsets.  Formally:

\begin{conjecture}[$R_{\epsilon}$-Conjecture]
For every $\epsilon >0$, every unit norm Riesz basic sequence
is a finite union of $\epsilon$-Riesz basic sequences.
\end{conjecture}

The $R_{\epsilon}$-Conjecture was first stated in \cite{CV} where
it was shown that KS implies this conjecture.
 It was recently shown in \cite{CT} that KS is equivalent to the
$R_{\epsilon}$-Conjecture.  We include this argument here since
it demonstrates a fundamental principle we will employ throughout
this paper.

\begin{theorem}\label{THS1}  The following are equivalent:

(1)  The Paving Conjecture.

(2) If $T:{\ell}_2 \rightarrow {\ell}_2$ is a bounded linear
operator with $\|Te_i\|=1$ for all $i\in I$, then
for every $\epsilon >0$, $\{Te_i\}_{i\in I}$
is a finite union of
$\epsilon$-Riesz basic sequences.

(3)  The $R_{\epsilon}$-Conjecture.
\end{theorem}

{\bf Proof}:
$(1) \Rightarrow (2)$:
Fix $\epsilon >0$.  Given $T$ as in (2), let $S=T^{*}T$.
Since $S$ has ones on its diagonal, 
the (infinite form of the) Paving Conjecture gives
$r=r(\epsilon ,\|T\|)$ and a
partition $\{{A}_{j}\}_{j=1}^{r}$
of $I$ so that for
every
$j=1,2,\ldots , r$ we have
$$
\|Q_{A_j}(I-S)Q_{A_j}\|\le {\delta}\|I-S\|
$$
where $\delta = \epsilon/(\|S\| +1)$.  Now, for all
$f=\sum_{i\in I}a_ie_i$
we have
\begin{eqnarray*}
\| \sum_{i\in A_j}a_i Te_i\|^2 &=& \|TQ_{A_j}f\|^2 \\
&=&
\langle TQ_{A_j}f,TQ_{A_j}f\rangle \\
&=& \langle T^{*}TQ_{A_j}f,Q_{A_j}f\rangle \\
&=& \langle Q_{A_j}f,Q_{A_j}f \rangle -
\langle Q_{A_j}(I-S)Q_{A_j}f,Q_{A_j}f\rangle \\
&\ge& \|Q_{A_j}f\|^2 - {\delta}\|I-S\|\|Q_{A_j}f\|^2\\
&\ge& (1-\epsilon )\|Q_{A_j}f\|^2 = (1-\epsilon)\sum_{i\in A_j}
|a_i|^2.
\end{eqnarray*}
Similarly, $\|\sum_{i\in A_j}a_iTe_i\|^2 \le (1+\epsilon )\sum_{i\in A_j}
|a_i|^2.$

$(2)\Rightarrow (3)$:  This is obvious.

$(3)\Rightarrow (1)$:  Let $T\in B({\ell}_2)$
 with $Te_i = f_i$ and $\|f_i\|=1$ for all
$i\in I$.  We need to show that the Gram operator {\bf G}
of $\{f_i\}_{i\in I}$
is pavable.  Fix $0<\delta < 1$ and let $\epsilon >0$.  Let
$g_i = \sqrt{1- {\delta}^2}f_i \oplus {\delta}e_i \in {\ell}_2\oplus
{\ell}_2$.  Then ,$\|g_i\|=1$ for all $i\in I$
 and for all scalars $\{a_i\}_{i\in I}$
\begin{eqnarray*}
{\delta}\sum_{i\in I}|a_i|^2 &\le& \|\sum_{i\in I}a_i g_i\|^2
= (1-{\delta}^2)\|\sum_{i\in I}a_iTe_i\|^2 +
{\delta}^2\sum_{i\in I}|a_i|^2 \\
&\le& \left [ (1-{\delta}^2)\|T\|^2 + {\delta}^2 \right ]
\sum_{i\in I}|a_i|^2.
\end{eqnarray*}
So $\{g_i\}_{i\in I}$ is a unit norm Riesz basic sequence and
$\langle g_i,g_k\rangle= (1-{\delta}^2)\langle f_i,f_k\rangle$
for all $i\not= k \in I$.  By the
$R_{\epsilon}$-Conjecture, there is a partition $\{A_j\}_{j=1}^{r}$
so that for all $j=1,2,\ldots ,r$ and all $f=\sum_{i\in I} a_i e_i$,
\begin{eqnarray*}
(1-\epsilon )\sum_{i\in A_j}|a_i|^2 &\le& \|\sum_{i\in A_j}a_ig_i\|^2
= \langle \sum_{i\in A_j}a_ig_i ,\sum_{k\in A_j}a_k g_k \rangle \\
&=& \sum_{i\in A_j}|a_i|^2 \|g_i\|^2 + \sum_{i\not= k\in A_j}
a_i \overline{a_k}\langle g_i,g_k\rangle \\
&=& \sum_{i\in A_j}|a_i|^2 + (1-{\delta}^2)\sum_{i\not= k \in A_j}
a_i \overline{a_k}\langle f_i,f_k \rangle \\
&=& \sum_{i\in A_j}|a_i|^2 + (1-{\delta}^2)\langle Q_{A_j}(G-D(G))Q_{A_j}
f, f\rangle \\
&\le& (1+\epsilon)\sum_{i\in A_{j}}|a_i|^2.
\end{eqnarray*}
Subtracting $\sum_{i\in A_j}|a_i|^2$ through the inequality yields,
$$
-{\epsilon}\sum_{i\in A_j}|a_i|^2 \le
(1-{\delta}^2)\langle Q_{A_j}(G-D(G))Q_{A_j}f,f\rangle
\le \epsilon \sum_{i\in A_j}|a_i|^2.
$$
That is,
$$
 (1-{\delta}^2)|\langle Q_{A_j}(G-D(G))Q_{A_j}f,f\rangle |
\le \epsilon \|f\|^2.
$$
Since $Q_{A_j}(G-D(G))Q_{A_j}$ is a self-adjoint operator, we have
$(1-{\delta}^2)\|Q_{A_j}(G-D(G))Q_{A_j}\|\le \epsilon$.  That is,
$(1-{\delta}^2)G$ (and hence $G$) is pavable.
\qed

\begin{remark}\label{R2}
The proof of $(3)\Rightarrow (1)$ of Theorem \ref{THS1} illustrates
a standard method for turning conjectures about unit norm Riesz
basic sequences $\{g_i\}_{i\in I}$ into conjectures about
unit norm Bessel sequences $\{f_i\}_{i\in I}$.
  Namely, given $\{f_i\}_{i\in I}$ and $0<\delta <1$,
let $g_i = \sqrt{1-{\delta}^2}f_i\oplus {\delta}e_i\in
{\ell}_2(I)\oplus {\ell}_2(I)$.  Then, $\{g_i\}_{i\in I}$ is a unit
norm Riesz basic sequence and for $\delta$ small enough, $g_i$
is close enough to $f_i$ to pass inequalities from $\{g_i\}_{i\in I}$
to $\{f_i\}_{i\in I}$.
\end{remark}

It follows from Remark \ref{R1} that we can finite-dimensionalize
the result in Theorem \ref{THS1}.

\begin{conjecture}\label{CHS1}
For every $\epsilon >0$ and every $T\in B({\ell}_2^n)$
with $\|Te_i\|=1$ for $i=1,2,\ldots ,n$ there is
an $r=r(\epsilon ,\|T\|)$
and a partition $\{A_j\}_{j=1}^{r}$ of $\{1,2,\ldots ,n\}$
so that for all $j=1,2,\ldots ,r$ and all scalars
$\{a_i\}_{i\in A_j}$ we have
$$
(1-\epsilon )\sum_{i\in A_j}|a_i|^2 \le
\|\sum_{i\in A_j}a_i Te_i \|^2 \le (1+\epsilon)
\sum_{i\in A_j}|a_i|^2.
$$
\end{conjecture}

By Remark \ref{R2}, we can reformulate Conjecture \ref{CHS1} into
a statement about unit norm Riesz basic sequences.

One advantage of the $R_{\epsilon}$-Conjecture is that it can be shown
to students right at the beginning of a course in Hilbert spaces.
We note that this conjecture fails for equivalent norms
on a Hilbert space.  For example, if we renorm ${\ell}_{2}$ by letting
$|\{a_i\}| = \|{a_i}\|_{{\ell}_2} +\sup_i |a_i|$, then the
$R_{\epsilon}$-Conjecture fails for this equivalent norm. To see this,
let $f_i = (e_{2i}+e_{2i+1})/(\sqrt{2}+1)$ where
$\{e_i\}_{i\in \N}$ is the unit vector basis of ${\ell}_2$.
This is now a
unit norm Riesz basic sequence, but no infinite subset satisfies the
$R_{\epsilon}$-Conjecture.  To check this, let
$J\subset \N$ with $|J|=n$ and $a_i = 1/\sqrt{n}$ for
$i\in J $.  Then,
$$
|\sum_{i\in J}a_i f_i| =  \frac{1}{\sqrt{2}+1}\left ( \sqrt{2} +
\frac{1}{\sqrt{n}}\right ).
$$
Since the norm above is bounded away from one for $n\ge 2$,
we cannot satisfy the requirements of the $R_{\epsilon}$-Conjecture.
It follows that a positive solution to KS would imply a fundamental new
result concerning ``inner products'', not just norms.  Actually, the
$R_{\epsilon}$-Conjecture is way too strong for proving KS.  As we will
see, having either the upper inequality or the lower inequality hold is
a sufficient enough assumption to prove KS - and for each of these we
just need a universal constant to work instead of $1-\epsilon$ or
$1+\epsilon$.

Using Conjecture \ref{C5} we can show that the following conjecture
is equivalent to KS:

\begin{conjecture}\label{CC1}
There is a universal constant $1\le D$ so that for all $T\in B({\ell}_2^n)$
with $\|Te_i\|=1$ for all $i=1,2,\ldots ,n$, there is an
$r=r(\|T\|)$ and a partition $\{A_j\}_{j=1}^{r}$ of $\{1,2,\ldots ,n\}$
so that for all $j=1,2,\ldots ,r$ and all scalars $\{a_i\}_{i\in A_j}$
$$
\|\sum_{i\in A_j}a_iTe_i \|^2 \le D\sum_{i\in A_j}|a_i|^2.
$$
\end{conjecture}

\begin{theorem}
Conjecture \ref{CC1} is equivalent to KS.
\end{theorem}

{\it Proof}:
Since Conjecture \ref{CHS1} clearly implies Conjecture \ref{CC1},
we just need to show that Conjecture \ref{CC1} implies Conjecture
\ref{C5}.  So, choose $D$ as in Conjecture \ref{CC1} and choose
$B\ge 4$ and $\epsilon > \sqrt{B}$ so that $D\le B-\epsilon$.
Let $\{f_i\}_{i\in I}$ be a unit norm $B$ tight frame for
${\ell}_2^n$.  If $Te_i = f_i$ is the synthesis operator for
this frame, then $\|T\|^2 = \|T^{*}\|^2 = B$.  So by Conjecture
\ref{CC1}, there is an $r=r(\|B\|)$ and a partition $\{A_j\}_{j=1}^{r}$
of $\{1,2,\ldots ,n\}$ so that for all $j=1,2,\ldots ,r$ and
all scalars $\{a_i\}_{i\in A_j}$
$$
\|\sum_{i\in A_j}a_iTe_i \|^2 = \|\sum_{i\in A_j}a_i f_i \|^2
\le D\sum_{i\in A_j}|a_i|^2 \le (B- \epsilon )\sum_{i\in A_j}|a_i|^2.
$$
So $\|TQ_{A_j}\|^2 \le B- \epsilon$ and for all $f\in {\ell}_2^n$
$$
\sum_{i\in A_j}|\langle f,f_i\rangle |^2 = \|(Q_{A_j}T)^{*}f\|^2 \le
\|TQ_{A_j}\|^2\|f\|^2 \le (B-\epsilon )\|f\|^2.
$$
This verifies that Conjecture \ref{C5} holds and so KS holds.
\qed

Remark \ref{R2} and Conjecture \ref{CC1} show that we only
need any universal upper bound in the $R_{\epsilon}$-Conjecture to hold
to get KS.

\section{Kadison-Singer in Banach space theory}\label{BST}
\setcounter{equation}{0}

In this section we state a fundamental theorem
of Bourgain and Tzafriri called the {\it restricted
invertibility principle}.  This theorem led to the
{\it (strong and weak) Bourgain-Tzafriri Conjectures}.
We will see that these conjectures are equivalent to KS.

In 1987, Bourgain and Tzafriri \cite{BT} proved a fundamental
result in Banach space theory known as the
{\bf restricted invertibility principle}.

\begin{theorem}[Bourgain-Tzafriri]\label{TBST1}
There are universal constants $A,c>0$ so that whenever
$T:{\ell}_{2}^{n}\rightarrow {\ell}_{2}^{n}$ is a linear
operator for which $\|Te_{i}\|=1$, for $1\le i\le n$,
then there exists
a subset ${\sigma}\subset \{1,2,\ldots , n\}$ of cardinality
$|{\sigma}|\ge{cn}/{\|T\|^{2}}$ so that for all
$j=1,2,\ldots ,n$ and for all
choices of scalars $\{a_{j}\}_{j\in {\sigma}}$,
$$
\|\sum_{j\in {\sigma}}a_{j}Te_{j}\|^{2} \ge
A\sum_{j\in {\sigma}}|a_{j}|^{2}.
$$
\end{theorem}

Theorem \ref{TBST1} gave rise
to a problem in the area which has received a great deal of
attention \cite{BT1,CT,CV}.

\begin{bourgain-tzafriri}[BT]\label{CBST1}
There is a universal constant $A>0$ so that
for every $B>1$ there is a natural number $r=r(B)$
satisfying:
For any natural number $n$,
if $T:{\ell}_2^n \rightarrow {\ell}_2^n$ is a linear operator
with $\|T\|\le B$ and $\|Te_{i}\|=1$ for all $i=1,2,\ldots , n$,
 then there is a partition
$\{A_{j}\}_{j=1}^{r}$ of $\{1,2,\ldots , n\}$ so that
for all $j=1,2,\ldots ,r$ and
 all choices of scalars $\{a_{i}\}_{i\in A_{j}}$ we have:
$$
\|\sum_{i\in A_{j}}a_{i}Te_{i}\|^{2}\ge A \sum_{i\in A_{j}}|a_{i}|^{2}.
$$
\end{bourgain-tzafriri}

It had been ``folklore'' for years that KS and BT must be equivalent.
But no one was quite able to actually give a proof of this fact.
Recently, Casazza and Vershynin \cite{CV} gave a formal proof of
the equivalence of KS and BT.
Sometimes BT is called {\bf strong BT}
since there is a weakening of
it called {\bf weak BT}.  In weak BT we allow $A$ to depend upon the
norm of the operator $T$.  A significant amount
of effort has been invested in trying to show that strong and
weak BT are
equivalent \cite{BCHL,CCLV,CV}.  Recently, Casazza and Tremain
\cite{CT} proved this equivalence by showing that these results
are all equivalent to yet another conjecture.

\begin{conjecture}\label{C27}
There exists a constant $A>0$ and a natural number $r$ so that for
all natural numbers $n$, if $T:{\ell}_2^n \rightarrow {\ell}_2^n$
with $\|Te_i\|=1$ for all $i=1,2,\ldots ,n$ and $\|T\|\le 2$, there
is a partition $\{A_j\}_{j=1}^{r}$ of $\{1,2,\ldots ,n\}$ so that
for all $j=1,2,\ldots ,r$ and all scalars $\{a_i\}_{i\in A_j}$ we have
$$
\|\sum_{i\in A_j}a_i Te_i \|^2 \ge A \sum_{i\in A_j}|a_i|^2.
$$
\end{conjecture}

The proof of the following theorem from \cite{CT} demonstrates
how we will use our two Fundamental Principles.

\begin{theorem}\label{BSTCT}
The following are equivalent:
\begin{enumerate}
\item The Kadison-Singer Problem.
\vspace*{.04cm}
\item The (strong) BT.
\vspace*{.04cm}
\item The (weak) BT.
\vspace*{.04cm}
\item Conjecture \ref{C27}
\vspace*{.04cm}
\item The Feichtinger Conjecture.
\end{enumerate}
\end{theorem}

{\it Proof}:
$(1) \Rightarrow (2)$:
By the $R_{\epsilon}$-Conjecture, KS implies (strong) BT.

It is clear that $ (2) \Rightarrow (3) \Rightarrow (4)$.

$(4) \Rightarrow (1)$:
It suffices to show that Conjecture \ref{C27}
implies Conjecture \ref{C2}.  
Let $r,A$ satisfy Conjecture \ref{C27}. 
Fix $0<\delta \le {3}/{4}$ and let $P$ be an
orthogonal projection on ${\ell}_2^n$ with ${\delta}(P)\le \delta$
(notation from Section \ref{Intro}).  Now, $\langle Pe_i,e_i \rangle =
\|Pe_i\|^2 \le \delta$ implies $\|(I-P)e_i\|^2 \ge 1- \delta \ge
\frac{1}{4}$.  Define $T:{\ell}_2^n \rightarrow {\ell}_2^n$ by
$Te_i = {(I-P)e_i}/{\|(I-P)e_i\|}$.  For any scalars
$\{a_i\}_{i=1}^{n}$ we have
\begin{eqnarray*}
\|\sum_{i=1}^{n}a_i Te_i\|^2 &=& \|\sum_{i=1}^{n}
\frac{a_i}{\|(I-P)e_i\|}(I-P)e_i\|^2 \\
&\le& \sum_{i=1}^{n}\left | \frac{a_i}{\|(I-P)e_i\|}\right | ^2 \\
&\le& 4\sum_{i=1}^{n}|a_i|^2.
\end{eqnarray*}
So $\|Te_i\|=1$ and $\|T\|\le 2$.  By Conjecture \ref{C27}, there is
a partition $\{A_j\}_{j=1}^{r}$ of $\{1,2,\ldots ,n\}$ so that
for all $j=1,2,\ldots ,r$ and all scalars $\{a_i\}_{i\in A_j}$ we have
$$
\|\sum_{i\in A_J}a_i Te_i \|^2 \ge A \sum_{i\in A_j}|a_i|^2.
$$
Hence,
\begin{eqnarray*}
\|\sum_{i\in A_j}a_i(I-P)e_i\|^2 &=& \|\sum_{i\in A_j}a_i \|(I-P)e_i\|Te_i\|^2 \\
&\ge& A \sum_{i\in A_j}|a_i|^2 \|(I-P)e_i\|^2 \\
&\ge& \frac{A}{4} \sum_{i\in A_j}|a_i|^2.
\end{eqnarray*}
It follows that for all scalars $\{a_i\}_{i\in A_j}$,
\begin{eqnarray*}
\sum_{i\in A_j}|a_i|^2 &=& \|\sum_{i\in A_j}a_iPe_i\|^2 +
\|\sum_{i\in A_j}a_i (I-P)e_i\|^2 \\
&\ge& \|\sum_{i\in A_j}a_iPe_i\|^2 + \frac{A}{4}\sum_{i\in A_j}|a_i|^2.
\end{eqnarray*}
Now, for all $f=\sum_{i=1}^{n}a_ie_i$
$$
\|PQ_{A_j}f\|^2 = \| \sum_{i\in A_j}a_i Pe_i\|^2 \le (1-\frac{A}{4})
\sum_{i\in A_j}|a_i|^2.
$$
Thus,
$$
\|Q_{A_j}PQ_{A_j}\| = \|PQ_{A_j}\|^2 \le 1-\frac{A}{4}.
$$
So Conjecture \ref{C2} holds.

$(1)\Rightarrow (5)$:  Since every unit norm frame $\{f_i\}_{i\in I}$
has the property that the operator $T(e_i) = f_i$ is bounded where
$\{e_i\}_{i\in I}$ is an orthonormal basis for $\H$, it follows that
from Theorem \ref{THS1} that PC implies FC.

$(5) \Rightarrow (4)$:  This arguement comes from \cite{CCLV}.
We will prove the contrapositive.  So we assume that (4) fails.
Then for every $M\in \mathbb
N$
and for every $A>0$ there is an $n=n(M,A)\in \mathbb N$,
a finite dimensional Hilbert space $H$ and a Bessel
sequence $\{f_{i}\}_{i=1}^{n}$ in $H$
with Bessel constant $2$ and $\|f_{i}\|=1$,
for all $1\le i\le n$, and whenever we partition $\{1,2,\ldots ,n\}$ 
into
sets $\{I_{j}\}_{j=1}^{M}$, then there exists some $1\le \ell \le M$ 
and
a set of scalars $\{a_{i}\}_{i\in I_{\ell}}$ with
$$
\|\sum_{i\in I_{\ell}}a_{i}f_{i}\|^{2}\le A\sum_{i\in
I_{\ell}}|a_{i}|^{2}.
$$
Now, for each $k\in \mathbb N$, 
we can choose a finite dimensional Hilbert 
space
$H_{k}$ of dimension, say $m_{k}$,  
and letting $M=k$ and $A=1/k$ above we can choose
$n_{k}=n(k,1/k)$ and $\{f_{i}^{k}\}_{i=1}^{n_{k}}$ satisfying
the above conditions.  Let $H=(\sum \oplus H_{k})_{{\ell}_{2}}$
and consider $\{f_{i}^{k}\}_{i=1,k=1}^{n_{k}\ ,\infty}$ as elements of
$H$.  For each $k\in {\mathbb N}$, let $\{e_{i}^{k}\}_{i=1}^{m_{k}}$
be an orthonormal basis for $H_{k}$ and consider  
$\{e_{i}^{k}\}_{i=1,k=1}^{ m_{k}\ ,\infty}$ as elements of $H$.
Since $\{e_{i}^{k}\}_{i=1,k=1}^{ m_{k}\ ,\infty}$ is an orthonormal
bais for $H$, the family
$\{f_{i}^{k}\}_{i=1,k=1}^{n_{k}\ ,\infty}\cup 
\{e_{i}^{k}\}_{i=1,k=1}^{\ m_{k} \ ,\infty}$ is a family of norm
one vectors in $H$ with Bessel bound $3$ and lower frame bound
$\ge 1$ and hence forms a frame for $H$.   
Fix $M,A>0$ and assume
we can partition this frame into
$M$ sets of Riesz basic sequences each with lower Riesz basis bound
$A$.  In particular, we can partition 
$\{f_{i}^{k}\}_{i=1,k=1}^{n_{k}\ ,\infty}$ into $M$ sets of Riesz
basic sequences each with lower Riesz basis bound $A$. 
But, for all $k$ with $k\ge M$ and $1/k\le A$,
$\{f_{i}^{k}\}_{i=1}^{n_{k}}$ cannot be partitioned into $M$ sets each
with lower Riesz basis bound $\ge A$, and hence
$\{f_{i}^{k}\}_{i=1,k=1}^{n_{k}\ ,\infty}$ cannot be partitioned this
way.
This shows that (5) fails. 
\qed

Finally, let us note that Remark \ref{R2} and
BT imply that KS is equivalent to just
the lower inequality in the $R_{\epsilon}$-Conjecture and even without
the lower constant having to be close to one.

\section{Kadison-Singer in harmonic analysis}\label{HA}
\setcounter{equation}{0}

In this section, we present a detailed study of the Paving Conjecture for Toeplitz operators, 
reducing this problem to an old and fundamental problem in Harmonic Analysis.
Given $\phi \in L^\infty([0,1])$, the corresponding {\bf Toeplitz operator} is $T_{\phi}:L^2[0,1]\rightarrow L^2[0,1]$,
$T_{\phi}(f) = f\cdot \phi$.  In the 1980's, much effort
was put into showing that the class of Toeplitz operators
satisfies the Paving Conjecture
(see Berman, Halpern, Kaftal and Weiss \cite{BHKW,HKW,HKW1,HKW2})
 during which time the uniformly
pavable operators were classified and it was shown
that $T_{\phi}$ is pavable if $\phi$ is Riemann
integrable \cite{HKW}.  In Section \ref{HA1} we will reduce the
conjecture to a fundamental question in harmonic analysis
and find the weakest conditions which need to be established
to verify PC.  In Section \ref{HA2} we give harmonic analysis
classifications of the Uniform Paving Property for Toeplitz
operators and for the Uniform Feichtinger Conjecture.  As a
consequence, we will discover a surprising universal identity
for all functions $f\in L^2[0,1]$.  Throughout this section
we will use the following notation.
\vskip10pt
\noindent {\bf Notation}:  If $I\subset \Z$, we let $S(I)$ denote the $L^2([0,1])$-closure of the span of the exponential functions with frequencies taken from $I$:
\begin{displaymath}
S(I)=\mathrm{cl}(\mathrm{span}\{\mathrm{e}^{2\pi\mathrm{i}nt}\}_{n\in I}).
\end{displaymath}

\subsection{The Paving Conjecture for Toeplitz operators}
\label{HA1}

A deep and fundamental question in Harmonic Analysis is to
understand the distribution of the norm of a function
$f\in S(I)$.  It is
known (Proposition \ref{PHA2}) if that if $[a,b] \subset [0,1]$ and $\epsilon >0$,
then there
is a partition of $\Z$ into arithmetic progressions
$A_j = \{nr+j\}_{n\in \Z}$, $0\le j\le r-1$
so that for all $f\in S(A_j)$ we have
$$
(1-\epsilon )(b-a)\|f\|^2
\le \|f\cdot {\chi}_{[a,b]}\|^2 \le (1+\epsilon)(b-a)\|f\|^2.
$$
What this says is that the functions in $S(A_j)$ have their norms nearly
uniformly distributed across $[a,b]$ and $[0,1]\setminus [a,b]$.  The
central question is whether such a result is true for arbitrary
measurable subsets of $[0,1]$ (but it is known that the partitions
can no longer be arithmetic progressions \cite{BS,HKW,HKW2}).
If $E$ is a measurable
subset of $[0,1]$,
 let $P_E$ denote the orthogonal projection
of $L^2[0,1]$ onto $L^2(E)$, that is, $P_E(f) = f\cdot {\chi}_E$.
The fundamental question here is then

\begin{conjecture}\label{C100}
If $E\subset [0,1]$ is measurable and $\epsilon >0$ is given,
there is a partition $\{A_j\}_{j=1}^{r}$ of $\Z$
so that for all $j=1,2,\ldots ,r$ and all $f\in S(A_j)$
\begin{equation}\label{E500}
(1-\epsilon )|E|\|f\|^2
\le \|P_{E}(f)\|^2 \le (1+\epsilon)|E|\|f\|^2.
\end{equation}
\end{conjecture}

Despite the many deep results in the field of Harmonic Analysis,
almost nothing is known about the distribution of the norms of functions
coming from the span of a finite subset of the characters,
 except that this question has deep connections to Number Theory \cite{BS}
(Also, see Theorem \ref{BST1}).
Very little progress has ever been made on Conjecture \ref{C100}
except for a specialized result of Bourgain and Tzafriri \cite{BT1}.
Any advance on
this problem would have broad applications throughout the field.

To this day, the
Paving Conjecture for Toeplitz operators remains a deep mystery.
The next theorem (from \cite{CT}) helps explain why so little progress has
been made on KS for Toeplitz operators --- this problem is in fact equivalent to Conjecture \ref{C100}.
Because this result is fundamental for the rest of this
section, we will give the proof from \cite{CT}.
To prove the theorem we will first look at the
decomposition of Toeplitz operators of the form $P_E$.

\begin{proposition}\label{P9}
If $E\subset [0,1]$ and $A\subset \Z$ then for every $f\in L^2 [0,1]$ we have
$$
\|P_{E}Q_{A}f\|^2 = |E|\|Q_{A}f\|^2 + \langle Q_{A}(P_{E}-D(P_{E}))
Q_{A}f,f\rangle,
$$
where $Q_A$ is the orthogonal projection of $L^2 [0,1]$ onto $S(A)$.
\end{proposition}

{\bf Proof}:
For any $f=\sum_{n\in \Z}a_n \mathrm{e}^{2\pi\mathrm{i}nt}\in L^2 [0,1]$ we have
\begin{eqnarray*}
\|P_{E}Q_{A}f\|^2 &=& \langle P_{E}Q_{A}f ,P_{E}Q_{A}f \rangle
= \langle \sum_{n\in A}a_n P_{E}(\mathrm{e}^{2\pi\mathrm{i}nt}),\sum_{m\in A}
a_m P_{E}(\mathrm{e}^{2\pi\mathrm{i}mt})\rangle \\
&=&
\sum_{n\in A}|a_n|^2 \|{\chi}_{E}\cdot \mathrm{e}^{2\pi\mathrm{i}nt}\|^2 +
\sum_{n\not= m\in A}a_n \overline{a_m}\langle P_{E}\mathrm{e}^{2\pi\mathrm{i}nt},
\mathrm{e}^{2\pi\mathrm{i}mt}\rangle \\
&=& |E|\sum_{n\in A}|a_n|^2 + \langle (P_{E}-D(P_{E}))\sum_{n\in A}a_n
\mathrm{e}^{2\pi\mathrm{i}nt}, \sum_{n\in A}a_n \mathrm{e}^{2\pi\mathrm{i}nt}\rangle \\
&=& |E|\|Q_{A}f\|^2 + \langle Q_{A}(P_{E}-D(P_{E}))Q_{A}f,f\rangle.
\end{eqnarray*}
\qed

Now we are ready for the theorem from \cite{CT}.

\begin{theorem}\label{THA5}
The following are equivalent:

(1)  Conjecture \ref{C100}.

(2)  For every measurable $E\subset [0,1]$, the Toeplitz operator
$P_E$ satisfies PC.

(3)  All Toeplitz operators satisfy PC.
\end{theorem}

{\bf Proof}:
$(2)\Leftrightarrow (3)$:  This follows from the fact that the class
of pavable operators is closed and the class of Toeplitz operators
are contained in the closed linear span of the Toeplitz operators
of the form $P_E$.
That is, an arbitrary bounded measurable function on
$[0,1]$ may be essentially uniformly approximated by simple functions.

$(1)\Leftrightarrow (2)$:  By Proposition \ref{P9}, 
Conjecture \ref{C100} holds if and only if for all $\epsilon>0$,
there exists a partition
$\{A_j\}_{j=1}^{r}$ such that
\begin{eqnarray*}
(1-\epsilon )|E|\| Q_{A_j}f\|^2 &\le&
|E|\|Q_{A_j}f\|^2 + \langle Q_{A_j}(P_E - D(P_E)Q_{A_j}f,f\rangle\\
&\le& (1+\epsilon)|E|\|Q_{A_j}f\|^2
\end{eqnarray*}
for all $j=1,2,\ldots,r$ and all $f\in L^2[0,1]$.

Subtracting like terms through the inequality yields
that this inequality is equivalent to
\begin{equation}\label{EE}
|\langle Q_{A_j}(P_E - D(P_E)Q_{A_j}f,f\rangle | \le
\epsilon |E|\|Q_{A_j}f\|^2.
\end{equation}
Since $Q_{A_j}(P_E - D(P_E)Q_{A_j}$ is a self-adjoint,
Equation \ref{EE} is equivalent to
$\|Q_{A_j}(P_E - D(P_E)Q_{A_j}\|\le \epsilon |E|$, and so $P_E$ is pavable.
\qed

Next we state a useful result from \cite{MV}.

\begin{proposition}\label{HAP7}
Suppose that ${\lambda}_{1},{\lambda}_{2}\ldots , {\lambda}_{N}$ are
distinct real numbers, and suppose that $\delta >0$ is chosen so that
$|{\lambda}_{m}-{\lambda}_{n}|\ge \delta$ whenever $n\not= m$.
Then for any coefficients $\{a_n \}_{n=1}^{N}$, and any $T>0$ we have
$$
\int_{0}^{T}|\sum_{n=1}^{N}a_n \mathrm{e}^{2\pi\mathrm{i}{\lambda}_{n}t}|^2dt
= \left ( T + \frac{\theta}{\delta}\right ) \sum_{n=1}^{N}|a_n |^2,
$$
for some $\theta$ with $-1\le \theta \le 1$.
\end{proposition}

As an immediate consequence of \ref{HAP7} we have

\begin{proposition}\label{PHA2}
Conjecture \ref{C100} holds for intervals.  Moreover, we can use
a partition of $\Z$ made up of arithmetic progressions.
\end{proposition}

We next show that a significantly weaker conjecture than Conjecture
\ref{C100} is equivalent to PC for Toeplitz operators.  It is clear
that this is the weakest inequality we can have and still get PC.

\begin{conjecture}\label{C14}
There is a universal constant $0<K$ so that for any measurable
set $E\subset [0,1]$ there is a partition $\Aj$ of $\Z$
so that for every $f\in S(A_j )$ we have $\|P_{E}f\|^2 \le
K|E|\|f\|^2$.
\end{conjecture}

Now we will see that Conjecture \ref{C14} is equivalent to PC
for Toeplitz operators.

\begin{proposition}\label{PP1}
Conjecture \ref{C14} is equivalent to Conjecture \ref{C100}.
\end{proposition}

{\it Proof}:
Conjecture \ref{C100} clearly implies Conjecture \ref{C14}.
Assuming Conjecture \ref{C14} holds,  fix $\epsilon >0$ and
a measurable set $E\subset [0,1]$.  
Choose intervals
$\{[a_k ,b_k ]\}_{k=1}^{\infty}$ so that
\begin{displaymath}
E\subset\bigcup_{k=1}^{\infty}[a_k ,b_k ],
\end{displaymath}
and
\begin{displaymath}
\left|\sum_{k=1}^{\infty}(b_k -a_k ) - |E|\right|< \frac{\epsilon |E|}{3K}.
\end{displaymath}
Next, choose $N$ so that
\begin{displaymath}
\sum_{k=N+1}^{\infty}(b_k -a_k ) < \frac{\epsilon |E|}{3K},
\end{displaymath}
and let
\begin{displaymath}
F = \bigcup_{k=1}^{N}[a_k ,b_k ],\ \
E_{1} = E\cap F,\ \ E_2 = E\setminus E_1 .
\end{displaymath}
Note that
\begin{displaymath}
|F\setminus E_1|
\le \left|(F\setminus E_1)\cup\left (\,\bigcup_{k=N+1}^{\infty}[a_k ,b_k ] \setminus E_{2}\right )\right|
\leq\left|\bigcup_{k=1}^{\infty}[a_k ,b_k ] \setminus E\right|<
\frac{\epsilon |E|}{3}.
\end{displaymath}
By Conjecture \ref{C14}, we can partition $\Z$ into $\Aj$ so that
for each $j$ and all $f\in S(A_j )$ we have
\begin{displaymath}
\|P_{E_2}f\|^2 \le K|E_2 |\|f\|^2,
\end{displaymath}
and
\begin{displaymath}
\|P_{F\setminus E_1}f\|^2 \le K|F\setminus E_1 |\|f\|^2.
\end{displaymath}
Since $\displaystyle{E_{2}}\subset \bigcup_{k=N+1}^{\infty}[a_k ,b_k ]$, for every
$f\in S(A_j )$ we have
$$
\|P_{E_2}f\|^2 \le K\frac{\epsilon |E|}{3K}\|f\|^2 =
\frac{\epsilon |E|}{3}\|f\|^2.
$$
Fix $1\le j\le r$.  By Proposition \ref{PHA2}, there is a partition
$\{B_k \}_{k=1}^{M_j}$ of $A_j$ so that for every $f\in S(B_k )$ we have

$$
(|F|-\frac{\epsilon |E|}{3})\|f\|^2 \le \|P_{F}f\|^2 \le
(|F|+\frac{\epsilon |E|}{3})\|f\|^2.
$$
Now, for every $f\in S(B_k )$ we have
\begin{eqnarray*}
\|P_{E}f\|^2 &\le& \|P_{E_1}f\|^2 +\|P_{E_2}f\|^2\\
&\le& \|P_{F}f\|^2 +
\frac{\epsilon |E|}{3}\|f\|^2\\
&\le& (|F|+\frac{\epsilon |E|}{3})\|f\|^2 + \frac{\epsilon |E|}{3}\|f\|^2\\
&\le& \sum_{k=1}^{N}(b_k -a_k ) \|f\|^2 + \frac{2\epsilon|E|}{3}\|f\|^2\\
&\le& (|E|+\frac{\epsilon |E|}{3K})\|f\|^2 + \frac{2\epsilon |E|}{3}\|f\|^2\\
&=& (1+\epsilon )|E|\|f\|^2,
\end{eqnarray*}
where, without loss of generality, we have assumed $K>1$.
For the other direction, we note that $P_{F}f= P_{E_1}f+P_{F\setminus E_1}f$
and $P_{E_1}f \perp P_{F\setminus E_1}f$ and so
\begin{eqnarray*}
\|P_{E}f\|^2 \ge \| P_{E_1}f\|^2 &=& \|P_F f\|^2 -
\|P_{F\setminus E_1}f\|^2\\
&\ge& (|F| - \frac{\epsilon |E|}{3})\|f\|^2 - \frac{\epsilon |E|}{3}\|f\|^2\\
&\ge& (|E| - \frac{\epsilon |E|}{3}-\frac{2\epsilon |E|}{3})\|f\|^2 \\
&=& (1- \epsilon )|E|\|f\|^2.
\end{eqnarray*}
\qed

We next present several equivalent formulations of a slightly weaker conjecture.

\begin{conjecture}\label{C11}
Suppose $E\subset [0,1]$ with $0<|E|$.  There is a partition
$\Aj$ of $\Z$ so that for all $j=1,2,\ldots , r$,
$P_{E}$ is an isomorphism of $S(A_j )$ onto its range.
\end{conjecture}

\begin{definition}
We say the Toeplitz operator $T_{\phi}$ satisfies the Feichtinger
Conjecture if there is a partition $\{A_j\}_{j=1}^{r}$ of $\Z$
so that $\{T_{\phi}\mathrm{e}^{2\pi\mathrm{i}nt}\}_{n\in A_j}$ is a Riesz basic
sequence for every $j=1,2,\ldots ,r$.
\end{definition}

Conjecture \ref{C11} is equivalent to all Toeplitz operators
satisfying the Feichtinger Conjecture.  That is,
for any Toeplitz operator $T_{\p}$, $\{T_{\p}\mathrm{e}^{2\pi\mathrm{i}nt}\}_{n\in \Z}$
is a finite union of Riesz basic sequences.

\begin{theorem}\label{HAT12}
The following are equivalent:

(1)  The Feichtinger Conjecture for Toeplitz operators.

(2)  The Feichtinger Conjecture for $P_E$ for every measurable
set $E\subset [0,1]$ with $0<|E|$.

(3)  Conjecture \ref{C11}.
\end{theorem}

{\it Proof}:
$(1)\Rightarrow (2)$:  This is obvious.

$(2)\Rightarrow (3)$:  If $E$ is a measurable subset of $[0,1]$
with $0<|E|$ and
we assume (2), then there is a partition $\{A_j\}_{j=1}^{r}$ of
$\Z$ and a constant $0<A$ so that $\{P_E \mathrm{e}^{2\pi\mathrm{i}nt}\}_{n\in A_j}$
is a Riesz basic sequence for all $j=1,2,\ldots ,r$ with lower Riesz
basis bound $A$.  Hence, for every $j=1,2,\ldots,r$ and every
$f=\sum_{n\in A_j}a_n \mathrm{e}^{2\pi\mathrm{i}nt}$ we have
$$
\|P_E f\|^2 = \|\sum_{n\in A_j}a_n P_E \mathrm{e}^{2\pi\mathrm{i}nt}\|^2
\ge A^2 \sum_{n\in A_j}|a_n|^2 = A^2\|f\|^2.
$$

$(3)\Rightarrow (1)$:  Let $T_{\phi}$ be a non-zero Toeplitz operator
on $L^2[0,1]$.  Choose $\epsilon >0$ and a measurable set
$E\subset [0,1]$ with $|E|>0$ so that $|{\phi}(t)|\ge \epsilon$ for
all $t\in E$.  By our assumption (3), there is a partition
$\{A_j\}_{j=1}^{r}$ of $\Z$ so that $P_E$ is an isomorphism of
$S(A_j)$ onto its range with, say, lower isomorphism
bound $A$.  That is, for all $j=1,2,\ldots ,r$ and
all $\{a_n\}_{n\in A_j}$ we have
\begin{eqnarray*}
\|\sum_{n\in A_j}a_nT_{\phi} \mathrm{e}^{2\pi\mathrm{i}nt}\|^2 &=&
\int_{0}^{1}|\sum_{n\in A_j}a_n\mathrm{e}^{2\pi\mathrm{i}nt}|^2|{\phi}(t)|^2\ dt\\
&\ge& \int_{0}^{1}|\sum_{n\in A_j}a_n\mathrm{e}^{2\pi\mathrm{i}nt}|^2\epsilon^2|{\chi}_E|^2\ dt\\
&\ge&
{\epsilon}^2\|P_E \sum_{n\in A_j}a_n \mathrm{e}^{2\pi\mathrm{i}nt}\|^2\\
&\ge&  {\epsilon}^2 A^2\sum|a_n|^2.
\end{eqnarray*}
\qed

At this time we do not know if the Feichtinger Conjecture for
Toeplitz operators is equivalent to PC for Toeplitz operators.
The main problem here is that we do not have an equivalent form
of Fundamental Principle I (Conjecture \ref{C2}) for Toeplitz operators.

There is some evidence for believing that Conjecture \ref{C11}
might be true since a weaker version of it holds as we see in
the next result.

\begin{proposition}\label{Prop6}
Assume $E\subset [0,1]$ with $0<|E|$.  Then, there is a partition
$\Aj$ of $\Z$ so that for every $j=1,2,\ldots ,r$ if $f\in S(A_j )$
and $f|_{E} = 0$ then $f=0$.
\end{proposition}

{\it Proof}:
We actually prove a stronger result,
namely that there is a single partition
that works for all measurable sets $E$.
In particular, let $A_1 = \N \cup \{0\}$ and $A_2 = \Z \setminus A_1$.
Now, if $f\in S(A_1 )$ then $\log|f| \in L^{1}[0,1]$ (See Duran \cite{D}, Page 16)
and so $f\not= 0$ on any set of positive
measure.  A similar argument (applied to the complex conjugate of $f$) applies in the case where $f\in S(A_2)$.
\qed

We end this section with one more equivalence of PC for Toeplitz
operators.

\begin{proposition}\label{PHA10}
For a Toeplitz operator $T_g$ the following are equivalent:

(1)  $T_g$ satisfies PC.

(2)  For every $\epsilon >0$ there is a partition $\Aj$ of
$\Z$ so that for every $j=1,2,\ldots ,r$ and for every
$f\in S(A_j)$ we have:
$$
\|f\|(\|g\|-\epsilon )\le \|f\cdot g \|\le \|f\|(\|g\|+\epsilon).
$$
\end{proposition}

{\it Proof}:
$(1)\Rightarrow (2)$:  If $T_g$ has the Paving Property,
then Conjecture \ref{C14} (and hence Conjecture \ref{C100}) holds.
Fix $\epsilon >0$.  Choose a simple function $h=\sum_{k=1}^{M}a_k
{\chi}_{E_k}$ such that $|g-h|<\epsilon$ almost everywhere,
where, without loss of generality, the sets $E_k$ are mutually disjoint.
 By Conjecture \ref{C100}, there is a partition $\Aj$ of $\Z$ so that for all
$j=1,2,\ldots ,r$ and for all $f\in S(A_j)$
we have
$$
\|f\|^2 \left ( |E_k|-\frac{\delta}{\sum_{k=1}^{M}|a_k|^2} \right )
 \le
\|f\cdot {\chi}_{E_k}\|^2 \le \|f\|^2 \left ( |E_k|+ \frac{\delta}
{\sum_{k=1}^{M}|a_k|^2}\right )
$$
for all $k=1,\dots,M$.
Now,
\begin{eqnarray*}
\|f\cdot h\|^2 &=& \| \sum_{k=1}^{M}a_k f\cdot {\chi}_{E_k}\|^2\\
&=& \sum_{k=1}^{M}|a_k|^2 \|f\cdot {\chi}_{E_k}\|^2\\
&\le& \sum_{k=1}^{M}|a_k|^2 \|f\|^2 \left ( |E_k| +
\frac{\delta}{\sum_{k=1}^{M}|a_k|^2} \right )\\
&\le& \|f\|^2 \left ( \sum_{k=1}^{M}|a_k |^2 |E_k | + \delta\right )\\
&=& \|f\|^2 (\|h\|^2 + \delta ).
\end{eqnarray*}
Similarly,
$$
\|f\cdot h\|^2 \ge \|f\|^2(\|h\|^2 - \delta ).
$$
Hence, for $\delta >0$ small enough we have
\begin{eqnarray*}
\|f\cdot g\| &\le& \|f\cdot h\| + \|f \cdot (g-h)\|\\
&\le& \|f\|\sqrt{\|h\|^2 + \delta} + \|f\| \delta \\
&\le& \|f\|\left ( \sqrt{(\|g\|+\delta )^2+\delta} + \delta \right )\\
&\le& \|f\|\sqrt{\|g\|^2 + \epsilon }.
\end{eqnarray*}
Thus,
$$
\|f\cdot g \|^2 \le \|f\|^2 (\|g\|^2 + \epsilon ).
$$
Similarly,
$$
\|f\cdot g\|^2 \ge \|f\|^2 (\|g\|^2 - \epsilon ).
$$

$(2)\Rightarrow (1)$:  This is immediate from Conjecture \ref{C100},
 and Theorem \ref{THA5}.
\qed

\subsection{The Uniform Paving Property}\label{HA2}

In this section we will classify the Toeplitz operators which
have the uniform paving property and the uniform Feichtinger property.

\begin{definition}
A Toeplitz operator $T_g$ has the {\it uniform paving property} if
for every $\epsilon >0$, there is a $K\in \mathbb{N}$ so that
if $A_k = \{nK+k\}_{n\in \Z}$ for $0\le k\le K-1$ then
$$
\|P_{A_k}(T_g-D(T_g))P_{A_k}\| < \epsilon.
$$
\end{definition}

\begin{definition}
A Toeplitz operator $T_g$ has the {\it uniform
Feichtinger property} if there is a $K\in \mathbb{N}$ so that
if $A_k = \{nK+k\}_{n\in \Z}$ for $0\le k\le K-1,$ then
$\{T_g \mathrm{e}^{2\pi\mathrm{i}nt}\}_{n\in A_k}$ is a Riesz basic sequence for all
$k=0,1,\ldots , K-1$.
\end{definition}

Halpern, Kaftal and Weiss \cite{HKW} made a detailed study of
the uniform paving property and in particular they showed that
$T_{\p}$ is uniformly pavable if $\p$ is a Riemann integrable
function.
As we saw in Section \ref{HA1}, the uniform paving property is really a
fundamental question in harmonic analysis.  In this section we will give
classifications of the Toeplitz operators having the uniform paving
property and those having the uniform Feichtinger property.  Our
approach will be of a harmonic analysis flavor and will lead to a
 new identity which holds for all $f\in L^2[0,1]$.

\begin{notation}
For all $g\in L^2[0,1]$, $K\in \N$
and any $0\le k\le K-1$ we let
$$
g_k^K(t) = \sum_{n\in \Z}\langle g,\mathrm{e}^{2\pi\mathrm{i}(nK+k)t}\rangle
\mathrm{e}^{2\pi\mathrm{i}(nK+k)t}.
$$
Also,
$$
g_{K}(t) = \frac{1}{K}\sum_{k=0}^{K-1}|g(t-\frac{k}{K})|^2.
$$
\end{notation}

The main theorems of this section are:

\begin{theorem}\label{HTT1}
Let $g\in L^{\infty}[0,1]$ and $T_{g}$ the Toeplitz operator
of multiplication by $g$.  The following are equivalent:

(1)  $T_{g}$ has uniform PC.

(2)  There is an increasing sequence of natural numbers $\{K_n\}$
so that
$$
\lim_{n\rightarrow \infty}\frac{1}{K_n}\sum_{k=0}^{K_n-1}
|g(t-\frac{k}{K_n})|^2 = \|g\|^2 \ \ a.e.
$$
uniformly over $t$.  That is, for every $\epsilon >0$ 
 there is
a $K\in \N$ so that
$$
|\frac{1}{K}\sum_{k=0}^{K-1}|g(t-\frac{k}{K})|^2 - \|g\|^2|< \epsilon
\ \ a.e.
$$
\end{theorem}

\begin{theorem}\label{HTT2}
Let $g\in L^{\infty}([0,1])$ and $T_{g}$ the Toeplitz operator
of multiplication by $g$.  The following are equivalent:

(1) $T_g$ has the uniform Feichtinger property.

(2)  There is a natural number $K\in \N$ and an $\epsilon >0$ 
so that 
$$
\frac{1}{K}\sum_{k=0}^{K-1}
|g(t-\frac{k}{K})|^2 \ge \epsilon \ \ a.e.
$$

(3) There is an $\epsilon >0$  and $K$ measurable sets
$$
E_k \subset [\frac{k}{K},\frac{k+1}{K}], \ \ \mbox{for all
$0\le k\le K-1$},
$$
satisfying:

\ \ \ \ \ \ (a)  The sets $\{E_{k}-\frac{k}{K}\}_{k=0}^{K-1}$ are disjoint
in $[0,\frac{1}{K}]$.

\ \ \ \ \ \ (b)  $\displaystyle\bigcup_{k=0}^{K-1}(E_k - k) = [0,\frac{1}{K}]$.

\ \ \ \ \ \ (c)  $|g(t)|\ge \epsilon$ on $\displaystyle\bigcup_{k=0}^{K-1}E_k$.
\end{theorem}

With a little effort we can recover the Halpern, Kaftal
and Weiss result \cite{HKW}.

\begin{corollary}\label{CoTO1}
For all Riemann integrable functions $\p$, the Toeplitz operator
$T_{\p}$ satisfies the uniform paving property.
If $|g|\ge \epsilon >0$ on an interval then $g$ has the uniform
Feichtinger property.  Hence, if $g$ is continuous at one point
and is not zero at that point, then $g$ has the uniform Feichtinger
property.
\end{corollary}
Examples of Toeplitz operators failing uniform pavability
were given in \cite{BS,HKW}.

\begin{example}\label{E1}
An example of a Toeplitz operator which fails the uniform Feichtinger
property.
\end{example}

{\it Proof}:
Choose $0< a_n$ so that $\sum na_n < 1$.  For each n, choose
$$
F_n \subset [0,\frac{1}{n}],\ \ \mbox{with}\ \ |F_n |= a_n.
$$
Let
$$
E_n = \bigcup_{k=0}^{n-1}(F_n+\frac{k}{n}), \ \ \mbox{and}\ \ E=\bigcup_{n=0}^{\infty}
E_n .
$$
Now,
$$
|E|\le \sum_{n}|E_n| = \sum na_n < 1.
$$
It is easily seen that $E^c$ contains no intervals.  If $g= {\chi}_{E^c}$,
then for all $K\in \Z$ we have
$$
\sum_{k=0}^{K-1}|g(t-\frac{k}{K})|^2 = 0,\ \ \mbox{on}\ \ E_K.
$$
So uniform Feichtinger fails.
\qed

\begin{example}
There is an open set $F\subset [0,1]$ so that for $g={\chi}_{F}$
the Toeplitz operator $T_g$ fails the uniform paving property.
\end{example}

{\it Proof}:
Let $g={\chi}_{F}$ where $F$ is the set given in Example \ref{E1}.
Let $F_n$ be the sets given in that example also.
Then $F$ is an open set with $|F|<1$ and for all $K\in \N$ we have
$$
\frac{1}{K}\sum_{k=0}^{K-1}|g(t-\frac{k}{K})|^2 = 1,\ \ \mbox{for all
$t\in F_n$}.
$$
Hence, $T_g$ fails uniform paving by Corollary \ref{CC1}.
\qed

\begin{example}
There is a $g\in L^2[0,1]$ so that the Toeplitz operator $T_{g}$
has the uniform Feichtinger property but fails the uniform
paving property.
\end{example}

{\it Proof}:
Choose a measurable set $F\subset [0,1/2]$ with ${\chi}_{F}$
failing uniform paving and let $E= F\cup [1/2,1]$.  By Corollary
\ref{CoTO1}
$P_E$ has the
uniform Feichtinger property but still fails uniform paving.
\qed

In order to prove Theorems \ref{HTT1} and \ref{HTT2}, we will need to do some
preliminary work.  Parts $(1), (2)$ of this proposition were
done originally by Halpern, Kaftal and Weiss \cite{HKW}, Lemma 3.4.

\begin{proposition}\label{LL1}
For any $g\in L^2([0,1])$ and any positive integer $K$,
\begin{enumerate}
\item $\displaystyle g_{k}^{K}(t) = \frac{1}{K}\sum_{j=0}^{K-1}g(t-\frac{j}{K})\mathrm{e}^{\frac{2\pi\mathrm{i}jk}{K}}$ for all $0\le k\le K-1$,
\item  For all $k, \ell \in \Z$ we
have:
$$
g_{k}^K(t-\frac{\ell}{K}) =\mathrm{e}^{-\frac{2{\pi}\mathrm{i}k\ell}{K}}g_{k}^K(t).
$$
\item If now $f=f_k^K$, then
$$
(f\cdot g )_{\ell}^K(t) = f(t)g_{\ell -k}^K(t).
$$
\end{enumerate}
\end{proposition}

{\it Proof}:
$(1)$:  For any $k=0,\dots,K-1$, consider $h\in L^2([0,1])$, 
\begin{displaymath}
h(t)=\frac{1}{K}\sum_{j=0}^{K-1}g(t-\frac{j}{K}).
\end{displaymath}
We want to show $h=g_k^K$, namely, 
\begin{displaymath}
h(t) = \sum_{n\in \Z}\langle g,\mathrm{e}^{2\pi\mathrm{i}(nK+k)t}\rangle
\mathrm{e}^{2\pi\mathrm{i}(nK+k)t}.
\end{displaymath}
To do this, note
\begin{eqnarray*}
\widehat{h}(n) 
&=& \int_{0}^{1}h(t)e^{-2{\pi}int}dt\\
&=& \frac{1}{K}\sum_{j=0}^{K-1}\mathrm{e}^{\frac{2{\pi}\mathrm{i}jk}{K}}\int_{0}^{1}
g(t-\frac{j}{K})e^{-2{\pi}int}dt,\ \ \mbox{let $s=t-\frac{j}{K}$}\\
&=& \frac{1}{K}\sum_{j=0}^{K-1}\mathrm{e}^{\frac{2{\pi}\mathrm{i}jk}{K}}\int_{0}^{1}
g(s)e^{-2{\pi}in(s+\frac{j}{K})}ds\\
&=& \frac{1}{K}\sum_{j=0}^{K-1}\mathrm{e}^{\frac{2{\pi}\mathrm{i}jk}{K}}\mathrm{e}^{-\frac{2{\pi}\mathrm{i}nj}{K}}
\int_{0}^{1}g(s)e^{-2{\pi}ins}ds\\
&=& \frac{1}{K}\sum_{j=0}^{K-1}\left [ \mathrm{e}^{-\frac{2{\pi}\mathrm{i}(n-k)}{K}}
\right ] ^{j} \hat{g}(n)\\
&=& \left\{\begin{array}{ll}\hat{g}(n)& \mathrm{e}^{-\frac{2{\pi}\mathrm{i}(n-k)}{K}} = 1\\0&\mathrm{e}^{-\frac{2{\pi}\mathrm{i}(n-k)}{K}}\neq1\end{array}\right.\\
&=& \left\{\begin{array}{ll}\hat{g}(n)& n\equiv{k}\mod{K}\\0&\mbox{else}\end{array}\right..
\end{eqnarray*}
Thus, $h$ and $g_k^K$ have the same Fourier coefficients, that is, $h=g_k^K$.

$(2)$:  We have:
\begin{eqnarray*}
g_{k}(t-\frac{\ell}{K}) &=& \frac{1}{K}\sum_{j=0}^{K-1}
g(t-\frac{\ell}{K}-\frac{j}{K})\mathrm{e}^{\frac{2{\pi}\mathrm{i}jk}{K}}\\
&=& \frac{1}{K}\sum_{j=0}^{K-1}g(t-\frac{\ell +j}{K})\mathrm{e}^{\frac{2{\pi}\mathrm{i}jk}{K}}\\
&=& \frac{1}{K}\sum_{j=0}^{K-1}g(t-\frac{\ell +j }{K})
\mathrm{e}^{\frac{2{\pi}\mathrm{i}(\ell + j)}{K}}\mathrm{e}^{-\frac{2{\pi}\mathrm{i}\ell k}{K}}\\
&=& \mathrm{e}^{-\frac{2{\pi}\mathrm{i}\ell k}{K}}g_{k}(t).
\end{eqnarray*}

$(3)$:  We have:
\begin{eqnarray*}
(f\cdot g )_{\ell}(t) &=& \frac{1}{K}\sum_{k=0}^{K}
f(t-\frac{j}{K})g(t- \frac{j}{K})
\mathrm{e}^{\frac{2{\pi}\mathrm{i}j{\ell}}{K}}\\
&=& \frac{1}{K}\sum_{j=0}^{K-1}f_{k}(t-\frac{j}{K})g(t-\frac{j}{K})
\mathrm{e}^{\frac{2{\pi}\mathrm{i}j\ell}{K}}\\
&=& \frac{1}{K}\sum_{j=0}^{K-1}\mathrm{e}^{-\frac{2{\pi}\mathrm{i}jk}{K}}f_{k}(t)
g(t-\frac{j}{K})\mathrm{e}^{\frac{2{\pi}\mathrm{i}j\ell}{K}}\\
&=& f_{k}(t)\frac{1}{K}\sum_{j=0}^{K-1}g(t-\frac{j}{K})
\mathrm{e}^{\frac{2{\pi}\mathrm{i}j(\ell -k)}{K}} = f(t)g_{\ell - k}(t).
\end{eqnarray*}
\qed

Now we can establish an important relationship between $g_k^K$ and $g$.

\begin{theorem}\label{TT3}
For all $g\in L^2 [0,1]$ we have
$$
\sum_{k=0}^{K-1}|g_{k}^{K}(t)|^2 =
\frac{1}{K}\sum_{k=0}^{K-1}
|g(t- \frac{k}{K})|^2.
$$
\end{theorem}

{\it Proof}:
With
$$
g_{k}(t) = \frac{1}{K}\sum_{j=0}^{K-1}g(t-\frac{j}{K})
\mathrm{e}^{\frac{2{\pi}\mathrm{i}jk}{K}},
$$
and for $t\in [0,1]$, let $h_{t}\in {\ell}({\mathbb Z}_{K})$ be
given by: $h_{t}(j) = g(t-\frac{j}{K})$.  Then
$$
g_{k}(t) = \frac{1}{K}\sum_{j=0}^{K-1}h_{t}(j)\mathrm{e}^{-\frac{2{\pi}\mathrm{i}jk}{K}}
= \frac{1}{\sqrt{K}}\widehat{h_{t}}(-k).
$$
Now,
\begin{eqnarray*}
\sum_{k=0}^{K-1}|g_{k}^{K}(t)|^2 &=&
\sum_{k=0}^{K-1}|g_{k}(t)|^2 = \sum_{k=0}^{K-1}|
\frac{1}{\sqrt{K}}\widehat{h_{t}}(-k)|^2\\
&=& \frac{1}{K}\sum_{k=0}^{K-1}|\widehat{h_{t}}(-k)|^2\\
&=& \frac{1}{K}\|\widehat{h_{t}}\|^2 = \frac{1}{K}\|h_{t}\|^2\\
&=& \frac{1}{K}\sum_{j=0}^{K-1}|h_{t}(j)|^2 = \frac{1}{K}
\sum_{j=0}^{K-1}|g(t-\frac{j}{K})|^2.
\end{eqnarray*}
Combined with our lemmas, this proves Theorem \ref{TT3}.
\qed

The next theorem gives an identity which holds for all $f\in L^2[0,1]$.
It says that {\it pointwise}, any Fourier series can be divided into
its subseries of arithmetic progressions so that the square sums of
the functions given by the subseries spreads the norm nearly equally
over the interval $[0,1]$ (Compare this to the discussion at the
beginning of Section \ref{HA1}).

\begin{theorem} \label{TT4}
For any $g\in L^\infty([0,1])$ there is an increasing sequence of
natural numbers $\{K_{n}\}_{n=1}^{\infty}$ so that
$$
\lim_{n\rightarrow \infty}\left ( \sum_{k=0}^{K_n-1}|g_{k}^{K_n}(t)|^2
\right )^{1/2} = \|g\|{\chi}_{[0,1]}\ \ a.e.
$$
\end{theorem}

{\it Proof}:
By Theorem \ref{TT3}, we can work with the functions $g_{K}$.
Also, we observe that it suffices to prove that this sum
converges in measure to $\|g\|^2{\chi}_{[0,1]}$.
We will do the proof in steps.

\noindent {\bf Step 1}:  The result holds for $g= {\chi}_{[S,T]}$.
And in fact, the convergence is uniform for these functions.

{\it Proof of Step 1}:  It suffices to assume $S=0$ and
$T< 1$.  Fix $K$ and
choose $1\le k\le K$ so that
$$
\frac{k-1}{K}\le T< \frac{k}{K}.
$$
So, $k-1 \le TK \le k$.  Now,
$$
\frac{1}{K}\sum_{j=0}^{K-1}|g(t-\frac{j}{K})|^2 \le
\frac{k}{K} = \frac{k-1}{K} + \frac{1}{K} \le T+ \frac{1}{K}.
$$
Similarly,
$$
T-\frac{1}{K}\le
\frac{1}{K}\sum_{j=0}^{K-1}|g(t-\frac{j}{K})|^2.
$$
This completes Step 1.

\noindent {\bf Step 2}:  If $\{g_j\}_{j=1}^{\ell}$ are disjointly
supported functions on $[0,1]$, $g=\sum_{j=1}^{\ell}g_j$ and all the
$g_j$ satisfy the theorem (respectively, satisfy the theorem
with uniform convergence), then $g$ satisfies the theorem (respectively,
with uniform convergence).

{\it Proof of Step 2}:  Since the $g_j$ are disjointly supported,
\begin{eqnarray*}
\lim_{n\rightarrow \infty}\frac{1}{K_n}
\sum_{k=0}^{K_n -1}|g(t-\frac{k}{K_n})|^2 &=&
\lim_{n\rightarrow \infty}\frac{1}{K_n}\sum_{k=0}^{K_n-1}
\sum_{j=1}^{\ell}|g_{j}(t-\frac{k}{K_n})|^2\\
&=& \sum_{j=1}^{\ell}\lim_{n\rightarrow \infty}\frac{1}{K_n}
|g_{j}(t-\frac{k}{K_n})|^2\\
&=& \sum_{j=1}^{\ell}\|g_j\|^2{\chi}_{[0,1]}\\
&=& \|g\|^2{\chi}_{[0,1]},\ \ \mbox{a.e. t (respectively, uniformly)}.
\end{eqnarray*}
This completes the proof of Step 2.

\noindent {\bf Step 3}:  Let $E$ be a measurable subset of $[0,1]$,
$\epsilon >0$, $K\in \N$ and set
$$
F = \{t\in [0,1]|\frac{1}{K}\sum_{k=0}^{K-1}|{\chi}_{E}(t-\frac{k}{K})|^2
\ge \epsilon \}.
$$
Let $F_1 = F\cap [0,\frac{1}{K}]$.  If $|F|\ge \epsilon$,
then $|E|\ge {\epsilon}^2$.

{\it Proof of Step 3}:
Let $E_{k} = E\cap [\frac{k}{K},
\frac{k+1}{K}]$ and compute
\begin{eqnarray*}
{\epsilon}^2 \le \epsilon |F| = \epsilon K|F_1| &\le&
\int_{0}^{\frac{1}{K}}|{\chi}_{E}(t-\frac{k}{K})|^2 dt\\
&=& \int_{0}^{\frac{1}{K}}\sum_{k=0}^{K-1}|{\chi}_{E}(t+\frac{k}{K})| dt\\
&=& \int_{0}^{\frac{1}{K}}\sum_{k=0}^{K-1}|{\chi}_{E_k}(t+\frac{k}{K})|dt\\
&=& \sum_{k=0}^{K-1}\int_{0}^{\frac{1}{K}}{\chi}_{E_k}(t+\frac{k}{K})dt\\
&=& \sum_{k=0}^{K-1}\int_{\frac{k}{K}}^{\frac{k+1}{K}}{\chi}_{E_k}(t)dt\\
&=& \sum_{k=0}^{K-1}|{\chi}_{E_k}| = |E|.
\end{eqnarray*}
This completes the proof of Step 3.

\noindent {\bf Step 4}:  If $E$ is a measurable set which is a
countable union of intervals, then the theorem holds.

{\it Proof of Step 4}:  Fix $\epsilon >0$ arbitrary and choose
$\delta >0$ to be specified later.  Assume that
$$
E = \cup_{n=1}^{\infty}(a_n ,b_n ),
$$
and choose a natural number $N$ so that
$$
\sum_{n=N+1}^{\infty}(b_n -a_n ) < {\delta}^2.
$$
Also let
$$
G= \cup_{n=N+1}^{\infty}(a_n,b_n ),
 \mbox{so that}\ \ |G|< {\delta}^2.
$$
Let
$$
H = \cup_{n=1}^{N} (a_n ,b_n ),
$$
and note that we may as well assume that $\{(a_n ,b_n )\}_{n=1}^{N}$
are disjoint.  Finally, let
$$
F = \{t\in [0,1]|\frac{1}{K}\sum_{k=0}^{K-1}
|{\chi}_{G}(t-\frac{k}{K})|^2 \ge \delta \},\ \ \mbox{so that}\ \
|F|< \delta .
$$
By Step 1,
$$
\lim_{K\rightarrow \infty}\frac{1}{K}\sum_{k=0}^{K-1}
|{\chi}_{G}(t-\frac{k}{K})|^2 = |G|{\chi}_{[0,1]},
$$
uniformly.  Choose a natural number $K_0$ so that for all
$K\ge K_0$ we have
$$
\frac{1}{K}\sum_{k=0}^{K-1}|{\chi}_{G}(t-\frac{k}{K})|^2
\le (|G|+\delta ){\chi}_{[0,1]}.
$$
Now we compute,
\begin{eqnarray*}
\sqrt{\frac{1}{K}\sum_{k=0}^{K-1}|{\chi}_{E}(t-\frac{k}{K})|^2}
&\le& \sqrt{\frac{1}{K}\sum_{k=0}^{K-1}|{\chi}_{G}(t-\frac{k}{K}) +
{\chi}_{H}(t-\frac{k}{K})|^2}\\
&\le& \sqrt{\frac{1}{K}\sum_{k=0}^{K-1}|{\chi}_{G}(t-\frac{k}{K})|^2}
+ \sqrt{\frac{1}{K}\sum_{k=0}^{K-1}|{\chi}_{H}(t-\frac{k}{K})|^2}\\
&\le& \sqrt{|H|+\delta}{\chi}_{[0,1]} +
\sqrt{\frac{1}{K}\sum_{k=0}^{K-1}|{\chi}_{G}(t-\frac{k}{K})|^2}\\
&\le& \sqrt{|H|+\delta }{\chi}_{[0,1]} + {\delta}{\chi}_{[0,1]},
\end{eqnarray*}
off of the set $F$ where $|F|< \delta < \epsilon$ and $\delta$ is
chosen so that
$$
\sqrt{|H|+\delta }{\chi}_{[0,1]}+ {\delta}{\chi}_{[0,1]}\le
\sqrt{|H|+\epsilon}{\chi}_{[0,1]}.
$$
Similarly, for $K\ge K_0$ we have
$$
\sqrt{\frac{1}{K}\sum_{k=0}^{K-1}|{\chi}_{E}(t-\frac{k}{K})|^2}
\ge \left ( \sqrt{|G|+\delta}-{\delta}  \right ) {\chi}_{[0,1]}
\ge \sqrt{|G|+\epsilon }{\chi}_{[0,1]}.
$$
This completes the proof of Step 4.

\noindent {\bf Step 5}:  The Theorem holds for ${\chi}_{E}$ for
every measurable set $E$.

{\it Proof of Step 5}:  Given a measurable set $E$ in $[0,1]$
and an $\epsilon >0$, choose intervals $\{(a_n,b_n)\}_{n=1}^{\infty}$
so that
$$
E\subset \cup_{n=1}^{\infty}(a_n ,b_n ) =: F,
$$
and
$$
||E|-\sum_{n=1}^{\infty}(b_n -a_n) |< \frac{\epsilon}{3}.
$$
Then, there is a measurable set $G$ with $|G|< {\epsilon}/3$ and
a natural number $K_0$ so that for every $K\ge K_0$ and for
all $t\notin G$ we have
\begin{eqnarray*}
\frac{1}{K}\sum_{k=0}^{K-1}|{\chi}_{E}(t-\frac{k}{K})|^2
&\le& \frac{1}{K}\sum_{k=0}^{K-1}|{\chi}_{F}(t-\frac{k}{K})|^2\\
&\le& (|F|+\frac{\epsilon}{3}){\chi}_{[0,1]} \le (|E|+\frac{2\epsilon}{3})
{\chi}_{[0,1]}.
\end{eqnarray*}
Similarly, there is a measurable set $G_1$
with $|G_1 |< {\epsilon}/3$ and a natural number $K_1$
so that for all $K\ge K_1$ and all $t\notin G_1$ we have
$$
\frac{1}{K}\sum_{k=0}^{K-1}|{\chi}_{E^c}(t-\frac{k}{K})|^2
\le (|E^c|+\epsilon ){\chi}_{[0,1]}.
$$
We next note that
\begin{eqnarray*}
{\chi}_{[0,1]}(t) &=& \frac{1}{K}\sum_{k=0}^{K-1}\left (
|{\chi}_{E}(t-\frac{k}{K})|^2 + |{\chi}_{E^c}(t-\frac{k}{K-1})|^2 \right )\\
&=& \frac{1}{K}\sum_{k=0}^{K-1}|{\chi}_{E}(t-\frac{k}{K})|^2
+ \frac{1}{K}\sum_{k=0}^{K-1} |{\chi}_{E^c}(t-\frac{k}{K-1})|^2.
\end{eqnarray*}
Now, for all $t\notin G_1$ we have
$$
\frac{1}{K}\sum_{k=0}^{K-1}|{\chi}_{E}(t-\frac{k}{K})|^2 \ge
1-\frac{1}{K}\sum_{k=0}^{K-1}|{\chi}_{E^c}(t-\frac{k}{K})|^2
\ge 1- (|E^c|+\epsilon ) = |E|- \epsilon.
$$
Hence, $|G\cup G_1| < \epsilon$ and for all $t\notin G\cup G_1$
we have
$$
|\frac{1}{K}\sum_{k=0}^{K-1}|{\chi}_{E}(t-\frac{k}{K})|^2 -
|E||< \epsilon.
$$
Hence,
$$
\{ \frac{1}{K}\sum_{k=0}^{K-1}|{\chi}_{E}(t-\frac{k}{K})|^2\}
$$
converges to $|E|{\chi}_{[0,1]}$ in measure.  This completes
the proof of Step 5.

\noindent {\bf Step 6}:  The general case for the theorem.

{\it Proof}:  If $g\in L^2 [0,1]$ and $\epsilon >0$ is given,
fix a $\delta >0$ (to be chosen later) and
choose a simple function
$$
h = \sum_{j=1}^{M}a_j {\chi}_{E_j},
$$
so that $|g - h|< \delta$ a.e.  Then,
\begin{eqnarray*}
\sqrt{\frac{1}{K}\sum_{k=0}^{K-1}|g(t-\frac{k}{K})|^2}
&\le& \sqrt{\frac{1}{K}\sum_{k=0}^{K-1}|h(t-\frac{k}{K})|^2}
+\sqrt{\frac{1}{K}\sum_{k=0}^{K-1}(\frac{\delta}{2})^2}\\
&\le& \sqrt{\frac{1}{K}\sum_{k=0}^{K-1}|h(t-\frac{k}{K})|^2}
+ \frac{\delta}{2}.
\end{eqnarray*}
Similarly,
$$
\sqrt{\frac{1}{K}\sum_{k=0}^{K-1}|g(t-\frac{k}{K})|^2}\ge
\sqrt{\frac{1}{K}\sum_{k=0}^{K-1}|h(t-\frac{k}{K})|^2}
- \frac{\delta}{2}.
$$
Now, there is a measurable set $|G|< \delta$ and a natural number
$K_0$ so that for all $K\ge K_0$ and all $t\notin G$ we have
\begin{eqnarray*}
\frac{1}{K}\sum_{k=0}^{K-1}|h(t-\frac{k}{K})|^2 &=&
\frac{1}{K}\sum_{k=0}^{K-1}\sum_{j=1}^{M}|a_j|^2|{\chi}_{E_j}
(t-\frac{k}{K})|^2\\
&=& \sum_{j=1}^{M}|a_j |^2 \frac{1}{K}\sum_{k=0}^{K-1}
|{\chi}_{E_j}(t-\frac{k}{K})|^2\\
&\le& \sum_{j=1}^{M}|a_j |^2 \left ( |E_j|+\frac{\delta}
{2M\sum |a_j |^2 } \right )\\
&=& \sum_{j=1}^{M}|a_j |^2 |E_j| + \frac{\delta}{2}=
\|h\|^2 + \frac{\delta}{2}.
\end{eqnarray*}
Hence,
\begin{eqnarray*}
\frac{1}{K}\sum_{k=0}^{K-1}|g(t-\frac{k}{K})|^2 &\le& (\|h\|^2 +
\frac{\delta}{2})+\frac{\delta}{2}\\
&\le& \|h\|^2 + \delta
\le (\|g\|+\frac{\delta}{2})^2 + \delta\\
&\le& \|g\|^2 + \frac{\delta}{2}\|g\| + \frac{{\delta}^2}{4}+\delta\\
&\le& \|g\|^2 + \epsilon,
\end{eqnarray*}
for an appropriately chosen $\delta>0$ and all $t\notin G$.
Similarly,
$$
\frac{1}{K}\sum_{k=0}^{K-1}|g(t-\frac{k}{K})|^2 \ge
\|g\|^2 - \epsilon,
$$
for all $t\notin G$.  This completes the proof of Step 6 and hence
of the theorem.
\qed

Now we proceed to the proof of the main theorems of this section.
For the proofs we will need a proposition.

\begin{proposition}
Fix $g\in L^{\infty}([0,1])$ and $K\in \N$.  For $0\le k\le K-1$
let
$$
g_{k}(t) = \sum_{n\in \Z}\langle g, \mathrm{e}^{2\pi\mathrm{i}(nK+k)t}\rangle
\mathrm{e}^{2\pi\mathrm{i}(nK+k)t}.
$$
Then, for every $f\in\mathrm{cl}(\mathrm{span}\{\mathrm{e}^{2\pi\mathrm{i}(nK+k)t}\}_{n\in \Z})$
we have:
$$
\|f\cdot g\|^2 =
\|f\cdot \left ( \sum_{k=0}^{K-1}|g_{k}|^2 \right ) ^{1/2} \|.
$$
\end{proposition}

{\it Proof}:
We will do the case $k=0$; the others require only notational
changes.  So, we compute:
\begin{eqnarray*}
\|f\cdot g\|^2 &=& \sum_{n\in \mathbb Z}|\langle f\cdot g,
\mathrm{e}^{2\pi\mathrm{i}nt} \rangle |^2\\
&=& \sum_{n\in \mathbb Z}|\langle f,\mathrm{e}^{2\pi\mathrm{i}nt}
\overline{g}\rangle |^2\\
&=& \sum_{k=0}^{K-1}\sum_{n\in \Z}|\langle f,\mathrm{e}^{2\pi\mathrm{i}(nK+k)t}
\overline{g}\rangle |^2\\
&=& \sum_{k=0}^{K-1}\sum_{n\in \Z}|\langle e^{-2{\pi}ikt}f,
\mathrm{e}^{2\pi\mathrm{i}nKt}\overline{g}\rangle |^2\\
&=& \sum_{k=0}^{K-1}\sum_{n\in \Z} |\langle e^{-2{\pi}inkt}f,\mathrm{e}^{2\pi\mathrm{i}nKt}
\overline{g_{k}}\rangle |^2\\
&=& \sum_{k=0}^{K-1}\sum_{n\in \Z}|\langle e^{-2{\pi}ikt}f,
\mathrm{e}^{2\pi\mathrm{i}nt}\overline{g_k}\rangle |^2\\
&=& \sum_{k=0}^{K-1}\sum_{n\in \Z}|\langle e^{-2{\pi}ikt}f\cdot g_k ,
\mathrm{e}^{2\pi\mathrm{i}nt}\rangle |^2\\
&=& \sum_{k=0}^{K-1}\|e^{-2{\pi}ikt}f\cdot g_{k}\|^2 =
\sum_{k=0}^{K-1}\|f\cdot g_k \|^2\\
&=& \sum_{k=0}^{K-1}\int_{0}^{1}|f(t)|^2|g_{k}(t)|^2dt\\
&=& \int_{0}^{1}|f(t)|^2\sum_{k=0}^{K-1}|g_{k}(t)|^2dt\\
&=& \|f \cdot \left ( \sum_{k=0}^{K-1} |g_k |^2 \right ) ^{1/2}\|^2.
\end{eqnarray*}
\qed

{\bf Proof of Theorem \ref{HTT1}}:  By Corollary \ref{CC1}, $T_g$ has
the uniform Kadison-Singer Property if and only if for every $\epsilon >0$
there is a natural number $K$ so that for all 
$f\in\mathrm{cl}(\mathrm{span}\{\mathrm{e}^{2\pi\mathrm{i}(nK+k)t}\}_{n\in \Z})$ we have:
$$
\|f\|^2 (\|g\|^2 - \epsilon ) \le \|f\cdot (\sum_{k=0}^{K-1}
|g_{k}^{K}(t)|^2 )^{1/2}\|^2
$$
$$
= \int_{0}^{1}|f(t)|^2 \frac{1}{K}\sum_{k=0}^{K-1}|g(t-\frac{k}{K})|^2dt
\le \|f\|^2(\|g\|^2 + \epsilon ).
$$

The proof of $(2)\Rightarrow (1)$ is immediate from here.

$(1)\Rightarrow (2)$:  If this implication fails, there is an
$\epsilon >0$ so that for all $K\in \N$ there is a measurable
set $E_0$ so that either
$$
\frac{1}{K}\sum_{k=0}^{K-1}|g(t-\frac{k}{K})|^2 \ge \|g\|^2 +
\epsilon ,
$$
or
$$
\frac{1}{K}\sum_{k=0}^{K-1}|g(t-\frac{k}{K})|^2 \le \|g\|^2 -
\epsilon .
$$
We will do the first case since the second is similar.  Choose
$\|h\|= \frac{1}{\sqrt{K}}$ so that $h = h\cdot {\chi}_{E_0}$.  Let
$E = \cup_{k=0}^{K-1}(E_0 + k)$ and let
$$
f = \frac{1}{K}\sum_{j=0}^{K-1}h(t-\frac{j}{K})\mathrm{e}^{\frac{2{\pi}\mathrm{i}jk}{K}}.
$$
Then, $g\in\mathrm{cl}(\mathrm{span}\{\mathrm{e}^{2\pi\mathrm{i}(nK+k)t}\}_{n\in \Z})$ and
\begin{eqnarray*}
\|f\cdot g\|^2 &=&
\int_{0}^{1}|f(t)|^2 \frac{1}{K}\sum_{k=0}^{K-1}|g(t-\frac{k}{K})|^2dt\\
&ge& \int_{E}|f(t)|^2 dt \left ( \|g\|^2 + \epsilon \right )\\
&=&
\|g\|^2 + \epsilon .
\end{eqnarray*}
Similarly, for the other case we have
$$
\|f\cdot g\|^2 =
\int_{0}^{1}|f(t)|^2 \frac{1}{K}\sum_{k=0}^{K-1}|g(t-\frac{k}{K})|^2dt
\le \|g\|^2 - \epsilon.
$$
So the Paving Conjecture fails for $T_g$ by Corollary \ref{CC1}.

{\bf Proof of Theorem \ref{HTT2}}: This is similar to the above.
The Toeplitz operator $T_g$ has the uniform Feichtinger property if
and only if there exists an $\epsilon >0$ and a natural number $K$
so that for all $f\in\mathrm{cl}(\mathrm{span}\{\mathrm{e}^{2\pi\mathrm{i}(nK+k)t}\}_{n\in \Z})$
we have
$$
\|f\cdot g\|^2 = \int_{0}^{1}|f(t)|^2 \frac{1}{K}\sum_{k=0}^{K}
g(t-\frac{k}{K})|^2 dt \ge \epsilon .
$$
As in the proof of Theorem \ref{HTT1}, this holds if and only if
there exists an $\epsilon >0$ and a $K$ so that
$$
\frac{1}{K}\sum_{k=0}^{K}|g(t-\frac{k}{K})|^2 \ge \epsilon \ \ a.e.
$$
This shows that $(1)\Leftrightarrow (2)$.

$(3)\Rightarrow (2)$:  By (3) we have that
$$
\frac{1}{K}\sum_{k=0}^{K}|g(t-\frac{k}{K})|^2 \ge \epsilon \ \ a.e.
$$

$(2)\Rightarrow (3)$:  If
$$
\frac{1}{K}\sum_{k=0}^{K}|g(t-\frac{k}{K})|^2 \ge \epsilon \ \ a.e.,
$$
then for all $0\le k\le K-1$ let
$$
F_k = \{t\in [\frac{k}{K},\frac{k+1}{K}]\ | \ |g(t)|\ge \epsilon \}.
$$
Now,
$$
\cup_{k=0}^{K}(F_k -k) = [0,\frac{1}{K}].
$$
Letting $E_0 = F_0$ and
$$
E_{k+1} = F_k \setminus \cup_{j=0}^{k-1}(F_j - j),
$$
produces the desired sets.  This completes the proof of Theorem \ref{HTT2}.
\qed

\section{Kadison-Singer in time-frequency analysis}\label{TFA}
\setcounter{equation}{0}

Although the Fourier transform has been a major tool
in analysis for over a century, it has a serious lacking
for signal analysis in that it hides in its phases information
concerning the moment of emission and duration of a signal.
What was needed was a localized time-frequency representation
which has this information encoded in it.  In 1946 Gabor
\cite{Ga} filled this gap and formulated a fundamental approach
to signal decomposition in terms of elementary signals.  Gabor's
method has become the paradigm for signal analysis in Engineering
as well as its mathematical counterpart: Time-Frequency Analysis.

To build our elementary signals, we choose a {\bf window function}
$g\in L^2 (\R)$.  For $x,y\in \R$
we define {\bf modulation by x} and {\bf translation by y} of $g$ by:
$$
M_{x}g(t) = \mathrm{e}^{2\pi\mathrm{i}xt}g(t),\ \ T_{y}g(t) = g(t-y).
$$
If $\Lambda \subset \R\times \R$
and $\{E_{x}T_{y}g\}_{(x,y)\in \Lambda}$ forms a frame for
$L^2(\R)$ we call this an (irregular) {\bf Gabor frame}.
Standard Gabor frames are the case where $\Lambda$
is a lattice $\Lambda = a\Z \times {b \Z}$ where $a,b>0$ and
$ab\le 1$.  For an introduction to time-frequency analysis
we recommend the excellent book of Grochenig \cite{G1}.

It was in his work on time-frequency analysis that Feichtinger
observed that all the Gabor frames he was working with could
be decomposed into a finite union of Riesz basic sequences.
This led him to formulate the Feichtinger Conjecture - which
we now know is equivalent to KS.
There is a significant amount of literature on the Feichtinger
Conjecture for Gabor frames as well as wavelet frames and
frames of translates \cite{BCHL,BCHL2,BS}.  It is known that Gabor frames
over rational lattices \cite{CCLV} and Gabor frames whose window function
is ``localized'' satisfy the Feichtinger Conjecture \cite{BCHL,BCHL2,G}.
But the
general case has defied solution.

Translates of a single function play a fundamental role in
frame theory, time-frequency analysis, sampling theory
and more \cite{Al,BS}.
If $g\in L^{2}(\R)$, ${\lambda}_n \in \R$ for $n\in \Z$ and
$\{T_{{\lambda}_n}g\}_{n\in \Z}$ is a frame for its closed
linear span, we call this a {\bf frame of translates}.  Although
considerable effort has been invested in the Feichtinger Conjecture
for frames of translates, little progress has been made.  One
exception is a surprising result from \cite{CCK}.

\begin{theorem}
Let $I \subset \Z$ be bounded below, $a>0$ and $g\in L^2 (\R)$.
Then $\{T_{na}g\}_{n\in I}$ is a frame sequence if and only if it is
a Riesz basic sequence.
\end{theorem}

Our next theorem will explain why the Feichtinger Conjecture for
frames of translates, wavelet frames and Gabor frames has proven
to be so difficult.  This is due to the fact that this problem 
is equivalent to a deep problem in harmonic analysis, 
namely Conjecture \ref{C11}, which in turn is equivalent to
having all Toeplitz operators satisfy the Feichtinger Conjecture
(Theorem \ref{HAT12}).

The proof our of our theorem is complicated,
and requires some preliminary work.
The main idea is to apply the Fourier transform to turn this
into a problem concerning functions of the form
$\{\mathrm{e}^{2\pi\mathrm{i}{\lambda}_nt}{\phi}\}_{n\in \Lambda}$
with $\phi \in L^2(\R)$.  Then we want to
use perturbation theory to reduce this problem into
 one with evenly spaced exponentials $\{\mathrm{e}^{2\pi\mathrm{i}nt}{\phi}\}_{n\in \Lambda}$.
There are two technical problems with this.  The first is that our
functions $\phi$ are no longer in $L^{\infty}[0,1]$ which causes
technicalities.  Second,
 perturbation theory fails miserably in the
frame setting if we perturb a frame by a sequence from outside the
space - as we have to do here (see Example
\ref{TFAE5} below).  What makes this all eventually work
is that perturbation theory {\it does work from outside the space}
for Riesz basic sequences and we are just trying to divide our family
of vectors into a finite number of Riesz basic sequences.

\begin{example}\label{TFAE5}
Let $\{e_i\}_{i=1}^{\infty}$ be an orthonormal basis for
${\ell}_2$.  For all $i\in \N$, define $g_{2i} = g_{2i+1} = e_{2i}$,
$f_{2i+1} = e_{2i}$ and
$$
f_{2i} = e_{2i}+ \frac{\epsilon}{2^{i+1}}e_{2i+1}.
$$
Then $\{g_{i}\}_{i=1}^{\infty}$ is clearly a 2-tight frame
for the span of $\{e_{2i}\}_{i=1}^{\infty}$.  Also, for any finitely
non-zero sequence of scalars $\{a_i\}_{i=1}^{\infty}$ we have
$$
\|\sum_{i=1}^{\infty}a_i(g_i-f_i)\|^2 =
\|\sum_{i=1}^{\infty} \frac{\epsilon}{2^{i+1}}a_{i+1}e_{i+1}\|^2
\le \epsilon  \sum_{i=1}^{\infty}|a_i|^2 \le \epsilon \sum_{i=1}^{\infty}
|a_i|^2.
$$
So $\{f_i\}_{i=1}^{\infty}$ is a small perturbation of $\{g_i\}_{i=1}^{\infty}$
but $\{f_i\}_{i=1}^{\infty}$ is not a frame for its span since for any
$j\in \N$ we have
$$
\sum_{i=1}^{\infty}|\langle e_{2j+1},f_i \rangle |^2 =
\frac{\epsilon}{2^{i+1}}.
$$
\end{example}

We will state
the main theorems here, then develop some theory for solving
them and give the proofs at the end.

\begin{theorem}\label{TFAP5}
The following are equivalent:

(1)  Conjecture \ref{HAT12}.

(2)  For every
$0\not= \phi \in L^{2}(\mathbb R )$ and ${\lambda}_{n}\in \R$
for $n\in \Lambda$, if
$\{T_{{\lambda}_n}\phi \}_{n\in \Lambda}$ is a
Bessel sequence, then it is a finite
union of Riesz basic sequences.

(3)  For every $\Lambda \subset \Z$ and every
$0\not= \phi \in L^{2}(\mathbb R )$, if
$\{T_{n}\phi \}_{n\in \Lambda}$ is a Bessel sequence, then it
is a finite union of Riesz basic sequences.

(4)  For every
$0\not= \phi \in L^{2}(\mathbb R )$ and ${\lambda}_{n}\in \R$
for $n\in \Lambda$, if
$\{T_{{\lambda}_n}\phi \}_{n\in \Lambda}$ is a
frame sequence, then it is a finite
union of Riesz basic sequences.

(5) For every $\Lambda \subset \Z$ and every
$0\not= \phi \in L^{2}(\mathbb R )$, if
$\{T_{n}\phi \}_{n\in \Lambda}$ is a Bessel sequence, then it
is a finite union of frame sequences.
\end{theorem}

Instead of proving Theorem \ref{TFAP5}, we will take the Fourier
transform of all this and prove the equivalent formulation given
in the next theorem.

\begin{theorem}\label{TFAP6}
The following are equivalent:

(1)  Conjecture \ref{HAT12}.

(2)  For every $\phi \in L^{2}(\R )$ and every
$\{{\lambda}_{n}\}_{n\in \Lambda}$, if
$\{\mathrm{e}^{2\pi\mathrm{i}{\lambda}_{n}t}\phi \}_{n\in \Lambda}$ is Bessel in
$L^{2}(\R )$, then it is a finite union of Riesz basic sequences.

(3)  For every $0\not= \phi \in L^{2}[0,1]$ and every
$\Lambda \subset \Z$, if $\{\mathrm{e}^{2\pi\mathrm{i}nt}\phi \}_{n\in \Lambda}$
is a Bessel sequence then it is a finite union of Riesz basic
sequences.

(4)  For every $\phi \in L^{2}(\R )$ and every
$\{{\lambda}_{n}\}_{n\in \Lambda}$, if
$\{\mathrm{e}^{2\pi\mathrm{i}{\lambda}_{n}t}\phi \}_{n\in \Lambda}$ is a frame sequence in
$L^{2}(\R )$, then it is a finite union of Riesz basic sequences.

(5)  For every $\Lambda \subset \Z$ and every
$0\not= \phi \in L^{2}[0,1]$, if
$\{\mathrm{e}^{2\pi\mathrm{i}nt}\phi \}_{n\in \Lambda}$ is a Bessel sequence, then it
is a finite union of frame sequences.
\end{theorem}

The first thing we will do is derive the perturbation theorem
we need for proving our results.  We start with a theorem due to
Christensen \cite{Ch2} which is a generalization of the Paley-Wiener
theorem \cite{K} (We state a slightly stronger conclusion at
the end which easily follows from the proof of \cite{Ch2}).

\begin{theorem}\label{CHRT1}
Let $\H$ be a Hilbert space and $\{f_i\}_{i\in I}$ a frame for
$\H$ with frame bounds $A,B$.  Let $\{g_i\}_{i\in I}$
be a sequence in $\H$. Assume there exists a ${\lambda},\mu>0$
with $\lambda + \frac{\mu}{\sqrt{A}} < 1$
and an increasing sequence of subsets $I_1 \subset I_2 \subset
\cdots \subset I$ with $\cup_{n=1}^{\infty}I_n = I$ so that
for all $n=1,2,\ldots$ and all families of scalars $\{a_i\}_{i\in I_n}$
we have
$$
\|\sum_{i\in I_n}a_i (f_i -  g_i)\| \le
{\lambda}\|\sum_{i\in I_n}a_i f_i\| + {\mu} \left (
\sum_{i\in I_n}|a_i|^2 \right )^{1/2}.
$$
Then $\{g_i\}_{i\in I}$ is a frame for $\H$ with frame bounds
$$
A(1-{\lambda} - \frac{\mu}{\sqrt{A}})^2,\ \ \ \
B(1+{\lambda}+ \frac{\mu}{\sqrt{B}})^2.
$$
Moreover, if $\{f_i\}_{i\in I}$ is a Riesz basic sequence, then
$\{g_i\}_{i\in I}$ is also a Riesz basic sequence.
\end{theorem}

We will need a variation of a result
proved independently by Balan \cite{Ba} and
Christensen \cite{Ch3}.  Since this
is a straightforward generalization
where we just insert a function into the
calculations of Balan \cite{Ba}, we will outline the proof.

\begin{theorem}\label{RBalan}
Let $\phi \in L^2[-\gamma, \gamma]$ and assume
 $\{e^{i{\lambda}_nt}\phi\}_{n\in \Z}$ is a Bessel sequence with
Bessel bound $B$.  Set
$$
L(\gamma) = \frac{\pi}{4\gamma}-\frac{1}{\gamma} \arcsin \left (
\frac{1}{\sqrt{2}}\left ( 1 - \sqrt{\frac{A}{B}}\right ) \right ).
$$
Suppose ${\mu}_n \in \R$ and $\sup_n |{\mu}_n
- {\lambda}_n| = \delta < 1/4$.  Then $\{e^{i{\mu}_nt}\phi \}_{n\in I}$
is a Bessel sequence with Bessel bound $B(2- \cos \ {\gamma}\delta
+ \sin\ {\gamma}\delta )^2$.  Moreover, if
$\{e^{i{\lambda}_n t}\phi\}_{n\in I}$
is a frame sequence with frame bounds $A,B$ and $\delta < L(\gamma)$ then
$\{e^{i{\mu}_nt}\phi \}_{n\in I}$ is a frame sequence with frame bounds
$$
A\left ( 1- \sqrt{\frac{A}{B}}(1- \cos\ {\gamma}\delta + \sin \ {\gamma}\delta
) \right )^2,\ \ B(2- \cos \ {\gamma}\delta
+ \sin\ {\gamma}\delta )^2.
$$
\end{theorem}

{\it Proof}:
By a change of scale ($\hat{\lambda}_n = \frac{\gamma}{\pi}{\lambda}_n$)
we may assume $\gamma = \pi$ and we need to show:
$$
L(\pi) = \frac{1}{4} - \frac{1}{\pi}\arcsin \left [ \frac{1}{\sqrt{2}}
\left ( 1- \sqrt{\frac{A}{B}}\right ) \right ].
$$
To prove this result, we rely on Kadec's classical estimations for
computing the Paley-Wiener constant \cite{Kad}.  Let $\{a_n\}_{n\in I}$
be scalars, $I_N \subset I$ with $|I_N| <\infty$ and let
${\delta}_n = {\mu}_n - {\lambda}_n$.  We compute:
\begin{equation}\label{Bal1}
U=: \big \| \sum_{n\in I_N}a_n \left ( \frac{1}{\sqrt{2\pi}}
e^{i{\lambda}_n t} - \frac{1}{\sqrt{2\pi}}e^{i{\mu}_nt} \right ) \phi (t)
\big \| = \frac{1}{\sqrt{2\pi}}\big \|
\sum_{n\in I_n}a_n e^{i{\lambda}_nt}(1-e^{i{\delta}_nt})\phi (t) \big \|.
\end{equation}
By expanding $1-e^{i{\delta}_nt}$ into a Fourier series relative to
the orthogonal system $\{1, \cos\ \nu t,\sin(\nu - \frac{1}{2})t\}$,
$\nu = 1,2,\ldots$ we have
\begin{eqnarray}\label{Bal2}
1- e^{i{\delta}_nt} &=& \left ( 1- \frac{\sin\ \pi{\delta}_n}
{{\pi}{\delta}_n} \right ) + \sum_{\nu =1}^{\infty}
\frac{(-1)^{\nu}2{\delta}_n\sin\ \pi{\delta}_n}{\pi({\nu}^2 -
{\delta}_n^2)}\cos\ (\nu t)\\
&+& i \sum_{\nu=1}^{\infty}\frac{(-1)^{\nu}2{\delta}_n \cos \ \pi{\delta}_n}
{\pi(({\nu}-\frac{1}{2})^2 - {\delta}_n^2)}\sin \left (
\left ( \nu - \frac{1}{2}\right ) t \right ).
\end{eqnarray}
We next insert \eqref{Bal2} into \eqref{Bal1}, change the order of summation,
apply the triangle inequality and use the bounds
$\|(\cos\ \nu t)f (t)\|\le \|f\|$ and
$\|(\sin(\nu - \frac{1}{2})t) f(t)\|\le \|f\|$ to arrive at
\begin{eqnarray*}
U &\le& \big \| \sum_{n\in I_N}\left ( 1- \frac{\sin\ \pi{\delta}_n}
{\pi{\delta}_n}\right ) a_n e^{i{\lambda}_n t} \phi (t)\big \|
+\sum_{\nu =1}^{\infty}  \big \| \sum_{n\in I_N}
\frac{2{\delta}_n \ \sin\ \pi{\delta}_n}{\pi ({\nu}^2 - {\delta}_n^2)}
a_n e^{i{\lambda}_n t}\phi (t) \big \| \\
&+& \sum_{\nu =1}^{\infty}
\big \| \sum_{n\in I_N}\frac{2{\delta}_n\ \cos\ \pi{\delta}_n}
{\pi (({\nu}- \frac{1}{2})^2 - {\delta}_n^2)}
a_ne^{i{\lambda}_n t} \phi (t)\big \| .
\end{eqnarray*}
Now we use the fact that $\{e^{i{\lambda}_n t}\phi(t)\}_{n\in I}$
is a $B$-Bessel sequence.  Therefore, each norm above can be bounded
as:
$$
\big \| \sum_{n\in I_N}c_na_n e^{i{\lambda}_nt}\phi (t)\big \|
\le \sqrt{B}\|\{c_n a_n\}_{n\in I_N}\| \le \sqrt{B} \sum_{n\in I_N}
|c_n| \|\{a_n\}_{n\in I_N}\|.
$$
Also,
$$
\left | 1- \frac{\sin\ \pi{\delta}_n}{\pi {\delta}_n} \right |
\le 1- \frac{\sin\ \pi{\delta}}{\pi \delta},
$$
$$
\left | \frac{2{\delta}_n \sin\ \pi{\delta}_n}{\pi ({\nu}^2 -
{\delta}_n^2)}\right | \le \frac{2{\delta}\sin\ \pi{\delta}}
{\pi ({\nu}^2 - {\delta}^2)},
$$
$$
\left | \frac{2{\delta}_n \cos\ \pi {\delta}_n}{\pi (({\nu}-
\frac{1}{2})^2 - {\delta}_n^2)}\right | \le
\frac{2 \delta \cos\ \pi \delta}{\pi((\nu - \frac{1}{2})^2
-{\delta}^2)},
$$
where the last inequality holds because $\delta < \frac{1}{4}$.
Thus,
$$
U\le \sqrt{B}(Re(1-e^{i{\pi}\delta})-Im(1-e^{i{\pi}\delta}))
\left ( \sum_{n\in I_N}|a_n|^2 \right )^{1/2}.
$$
That is,
$$
U \le \sqrt{B}(1- \cos\ \pi \delta + \sin\ \pi \delta)\left (
\sum_{n\in I_n}|a_n|^2 \right )^{1/2}.
$$
Now we apply Theorem \ref{CHRT1} with $\lambda = 0$ and $\mu =
\sqrt{B}(1-\cos\ \pi \delta + \sin\  \pi \delta )$.
The condition of that theorem becomes $\mu < \sqrt{A}$ or
$1- \cos\ \pi \delta + \sin\ \pi \delta < \sqrt{\frac{A}{B}}$.
Standard trigonometry yields
$$
\delta < L = \frac{1}{4} - \frac{1}{\pi}\arcsin\
\left ( \frac{1}{\sqrt{2}}\left ( 1- \sqrt{\frac{A}{B}}\right ) \right ).
$$
The frame bounds come from $A(1-\frac{\mu}{\sqrt{A}})^2$
and $B(1+\frac{\mu}{\sqrt{B}})^2$.  This completes the proof.
\qed

We also need a simple observation.

\begin{lemma}\label{LLL}
Suppose $\{{\phi}_n \}_{n\in \Lambda}$ is a Riesz basic sequence
in $L^2(I)$, $I\subset \mathbb R$ with Riesz basis bounds $A,B$.
If $|\phi |=1$ a.e. then
$\{{\phi}_n \phi\}_{n\in \Lambda}$ is a
 Riesz basic sequence with Riesz basis bounds
$A,B$.
\end{lemma}

{\it Proof}:
For any sequence of scalars $\{a_n\}_{n\in \Lambda}$ we have
\begin{eqnarray*}
\|\sum_{n\in \Lambda} a_n {\phi}_n \phi \|^2 &=&
\int_{I}|\sum_{n\in \Lambda} a_n {\phi}_{n}(t) \phi (t)|^2 dt \\
&=& \int_{I}|\sum_{n\in \Lambda} a_n {\phi}_{n}(t)|^2 dt\\
&=& \|\sum_{n\in \Lambda} a_n {\phi}_{n}\|^2.
\end{eqnarray*}
\qed

\begin{corollary}\label{CCC2}
Suppose $\phi \in L^2(\R)$,
$\{\mathrm{e}^{2\pi\mathrm{i}nt}\phi \}_{n\in \Lambda}$ is a Riesz basic
sequence
in $L^2 [0,1]$ and ${\lambda}_{n}\in \mathbb R$ with $|n-
{\lambda}_{n}|<1$
for every $n\in \Lambda$.  Then
$\{\mathrm{e}^{2\pi\mathrm{i}{\lambda}_{n}t}\phi \}_{n\in \Lambda}$ is a finite union of
Riesz basic sequences.
\end{corollary}

{\it Proof}:
Suppose $\{\mathrm{e}^{2\pi\mathrm{i}nt}\phi \}_{n\in \Lambda}$ has Riesz basis bounds
$A,B$.  Choose $N$ so that $1/(2{\pi}N) < 1/4$ and
$$
\sqrt{\frac{B}{A}}(1- \cos \ \frac{1}{N} + \sin\ \frac{1}{N} ) < 1.
$$
For $j=0,1,\ldots , N-1$ let
$$
{\Lambda}_{j} = \{ n\in \Lambda | n+\frac{j}{N} \le {\lambda}_{n} \le
n + \frac{j+1}{N} \},
$$
and
$$
{\mu}_{n} = n+ \frac{j}{N},\ \ \mbox{for $n\in {\Lambda}_n$}.
$$
Then
$$
\{\mathrm{e}^{2\pi\mathrm{i}{\mu}_{n}t}\phi \}_{n\in {\Lambda}_{j}} =
\{\mathrm{e}^{2\pi\mathrm{i}nt}{\phi}\mathrm{e}^{2\pi\mathrm{i}\frac{j}{N}t}\}_{n\in {\Lambda}_{j}},
$$
is a Riesz basic sequence with Riesz basis bounds $A,B$, by Lemma
\ref{LLL}.  Since $|{\mu}_{n}-{\lambda}_{n}|< \frac{1}{N}$
by Theorem \ref{RBalan} (rescaled to this setting) we have
$$
\|\sum a_k (\mathrm{e}^{2\pi\mathrm{i}{\mu}_{n}t}-\mathrm{e}^{2\pi\mathrm{i}{\lambda}_{n}t})\phi \|
\le \sqrt{B}(1- \cos\ \frac{1}{N} + \sin\ \frac{1}{N}),
$$
and
$$
\sqrt{\frac{B}{A}}(1- \cos\ \frac{1}{N} + \sin\ \frac{1}{N}) < 1.
$$
So by Corollary \ref{CCC2},
$$
\{\mathrm{e}^{2\pi\mathrm{i}{\lambda}_{n}t}\phi\}_{n\in {\Lambda}_{j}}
$$
is a Riesz basic sequence for $j=0,1,\ldots , N-1$.
\qed

We will need a little more notation.  If $\phi \in L^{2}(\mathbb R )$ we
define
${\Phi}_{b}:[0,1] \rightarrow \mathbb R$ by
$$
{\Phi}_{b}(t ) = \sum_{n\in \Z}\left | \hat{\phi}\left (
\frac{t + n}{b}\right ) \right | ^2.
$$
If $\Lambda \subset \Z$ we let $S({\Lambda})$ be the closed subspace of
$L^{2}([0,1])$ generated by the characters $\mathrm{e}^{2\pi\mathrm{i}n{t}}$
for $n\in \Lambda$.  We let $E_{\Lambda}$ be the closed subspace of
$S({\Lambda})$ of all $f$ such that ${\Phi}_{b}(t )f(t ) = 0$
a.e.  If $f\in S({\Lambda})$ we denote by $d(f,E_{\Lambda})$
the distance of $f$ from the subspace $E_{\Lambda}$.  We denote
$T_{x}$ the translation operator on $L^{2}(\mathbb R )$ by $x$.
Now we can state
the result from \cite{CCK}.

\begin{theorem}\label{CCK5}
Suppose $\phi \in L^{2}(\mathbb R )$ and $b>0$.  If $\Lambda \subset \Z$
then $\{T_{nb}\phi \}_{n\in \Lambda}$ is a frame sequence with frame
bounds $A,B$ if and only if for every $f\in S({\Lambda})$ we have
$$
Ad(f,E_{\Lambda})^2 \le \frac{1}{b}\int_{0}^{1}|f(t )|^2
{\Phi}_{b}(t )d{t} \le B\|f\|^2,
$$
or equivalently, for all $f\in S({\Lambda})\cap E_{\Lambda}^{\perp}$,
$$
A\|f\|^2   \le \frac{1}{b}\int_{0}^{1}|f( t )|^2
{\Phi}_{b}(t )d{t} \le B\|f\|^2 .
$$
Furthermore, if this condition is satisfied,
$\{T_{nb}{\phi}\}_{n\in \Lambda}$ is a Riesz basic sequence with
the same frame bounds if and only if $E_{\Lambda} = \{0\}$.
\end{theorem}

Now we are ready to prove our main theorem.

\noindent {\bf Proof of Theorem \ref{TFAP6}}:

$(3)\Rightarrow (2)$:  We first note the existence
of a natural number $N$ so that any interval of length
one in $\R$ contains at most $N$ of the ${\lambda}_n's$.
If not, then for every natural number $N$ there is a set
$I\subset \Lambda$ with $|I|=N$ and
$$
|\langle \mathrm{e}^{2\pi\mathrm{i}{\lambda}_{n}t}\phi , \mathrm{e}^{2\pi\mathrm{i}{\mu}_{n}t}\phi
\rangle |^2 \ge \frac{1}{2}\|\phi\|^2,
$$
for all ${\lambda}_n ,{\mu}_n \in I$.  Now, for ${\lambda}_n \in I$
fixed we have
$$
\sum_{{\mu}_n \in I}
|\langle \mathrm{e}^{2\pi\mathrm{i}{\lambda}_{n}t}\phi , \mathrm{e}^{2\pi\mathrm{i}{\mu}_{n}t}\phi
\rangle |^2 \ge \frac{N}{2}\|\phi\|^2.
$$
So $\{\mathrm{e}^{2\pi\mathrm{i}{\lambda}_{n}t}\}_{n\in \Lambda}$ is not Bessel which
contradicts our assumption.  It follows that we
can write $\Lambda$ as a finite union of sets so that
$|{\lambda}_{n}-{\lambda}_{m}|\ge 1$ for all $n\not= m$.  So we
may just assume that $\Lambda$ has this property.  By reindexing,
we may assume there is some $\Lambda \subset \Z$ so that for
$n\in \Lambda$ we have $|{\lambda}_{n} - n|<1$.  Let
$\{{\Lambda}_{j}\}_{j=0}^{4}$ be:
$$
{\Lambda}_{i} = \{n\in \Lambda | n+\frac{i}{5}\le {\lambda}_{n} <
n+\frac{i+1}{5}\}.
$$
Fix $0\le j \le 4$.  Since
$\{\mathrm{e}^{2\pi\mathrm{i}{\lambda}_{n}t}\phi \}_{n\in {\Lambda}_{j}}$ is Bessel with
Bessel bound $B$, by Theorem \ref{RBalan} we have that
$\{\mathrm{e}^{2\pi\mathrm{i}nt}\phi \}_{n\in {\Lambda}_{j}}$ is Bessel with Bessel
bound
$$
B_1 = B[2-\cos\ \frac{\pi}{5} + \sin \frac{\pi}{5} ].
$$
By our assumption (3), we can partition ${\Lambda}_{j}$ into a finite
number of sets $\{{\Lambda}_{jk}\}_{k=1}^{M_j}$ so that for every
$k=1,2,\ldots , M_j$, the family
$\{\mathrm{e}^{2\pi\mathrm{i}nt}\phi \}_{n\in {\Lambda}_{jk}}$ is a Riesz basic sequence
with some lower Riesz basis bound $A>0$.  By Corollary \ref{CCC2}, we have that
$\{\mathrm{e}^{2\pi\mathrm{i}{\lambda}_{n}t}\phi\}_{n\in {\Lambda}_{jk}}$ is a finite
union of Riesz basic sequences.

$(2)\Rightarrow (4)$:  This is obvious.

$(4)\Rightarrow (1)$:  Since $\{\mathrm{e}^{2\pi\mathrm{i}nt}{\chi}_{E}\}_{n\in \Z}$
is a Parseval frame, by (4) there is a partition $\Aj$ of $\Z$
so that $\{\mathrm{e}^{2\pi\mathrm{i}nt}\phi \}_{n\in A_j}$ is a Riesz basic sequence
(with lower Riesz basis bound $A>0$) for all $j=1,2,\ldots ,M$.
Hence, for any $f=\sum_{n\in A_j}a_n \mathrm{e}^{2\pi\mathrm{i}nt}$, we have that
$\|f\|^2 = \sum_{n\in A_j}|a_n |^2$ and
$$
\|P_{E}f\|^2 = \|\sum_{n\in A_j}a_n \mathrm{e}^{2\pi\mathrm{i}nt}{\chi}_{E}\|^2
\ge A^2
\sum_{n\in A_j}|a_n |^2  = A^2\|f\|^2.
$$
That is, $P_{E}$ is an isomorphism onto its range.

$(1)\Rightarrow (3)$:  Suppose $\{\mathrm{e}^{2\pi\mathrm{i}nt}\phi \}_{n\in \Lambda}$
is Bessel in $L^2 (\R )$.  So there exists a $B>0$ so that
for all $f\in H_{\Lambda}$ we have
$$
\int_{0}^{1}|f(t)|^2 {\Phi}(t )d{t} \le B\|f\|^2.
$$
Since $\phi \not= 0$, ${\Phi}\not= 0$.  So there is a measurable set
$E\subset [0,1]$ with $0<|E|$ and an $\epsilon >0$ so that
$|{\Phi}(t )|\ge \epsilon$ for all $t \in E$.  By the above,
$\{\mathrm{e}^{2\pi\mathrm{i}nt}\Phi \}_{n\in \Lambda}$ is a bounded Bessel sequence
in $L^2 [0,1]$.  By (1), there is a partition of $\Z \cap \Lambda$ of
the form $\{A_j\}_{j=1}^{M}$ so that $P_{E}$ is an isomorphism on $S(A_j )$, for every
$j=1,2,\ldots , M$ with lower isomorphism bound $A>0$.  Hence, for every
$\{a_n \}_{n\in A_j}$ we have
$$
\|\sum_{n\in A_j}a_n \mathrm{e}^{2\pi\mathrm{i}nt}\Phi \|_{L^{2}[0,1]} \ge
\|\sum_{n\in A_j}a_n \mathrm{e}^{2\pi\mathrm{i}nt}{\chi}_{E}\|_{L^2 [0,1]}
\ge A \left ( \sum_{n\in A_j }|a_n |^2 \right ) ^{1/2}.
$$
So $\{\mathrm{e}^{2\pi\mathrm{i}nt}\Phi \}_{n\in A_j}$ is a Riesz basic sequence for
all $j=1,2,\ldots ,M$.

$(5)\Rightarrow (3)$:  Assume $\{\mathrm{e}^{2\pi\mathrm{i}nt}\phi \}_{n\in \Lambda}$
is Bessel.  Let
$$
{\Lambda}^{+}= \{n\in \Lambda | 0\le n\},\ \
{\Lambda}^{-} = \{n\in \Lambda | n< 0\}.
$$
Now, $\{\mathrm{e}^{2\pi\mathrm{i}nt}\phi \}_{n\in {\Lambda}^{+}}$ is Bessel with the
same Bessel bound.  So by (5), we can write it as a finite union
of frame sequences.  By Proposition
\ref{Prop6} and by Theorem \ref{CCK5}, these frame
sequences are all Riesz basic sequences.

$(3)\Rightarrow (5)$:  This is obvious.

This completes the proof of the theorem.
\qed

We end this section with a result of Bownik and Speegle \cite{BS} which
makes a connection between number theory and PC for Toeplitz
operators.  This is related to a possible generalization of
van der Waerden's theorem \cite{Mont,MV}.

\begin{definition}
Let $g:{\R}^2 \rightarrow \R$.  We say that $I\subset \Z$ satisfies the
$g(\ell ,N)$-arithmetic progression condition if for every $\delta >0$
there exists $M,N,\ell \in \Z$ such that

(1)  $g(\ell , N)<\delta$, and

(2)  $\{M,M+\ell ,\ldots ,M+N{\ell}\}\subset I$.
\end{definition}

Taking the Fourier transform through
 theorem  4.1.2 in \cite{BS} yields:

\begin{theorem}\label{BST1}
A positive solution to the Feichtinger Conjecture for Toeplitz operators
implies there is a partition $\{A_j\}_{j=1}^{r}$ of $\Z$ so that
each $A_j$ fails the $g(\ell ,N) = {\ell}N^{-1/2}log^3\ N$
arithmetic progression condition.
\end{theorem}

In \cite{BS}, it is observed then
if we randomly assign each integer to one of $r$ subsets
$\{A_j\}_{j=1}^{r}$ of $\Z$, then with probability one,
for each $j$ and $r$ there will exist an $M_j$ such that
$$
\{M_j,M_j+1,\ldots ,M_j+L\} \subset A_j.
$$
Now, if $\Z$ is partitioned as $\{A_j\}_{j=1}^{r}$, the probability
that $\{T_n \phi\}_{n\in A_j}$
(and hence $\{\mathrm{e}^{2\pi\mathrm{i}nt}\hat{\phi}\}_{n\in A_j}$)
is a Riesz basic sequence is zero.

\section{Kadison-Singer in Engineering}\label{Eng}
\setcounter{equation}{0}

Frames have traditionally been used in signal processing because
of their resilience to additive noise, resilience to quantization,
numerical stability of reconstruction and the fact that they give
greater freedom to capture important signal characteristics
\cite{D,GKV}.  Recently, Goyal, Kova\v{c}evi\'c
 and Vetterli \cite{GKV} (see also \cite{Go,GK,GKV2,GKVT}) proposed
using the redundancy of frames to mitigate the losses in packet
based transmission systems such as the internet.  These systems
transport packets of data from a ``source'' to a ``recipient''.
These packets are sequences of information bits of a certain length
surrounded by error-control, addressing and timing information that
assure that the packet is delivered without errors.  It accomplishes
this by not delivering the packet if it contains errors.  Failures here
are due primarily to buffer overflows at intermediate nodes in the
network.  So to most users, the behavior of a packet network is not
characterized by random loss but rather by {\it unpredictable transport
time}.  This is due to a protocol, invisible to the user, that
retransmits lost or damaged packets.  Retransmission of packets takes
much longer than the original transmission and in many applications
retransmission of lost packets is not feasible.  If a lost packet
is independent of the other transmitted data, then the information is truly
lost.  But if there are dependencies between transmitted packets,
one could have partial or complete recovery despite losses.  This leads
us to consider using frames for encoding.  But which frames?
In this setting,
when frame coefficients are lost we call them {\bf erasures}.
It was shown in \cite{GKK} that an equal norm frame minimizes
mean-squared error in reconstruction with erasures
if and only if it is tight.   So
a fundamental question is to identify the optimal classes of
equal norm Parseval frames for doing reconstruction with erasures.
Since the lower frame bound of a family of vectors determines the computational
complexity of reconstruction, it is this constant we need to control.
Formally, this is a max/min problem which looks like:

\begin{problem}\label{Pr}
Given natural numbers $k,K$ find the class of equal norm Parseval
frames $\{f_i\}_{i=1}^{Kn}$ in ${\ell}_2^n$ which maximize the minimum below:
$$
min\ \{A_J : J\subset \{1,2,\ldots n\},\  |J|=k,\ A_J \ \mbox{the
lower frame bound of}\ \{f_i\}_{i\in J^c}\}.
$$
\end{problem}

This problem has proved to be very difficult.  We only
 have a complete solution to the problem
for two erasures \cite{BP,CK,HP}.  Recently, Bodemann and Paulsen
\cite{BP} have given sharp error bounds for an arbitrary number
of erasures and, more importantly, have characterized when we have
equality in these bounds.  In some settings, this proves that
equal norm tight frames are optimal.  Vershynin \cite{V2} 
shows that for any n-dimensional frame, 
any source can be linearly reconstructed from 
only n log n randomly chosen frame coefficients, 
with a small error and with high probability. Thus 
every frame expansion withstands random erasures better 
(for worst case sources) than the orthogonal basis expansion, 
for which the n log n bound is attained.
It was hoped that some special cases of the problem would be
more tractable and serve as a starting point for the classification
since the frames we are looking for are contained in this class.

\begin{conjecture}\label{CC}
There exists an $\epsilon >0$ so that for large $K$,
for all $n$ and all equal norm Parseval frames $\{f_i\}_{i=1}^{Kn}$
for ${\ell}_2^n$, there is a $J\subset \{1,2,\ldots ,Kn\}$ so that
both $\{f_i\}_{i\in J}$ and $\{f_i\}_{i\in J^c}$ have lower frame
bounds which are greater than $\epsilon$.
\end{conjecture}

The ideal situation would be for Conjecture \ref{CC} to hold for
all $K\ge 2$.
In order for $\{f_i\}_{i\in J}$ and $\{f_i\}_{i\in J^c}$ to both be
frames for ${\ell}_2^n$, they at least have to span ${\ell}_2^n$.  So the first
question is whether we can partition our frame into spanning sets.
This will follow from the Rado-Horn theorem \cite{H,R}.
For a generalization of the theorem see \cite{CKS}.

\begin{theorem}[Rado-Horn]\label{radohorn} Let $I$ be a
finite or countable index set and let
$\{f_i\}_{i\in I}$ be a collection of vectors in a vector
space.  There is a partition $\{A_j\}_{j=1}^{r}$
such that for each $j=1,2,\ldots , r$, $\{f_i\}_{i\in A_j}$
is linearly independent if and only if for all finite $J \subset I$
\begin{equation}\label{maineq}
\frac {|J|} {\dim\mathrm{span}\{f_i\}_{i\in J}} \le r.
\end{equation}
\end{theorem}

The terminology ``Rado-Horn Theorem" was introduced, to our knowledge,
in the paper \cite{B}.  This theorem has had several interesting
applications in analysis,  for one, a characterization of
Sidon sets in $\Pi_{k = 1}^\infty \Z_p$ due to Bourgain and Pisier
\cite{B, P}.  Another application is a proof that the Feichtinger
Conjecture
is equivalent Conjecture \ref{Conj101} \cite{CKS}.
In \cite{CT} it was shown that
the Rado-Horn Theorem will decompose our frames for us.

\begin{proposition}\label{1P}
Every equal norm Parseval frame $\{f_i\}_{i=1}^{Kn}$ for ${\ell}_2^n$ can
be partitioned into $K$ linearly independent spanning sets.
\end{proposition}

{\bf Proof}:  If $J\subset \{1,2,\ldots ,Kn\}$, let $P_J$ be the
orthogonal projection of ${\ell}_2^n$ onto $\mathrm{span}\{f_i\}_{i\in J}$.  Since
$\{f_i\}_{i=1}^{Kn}$ is an equal norm
Parseval frame (see Section \ref{FT})
$\sum_{i=1}^{Kn}\|f_i\|^2 = Kn\|f_1\|^2 = n$.  Now,
$$
dim(\mathrm{span}\{f_i\}_{i\in J}) = \sum_{i=1}^{Kn}\|P_Jf_i\|^2
\ge \sum_{i\in J}\|P_Jf_i\|^2 = \sum_{i\in J}\|f_i\|^2 =
\frac{|J|}{K}.
$$
So the Rado-Horn conditions hold with constant $r=K$.  If we divide
our family of Kn vectors into K linearly independent sets, since
each set cannot contain more than $n$-elements, it follows that each
has exactly $n$-elements.
\qed

If we are going to be able to erase arbitrary $k$-element subsets of our
frame, then the frame must be a union of erasure sets.  So a
generalization of Conjecture \ref{CC} which is a class containing
the class given in Problem \ref{Pr} is

\begin{conjecture}\label{CCC}
There exists $\epsilon >0$ and a natural number $r$ so that for all
large $K$ and all equal norm
Parseval frames $\{f_i\}_{i=1}^{Kn}$
in ${\ell}_2^n$
there is a partition $\{A_j\}_{j=1}^{r}$ of $\{1,2,\ldots, Kn\}$
so that for all $j=1,2,\ldots,r$
the Bessel bound of $\{f_i\}_{i\in A_j}$ is $\le 1-\epsilon$.
\end{conjecture}

Little progress has been made on this list of
problems.  But before we discuss why, let us turn to another
setting where these problems arise.
For many years engineers
have believed that it should be possible to do signal reconstruction
without phase.  Recently, Balan, Casazza and Edidin \cite{BCE} verified
this longstanding conjecture of the signal processing community
by constructing new classes of equal norm Parseval frames.
This problem comes from a fundamental problem in speech recognition
technology called the ``cocktail party problem''.

\begin{cocktail}
We have a tape recording of a group of people talking at a cocktail party.
Can we recover each individual voice with all of its voice
characteristics?
\end{cocktail}

 As we will see, the main problem here is ``signal reconstruction
with noisy phase''.
One standard technique for removing noise from  a \textit{signal}
$f\in L^2(\R)$ is to digitalize
$f$ by sending it through the {\bf fast Fourier transform}
\cite{BCE}.  This proceedure just computes the frame coefficients
of $f$ with respect to a Gabor frame (see Section \ref{TFA}),
say $\{\langle f,f_i\rangle \}_{i\in I}$.  Next, we take the
absolute values of the frame coefficients to be {\bf processed}
and {\bf store} the phases
$$
X_i(f) =\frac{\langle f,f_i\rangle}{|\langle f,f_i \rangle |}.
$$
There are countless methods for processing a signal.  One of the
simplest is {\bf thresholding}.  This is a process of deleting
any frame coefficients whose moduli fall outside of a ``threshold
interval,'' say $[A,B]$, where $0<A<B$.  The idea is that if our frame
is chosen carefully enough then the deleted coefficients will
represent the ``noise'' in the signal.  Now it is time to
reconstruct a clear signal.  This is done by passing our signal
back through the inverse fast Fourier transform (that is, we are
inverting the frame operator).  But to do this we need phases for
our coefficients.  So we take our stored $X_i(f)$ and put them
back on the processed frame coefficients which are at this time
all non-negative real numbers.  This is where the problem arises.
If the noise in the signal is actually in the phases
(which occurs in speech recognition), then we
just put the noise back into the signal.
The way to avoid this is
to construct frames for which reconstruction can be done directly
from the absolute value of the frame coefficients and not
needing the phases.  This was
done in \cite{BCE}.

\begin{theorem}\label{T44}
For a generic real frame on ${\ell}_2^n$ with at
 least $(2n-1)$-elements,  the mapping
$\pm f\rightarrow \{|\langle f,f_i \rangle |\}_{i\in I}$
is one-to-one.

For a generic complex frame on ${\ell}_2^n$
with at least $(4n-2)$-elements,
the mapping $c f\rightarrow \{|\langle f,f_i \rangle |\}_{i\in I}$, $|c|=1$,
is one-to-one.

``Generic'' here means that the set of frames with this property
is dense in the class of all frames in the Zariski topology
on the Grassman manifold \cite{BCE}.
\end{theorem}

In the process of looking for algorithms for doing reconstruction
directly from the absolute value of the frame coefficients,
it was discovered in the real case (the complex case
is much more complicated) that the standard algorithms
failed when the vector was getting approximately half of its norm
from the positive frame coefficients and half from the negative
coefficients \cite{BCE2}.  The algorithms behave as if one of
these sets has been ``erased''.  The necessary conditions for
reconstruction without phase in \cite{BCE} help explain why.
These conditions imply that every vector in the space must be
reconstructable from either the positive frame coefficients or
the negative ones.  It is also shown in \cite{BCE2} that signal
reconstruction without phase is equivalent to a $(P_0)$ problem
with additional constraints (See equation \ref{EE1} below).
So once again we have bumped into Problem
\ref{Pr} and Conjectures \ref{CC} and \ref{CCC}.

The next theorem (from \cite{CT}) helps to explain why
all of these reconstruction
problems have proved to be so difficult.  Namely, because KS has
come into play again.

\begin{theorem}\label{TTT}
(1)  Conjecture \ref{CC} implies Conjecture \ref{CCC}.

(2)  Conjecture \ref{CCC} is equivalent to KS.
\end{theorem}

{\bf Proof}:
$(1)$:  Fix $\epsilon >0$, $r,K$ as
 in Conjecture \ref{CC}.  Let $\{f_i\}_{i=1}^{Kn}$
be an equal norm Parseval frame for an $n$-dimensional
Hilbert space $\HN$. By Theorem \ref{T3} there is an
orthogonal projection $P$ on ${\ell}_2^{Kn}$ with
$Pe_i = f_i$ for all $i=1,2,\ldots ,Kn$.  By Conjecture \ref{CC},
there is a $J\subset \{1,2,\ldots ,Kn\}$ so that $\{Pe_i\}_{i\in J}$
and $\{Pe_i\}_{i\in J^c}$ both have a lower frame bound of $\epsilon>0$.
Hence, for $f\in \HN = P({\ell}_{2}^{Kn})$,
\begin{eqnarray*}
\|f\|^2 &=& \sum_{i=1}^{n}|\langle f,Pe_i \rangle |^2
= \sum_{i\in J}|\langle f,Pe_i \rangle |^2 +
\sum_{i\in J^c}|\langle f,Pe_i \rangle |^2 \\
&\ge&  \sum_{i\in J}|\langle f,Pe_i
\rangle|^2 + {\epsilon} \|f\|^2.
\end{eqnarray*}
That is, $\sum_{i\in J} |\langle f,Pe_i \rangle |^2 \le (1-\epsilon )\|f\|^2$.
So the upper frame bound of $\{Pe_i\}_{i\in J}$ (which is the norm
of the analysis operator $(PQ_J)^{*}$ for this frame) is $\le 1-\epsilon$.
Since $PQ_J$ is the synthesis operator for this frame, we have
that $\|Q_J PQ_J\|= \|PQ_J\|^2 = \|(PQ_J)^{*}\|^2 \le 1-\epsilon$.
Similarly, $\|Q_{J^c}PQ_{J^c}\| \le 1-\epsilon$.  So Conjecture
\ref{CCC} holds for $r=2$.

$(2)$:  We will show that Conjecture \ref{CCC} implies Conjecture \ref{C5}.
  Choose an integer $K$ and an $r, \epsilon >0$ with
$\frac{1}{\sqrt{K}}< \epsilon$.
  Let $\{f_i\}_{i=1}^{M}$ be a unit norm $K$-tight frame for an
$n$-dimensional Hilbert space $\HN$.
Then (see Section \ref{FT}) $M=\sum_{i=1}^{M}\|f_i\|^2 = Kn$.
Since $\{\frac{1}{\sqrt{K}}f_i \}_{i=1}^{M}$ is an
equal norm Parseval frame,  by Theorem \ref{T3}, there is an
orthogonal projection $P$ on ${\ell}_2^{M}$
with $Pe_i = \frac{1}{\sqrt{K}}f_i$, for $i=1,2,\ldots,M$.  By
Conjecture \ref{CCC},
 we have universal $r,\epsilon>0$ and
a partition $\{A_j\}_{j=1}^{r}$ of $\{1,2,\ldots ,M\}$
so that the Bessel
bound $\|(PQ_{A_j})^{*}\|^2$
for each family $\{f_i\}_{i\in A_j}$ is $\le 1-\epsilon$.
So for $j=1,2,\ldots ,r$ and any $f\in {\ell}_2^n$ we have
\begin{eqnarray*}
\sum_{i\in A_j}|\langle f,\frac{1}{\sqrt{K}}f_i \rangle |^2
&=&
\sum_{i\in A_j}|\langle f,PQ_{A_j}e_i \rangle |^2 =
\sum_{i\in A_j}|\langle Q_{A_j}Pf,e_i\rangle |^2
\le \|Q_{A_j}Pf\|^2  \\
&\le& \|Q_{A_j}P\|^2\|f\|^2 = \|(PQ_{A_j})^{*}\|^2\|f\|^2 \le
(1-\epsilon) \|f\|^2.
\end{eqnarray*}
Hence,
$$
\sum_{i\in A_j}|\langle f,f_i\rangle |^2 \le K(1-\epsilon)\|f\|^2
= (K-K\epsilon)\|f\|^2.
$$
Since $K\epsilon >\sqrt{K}$, we have verified Conjecture \ref{C5}.

For the converse, choose $r,{\delta},\epsilon$ satisfying Conjecture
\ref{C2}.  If $\{f_i\}_{i=1}^{Kn}$ is an equal norm Parseval frame
for an $n$-dimensional Hilbert space $\HN$
with $\frac{1}{K}\le \delta$, by Theorem \ref{T3} we have an
orthogonal projection $P$ on ${\ell}_2^{Kn}$ with $Pe_i = f_i$ for
$i=1,2,\ldots , Kn$.  Since ${\delta}(P) = \|f_i\|^2 \le
\frac{1}{K}\le
\delta$ (see the proof of Proposition \ref{1P}), by Conjecture \ref{C2}
there is a partition $\{A_j\}_{j=1}^{r}$ of $\{1,2,\ldots ,
Kn\}$ so that for all $j=1,2,\ldots ,r$,
$$
\|Q_{A_j}PQ_{A_j}\|=\|PQ_{A_j}\|^2 = \|(PQ_{A_j})^{*}\|^2 \le 1- \epsilon.
$$
Since $\|(PQ_{A_j})^{*}\|^2$ is the Bessel bound for
$\{Pe_i\}_{i\in A_j}=\{f_i\}_{i\in A_j}$,
we have that Conjecture \ref{CCC} holds.
\qed

Theorem \ref{TTT} yields yet another equivalent form of KS.
That is, KS is equivalent to finding a quantitative version of
the Rado-Horn Theorem.

We end this section with a class of Conjectures which were thought
to be equivalent to KS.  But, we will show that these conjectures are 
just weak enough to have a positive solution.
There is currently a flury of activity surrounding sparse solutions to vastly
underdetermined systems of linear equations.  This has applications
to problems in signal processing (recovering signals from highly
incomplete measurements), coding theory (recovering an input
vector from corrupted measurements) and much more.  If $A$ is an
$n\times m$ matrix with $n<m$, the sparsest solution to
$Af=g$ is
\begin{equation}\label{EE1}
(P_0)\ \ \ \min_{f\in {\R}^m}\|f\|_{{\ell}_0}\ \ \ \mbox{subject to
$Af=g$},
\end{equation}
where $\|f\|_{{\ell}_0} = |\{i: f(i)\not= 0\}|$.  The problem
with $(P_0)$ is that it is NP hard in general \cite{Do,N}.
This has led researchers
to consider the ${\ell}_1$ version of the problem known as
{\it basis pursuit}.
$$
(P_1)\ \ \ \min_{f\in {\R}^m}\|f\|_{{\ell}_1}\ \ \ \mbox{subject to
$Af=g$},
$$
where $\|f\|_{{\ell}_1} = \sum_{i=1}^{m}|f(i)|$.  Building on the
groundbreaking work of Donoho and Huo \cite{DH}, it has now been
shown \cite{CRT,CTa,CTa1,CTRV,Do,DE,E,Tr} that there are classes of matrices
for which the problems $(P_0)$ and $(P_1)$ have the same unique
solutions.  Since $(P_1)$ is a convex program, it can be solved
by its classical reformulation as a linear program.  A recent approach
to these problems involves {\it restricted isometry constants}
\cite{CTa1}.  If $A$ is
a matrix with column vectors $\{v_j\}_{j\in J}$, for all $1\le S\le |J|$
we define the {\it S-restricted isometry constant} ${\delta}_S$ to be
the smallest constant so that for all $T\subset J$ with $|T|\le S$ and
for all $\{a_j\}_{j\in T}$,
$$
(1-{\delta}_S )\sum_{j\in T}|a_j|^2 \le
\|\sum_{j\in T}a_j v_j \|^2 \le (1+{\delta}_S )\sum_{j\in T}|a_j|^2.
$$
The fundamental principle here is the construction of (nearly)
unit norm frames for which subsets of a fixed size are (nearly)
Parseval (or better, nearly orthogonal).  The classification of
these frames is out of our grasp at this time.  But this did
lead to a natural conjecture.

\begin{conjecture}\label{EC20}
For every $S\in \N$, for every $0< \delta <1$ and
for every unit norm tight frame $\{f_i\}_{i=1}^{\infty}$, there is a
partition $\{A_j\}_{j=1}^{r}$ of $\N$ so that for all $j=1,2,\ldots ,r$,
$\{f_i\}_{i\in A_j}$ is a frame sequence with $S$-restricted isometry
constant ${\delta}_S \le \delta$.
\end{conjecture}

There is also a finite dimensional version of Conjecture \ref{EC20}.

\begin{conjecture}\label{EC21}
For every $S\in \N$ and B and every $0< \delta <1$,  there is a natural number
$r = r(\delta, S,B)$ so that for every $n$ and every
 unit norm $B$-tight frame $\{f_i\}_{i=1}^{M}$ for ${\ell}_2^n$ there is a
partition $\{A_j\}_{j=1}^{r}$ of $\{1,2,\ldots ,M\}$
so that for all $j=1,2,\ldots ,r$,
$\{f_i\}_{i\in A_j}$ is a frame sequence with $S$-restricted isometry
constant ${\delta}_S \le \delta$.
\end{conjecture}

We would like to invoke Remark \ref{R1} here to see that Conjectures
\ref{EC20} and \ref{EC21} are equivalent.  It is easily seen by that
remark that Conjecture \ref{EC20} implies Conjecture \ref{EC21}.
Unfortunately, this approach
does not directly work for the converse since we are working with unit norm
tight frames and if we take finite ``parts'' of these, say
$\{f_i\}_{i=1}^{M}$, then these are not tight frames.  We will not
prove in detail that these are equivalent but instead just point
out that combined with the following result of Balan, Casazza,
Edidin and Kutyniok \cite{BCEK},
Remark \ref{R1} will work.

\begin{theorem}\label{ET3}
If $\{f_i \}_{i\in I}$ is a unit norm Bessel sequence in
$\HN$ with Bessel bound $B$, then there is a unit norm
family $\{g_j\}_{j\in J}$ so that $\{f_i\}_{i\in I}\cup
\{g_j \}_{j\in J}$ is a unit norm tight frame for $\HN$ with tight
frame bound ${\lambda}\le B+2$.
\end{theorem}

\begin{remark}
By Theorem \ref{ET3} and Remark \ref{R2}, Conjectures \ref{EC20} and
\ref{EC21} are equivalent to the same conjectures with ``unit norm
tight frame'' replaced by ``unit norm frame'' or replaced by
``unit norm Bessel sequence''.
\end{remark}

  A particularly interesting place to
look for frames with good restricted isometry constants
is in $L^2 [0,1]$.  

\begin{conjecture}\label{EC5}
For every measurable set $E\subset [0,1]$
with $0<|E|$, for every natural number $S$ and
for every $0< \delta <1$
there is a
a partition $\{A_j\}_{j=1}^{r}$ of $\Z$ so that for every
$j=1,2,\ldots ,r$ the family
$$
\{\frac{1}{\sqrt{|E|}}\mathrm{e}^{2\pi\mathrm{i}nt}{\chi}_E\}_{n\in A_j},
$$
is a frame sequence with S-restricted isometry constant
${\delta}_S \le \delta$.
\end{conjecture}

These conjectures deal directly with the frame.  If we want to
deal with the columns of the frame vectors we have the following
conjecture.

\begin{conjecture}\label{W6}
Let $\{f_i\}_{i=1}^{M}$ be a frame for ${\ell}_2^n$ with frame
bounds $A,B$.  Let $\{e_k\}_{k=1}^{n}$ be the unit vector basis
of ${\ell}_2^n$ and assume the column vectors $\{v_k\}_{k=1}^{n}$ are
norm one.  That is, assume $\|v_k\|^2 = \sum_{i=1}^{M}
|\langle f_i ,e_k\rangle |^2 = 1$ for every $1\le k\le n$.  
For every $\delta >0$ and for every
$S\le n$ there exists a natural number $r = r(\delta ,S,B)$
and a partition $\{A_j\}_{j=1}^{r}$ of $\{1,2,\ldots ,n\}$ 
so that for all $j=1,2,\ldots ,r$,
the family $\{v_k\}_{k\in A_j}$ has $S$-restricted isometry constant
${\delta}_S \le \delta$.
\end{conjecture}

We do not need to assume the column vectors are norm one in Conjecture
\ref{W6} but rather that they are within $\epsilon$ of being one.
The following result of Casazza, Kutyniok and Lammers \cite{CKL}
yields that Conjecture \ref{W4} is equivalent to Conjecture \ref{W6}.

\begin{theorem}
A family $\{f_i\}_{i=1}^{M}$ is a frame for ${\ell}_2^n$ with
frame bounds $A,B$ if and only if the column vectors of the
frame vectors form a Riesz basic sequence in ${\ell}_2^M$
with Riesz basis
bounds $\sqrt{A},\sqrt{B}$.
\end{theorem}

It is immediate from the $R_{\epsilon}$-Conjecture
(See section \ref{HST}) that a positive solution to KS would
imply a positive solution to Conjectures \ref{EC20}, \ref{EC21},
\ref{EC5}, and \ref{W6}.  Actually, all these conjectures are true as
we will now see.  For this we need to recall a result of 
Berman, Halpern, Kaftal and Weiss \cite{BHKW}.

\begin{theorem}\label{WKHB}
There is a natural number $r=r(B)$ satisfying the following.
Let $(a_{ij})_{i,j=1}^{n}$ be a self-adjoint matrix with non-negative
entries and with zero diagonal so that
$$
\sum_{m=1}^{n} a_{im} \le B,\ \ \mbox{for all $i = 1,2,\ldots ,n$}. 
$$
Then for every
 $r\in \N$ there is a partition $\{A_{j}\}_{j=1}^{r}$ of
$\{1,2,\ldots ,n\}$ so that for every $j=1,2,\ldots ,r$,
\begin{equation}\label{ETE}
\sum_{m\in A_j} a_{im} \le \frac{1}{r} \sum_{m\in A_{\ell}}a_{im},\ \ 
\mbox{for every $i\in A_j$ and $\ell \not= j$}.
\end{equation} 
\end{theorem}

Now we can prove our conjectures hold true.

\begin{theorem}\label{TP1}
Conjectures \ref{EC20}, \ref{EC21},
\ref{EC5}, and \ref{W6} have a positive solution. 
\end{theorem}

{\it Proof}:
We will prove Conjecture \ref{EC21} for unit norm Bessel sequences.
Let $\{f_i\}_{i=1}^{M}$ be a unit norm $B$-Bessel sequence in
${\ell}_2^n$.  Let $H$ be the matrix 
$$
 H =  \left ( \langle f_i,f_m \rangle  \right )_{i,m=1}^{M}.
$$
For each $1\le i\le M$,
$$
\sum_{m=1}^{M}|\langle f_i ,f_m \rangle |^2 \le B.
$$
Fix a $k\in \N$ with $\sqrt{\frac{BS}{k}}\le \delta$
and fix $S$ as in Conjecture
\ref{EC21}.  By Theorem \ref{WKHB}, there is a natural number
$r=r(B,S,k) \in \N$ and a partition $\{A_j\}_{j=1}^{r}$ so that
$H-D(H)$ satisfies Equation \ref{ETE}.  Fix $1\le j\le r$, let
$T\subset A_j$ with $|T|\le S$ and let $(a_{i})_{i\in T}$ be
scalars.  Then,
\begin{eqnarray*}
\|\sum_{i\in T} a_if_i\|^2 &=& \sum_{i\in T}|a_i|^2 \|f_i\|^2 +
\sum_{i\not= m \in T}a_i\overline{a_m}\langle f_i ,f_m\rangle \\
&\le& \sum_{i\in T}|a_i|^2 + \left ( \sum_{i\not= m\in T}
|a_i|^2 |a_m|^2 \right )^{1/2} \left ( \sum_{i\not= m \in T}
|\langle f_i ,f_m \rangle |^2 \right )^{1/2} \\
&\le& \sum_{i\in T}|a_i|^2 + \left [ \left ( \sum_{i\in T} |a_i|^2 
\right )^{2} \right ]^{1/2} \left (\sum_{i\in T}\frac{B}{k} \right )^{1/2}\\
&\le& \sum_{i\in T}|a_i|^2 +  \left ( \sum_{i\in T} |a_i|^2 
\right ) \sqrt{\frac{BS}{k}} \\
&\le& \sum_{i\in T}|a_i|^2 +  \delta \sum_{i\in T} |a_i|^2 \\
&=& (1+\delta)\sum_{i\in T}|a_i|^2.
\end{eqnarray*}
Similarly we have
$$
\|\sum_{i\in T} a_if_i\|^2 \ge (1-\delta)\sum_{i\in T}|a_i|^2.
$$
It follows that $\{f_i\}_{i\in A_j}$ has $S$-restricted
isometry constant ${\delta}_{S} \le \delta$.  
\qed

What this section is trying to tell us is the following.  In applied
mathematics and engineering problems we are generally looking for
the best examples we can find to use in practice.  However, if we
instead ask the question:  {\it Let's classify all objects which
satisfy our requirements}, then we have entered the world of the
deepest unsolved problems in pure mathematics.

\section{Towards a counter-example to Kadison-Singer}\label{CE}
\setcounter{equation}{0}

In this section we will give some more equivalents of Kadison-Singer
which lend themselves to viable approaches for constructing a
counterexample to KS.  Throughout this section we will use the
notation:
\vskip10pt
\noindent {\bf Notation}:  If $E\subset I$ we let $P_E$ denote the
orthogonal projection of ${\ell}_2(I)$ onto ${\ell}_2(E)$.
Also, recall that we write $\{e_i\}_{i\in I}$ for the standard
orthonormal basis for ${\ell}_2(I)$.
\vskip10pt
For results on frames, see Section \ref{FT}.

\begin{definition}
A subspace $\H$ of ${\ell}_2(I)$ is {\it A-large} for $A>0$
if it is closed and for each $i\in I$, there is a vector
$f_i \in \H$ so that $\|f_i\|=1$ and $|f_i(i)|\ge A$.  The
space $\H$ is {\it large} if it is A-large for some $A>0$.
\end{definition}

We are going to classify PC in terms of A-large subspaces of ${\ell}_2(I)$.
To do this we need some preliminary results.

\begin{lemma}
 Let $T^{*}:\H \rightarrow
{\ell}_2(I)$ be the analysis operator for a frame
$\{f_i\}_{i\in I}$ for $\H$ and let $P$ be the orthogonal
projection of ${\ell}_2(I)$ onto $\H$.  Then $\{Pe_i\}_{i\in I}$
is a Parseval frame for $T^{*}(\H)$ which is equivalent to $\{f_i\}_{i\in I}$.
\end{lemma}

{\it Proof}:
Note that $\{Pe_i\}_{i\in I}$ is a Parseval frame (Theorem \ref{T3})
with synthesis operator $P$ and analysis operator $T_{1}^{*}$ satisfying
$T_{1}^{*}(\H) = P({\ell}_2(I)) = T^{*}(\H)$.  By Proposition \ref{FTP10},
$\{Pe_i\}_{i\in I}$
is equivalent to $\{f_i\}_{i\in I}$.
\qed

\begin{proposition}\label{CEP2}
Let $\H$ be a subspace of ${\ell}_2(I)$.  The following are equivalent:

(1)  The subspace $\H$ is large.

(2)  If $P$ is the orthogonal projection of ${\ell}_2(I)$ onto
$\H$ then there is an $A>0$ so that $\|Pe_i\|\ge A$, for all $i\in I$.

(3)  The subspace $\H$ is the range of the analysis operator of
some bounded frame $\{f_i\}_{i\in I}$.
\end{proposition}

{\it Proof}:
$(1)\Rightarrow (2)$:  Suppose $\H$ is large.  So, there exists an
$A>0$ such that for each $i\in I$, there exists a vector $f_i \in \H$
with $\|f_i\|=1$ and $|f_i(i)|\ge A$.  Given the projection $P$ of
(2) we have
$$
A \le |f_i(i)| = |\langle e_i,f_i\rangle | = |\langle Pe_i,f_i\rangle |
\le \|Pe_i\|\|f_i\| = \|Pe_i\|.
$$

$(2)\Rightarrow (3)$:  By (2), $\{Pe_i\}_{i\in I}$ is a bounded
sequence which is a Parseval frame by Theorem \ref{T3} and having
$\H$ as the range of its analysis operator.

$(3) \Rightarrow (1)$:  Assume $\{f_i\}_{i\in I}$ is a bounded
frame for a Hilbert space $\K$ with analysis operator $T^{*}$ and
$T^{*}(\K) = \H$.  Now, $\{Pe_i\}_{i\in I}$ is a Parseval frame
for $\H$ which is the range of its own analysis operator.  Hence,
$\{f_i\}_{i\in I}$ is equivalent to $\{Pe_i\}_{i\in I}$ by
Proposition \ref{FTP10}.  Since $\{f_i\}_{i\in I}$ is bounded,
so is $\{Pe_i\}_{i\in I}$.  Choose $A>0$ so that $A\le \|Pe_i\| \le 1$,
for all $i\in I$.  Then
$$
A\le |\langle Pe_i,Pe_i \rangle | = |\langle Pe_i,e_i\rangle | =
|Pe_i (i)|.
$$
So $\H$ is a large subspace.
\qed

Now we need to learn how to decompose the range of the analysis operator
of our frames.

\begin{definition} \label{D:rD}
A closed subspace $\H$ of ${\ell}_2(I)$ is {\it r-decomposable}
if for some natural number $r$ there exists a partition
$\{E_j\}_{j=1}^{r}$ of $I$
 so that $P_{E_j}(\H) = {\ell}_2(E_j)$, for all
$j=1,2,\ldots ,r$.  The subspace $\H$ is {\it finitely decomposable}
if it is r-decomposable for some r.
\end{definition}

For the next proposition we need a small observation.

\begin{lemma}\label{TFAL1}
Let $\{f_i\}_{i\in I}$ be a Bessel sequence in $\H$ having
synthesis operator $T$ and
analysis operator $T^{*}$, let
$E\subset I$, and let $\{f_i\}_{i \in E}$ have analysis operator $(T|_E)^{*}$.  Then
$$
P_ET^{*} = (T|_E)^{*}.
$$
\end{lemma}

{\it Proof}:
For all $f\in \H$,
$$
P_E T^{*}(f) = P_E\left ( \sum_{i\in I}\langle f,f_i\rangle e_i \right )
= \sum_{i\in E}\langle f,f_i \rangle e_i = (T|_E)^{*}(f).
$$
\qed

We now have

\begin{proposition}\label{CEP1}
A frame $\{f_i\}_{i\in I}$ for $\K$ satisfies the Feichtinger Conjecture
if and only if $\H = T^{*}(\K)$ is finitely decomposable.
\end{proposition}

{\it Proof}:
We can partition $I$ into $\{E_j\}_{j=1}^{r}$ so that each
$\{f_i\}_{i\in E_j}$ is a Riesz basic sequence if and only if
(see the discussion after Theorem \ref{FTTT})
$(T|_{E_j})^{*}$ is onto for every
$j=1,2,\ldots ,r$ if and only if (by Lemma \ref{TFAL1})
$P_{E_j}T^{*}$ is onto for all $j=1,2,\ldots ,r$.
\qed

Now we can put this altogether.

\begin{theorem}
The following are equivalent:

(1)  The Kadison-Singer Problem.

(2)  Every large subspace of ${\ell}_2(I)$ is finitely decomposable.

(3)  For every $0<A<1$ there is a natural number $r = r(A)$ so
that every A-large subspace of ${\ell}_2(I)$ is r-decomposable.
\end{theorem}

{\it Proof}:
$(1)\Leftrightarrow (2)$:  This is immediate from
Propositions \ref{CEP2} and \ref{CEP1}.

$(2)\Rightarrow (3)$:  We prove the contrapositive.  If (3)
fails, then there is an $0<A<1$ and a sequence of subspaces
${\H}_j$ $j=1,2,\ldots$ of ${\ell}_2(I)$ so that each
${\H}_j$ is A-large but not j-decomposable.  But now,
$(\oplus_{j\in \N} {\H}_j )_{{\ell}_2}$ is an A-large
subspace of $(\oplus_{j\in \N} {\ell}_2(I))_{{\ell}_2}$
which fails to be decomposable.

$(3)\Rightarrow (2)$:  This is obvious.
\qed

Now we want to give quite explicit information about the existence
of certain families of vectors in every large subspace of
${\ell}_2(I)$.  We will see that this gives us an approach to
producing a counterexample to KS.

\begin{proposition}\label{CEP3}
Let $E\subset I$, and assume for every $i\in E$ there are vectors
$$
f_i = e_i + g_i \in {\ell}_2(I),
$$
where each $g_i \in {\ell}_2(E^c)$, and the collection $\{g_i\}_{i\in E}$ is a Bessel sequence.
Then, $\{f_i\}_{i\in E}$ is a Riesz basic sequence.  Moreover, if $\mathbb{K}$ is the closed span of $\{f_i\}_{i \in E}$, then $P_{E} \mathbb{K} = \ell_{2}(E)$.
\end{proposition}

{\it Proof}:  That $\{f_i\}$ is a Bessel sequence is obvious, and so $\{f_i\}$ possesses an upper basis bound.

We establish a lower basis bound.  For all sequences of scalars $\{a_i\}_{i\in E}$ we have:
\begin{eqnarray}
\| \sum_{i\in E}a_i f_i \| &=& \|\sum_{i\in E}a_i e_i + \sum_{i\in E}a_i g_i\| \notag \\
&\ge& \|\sum_{i\in E}a_i e_i \| \label{E:lb} \\
&=&  \left ( \sum_{i\in E}|a_i|^2 \right )^{1/2}, \notag
\end{eqnarray}
where the estimate in (\ref{E:lb}) follows by virtue of the orthogonality of $\sum_{i\in E} a_i e_i$ and $\sum_{i\in E} a_i g_i$.

If $\sum_{i \in E} a_i e_i \in \ell_{2}(E)$, then $\sum_{i \in E} a_i f_i \in \mathbb{K}$ and, since the $g_i$'s are supported outside of $E$,
\[  P_{E} \left( \sum_{i \in E} a_i f_i \right) = \sum_{i \in E} a_i e_i . \]
\qed


The following is a converse to Proposition \ref{CEP3}.

\begin{theorem}\label{ABC}
Let $\H$ be a closed subspace of ${\ell}_2(I)$.  The following are
equivalent:

(1)  $\H$ is finitely decomposable.

(2)  We can partition $I$ into subsets $\{E_j\}_{j=1}^{r}$
so that for every $j=1,2,\ldots ,r$ and all $i\in E_j$ we can
find vectors
$$
f_{ji} = e_i + g_{ji}\in \H,
$$
so that $g_{ji}\in\mathrm{span}_{k\notin E_j}e_k$ and $\{g_{ji}\}_{i\in E_j}$ is Bessel.
\end{theorem}

{\it Proof}:
$(1)\Rightarrow (2)$:  Assume $\H$ is finitely decomposable.
Let $\{E_j\}_{j=1}^{r}$ be a partition of $I$ which satisfies Definition \ref{D:rD}.
Fix $1\le j\le r$.  Since $P_{E_j}:\H \rightarrow {\ell}_2(E_j)$
is bounded, linear and onto, it follows that $P^{*}_{E_j}$
is an (into) isomorphism.  Therefore, $\{P^{*}_{E_j}e_i\}_{i\in E_{j}}$
is a Riesz basis for its span.  Let $\{f_{ji}\}_{i\in E_j}$
be the dual functionals for this Riesz basis.
Now, for all $i,\ell\in E_j$
we have
$$
{\delta}_{{\ell}i} = \langle P^{*}_{E_j}e_{\ell} ,f_{ji}\rangle
= \langle e_{\ell}, P_{E_j}f_{ji}\rangle.
$$
It follows that $f_{ji}(\ell) = 0$ if $i\in E_j$ and $i\not= \ell$,
and $f_{ji}(i) = 1$.  Hence, $f_{ji} = e_i + g_{ji}$ where
$g_{ji}\in\mathrm{span}_{i\notin E_j}e_i$.  Finally, since $\{f_{ji}\}_{i\in E_j}$
is a Riesz basis, it follows that $\{g_{ji}\}_{i\in E_j}$ is Bessel.

$(2) \Rightarrow (1)$:  This is immediate from Proposition
\ref{CEP3}.
\qed

\begin{remark}
The vectors which arise in Theorem \ref{ABC} are unique.  That is,
if
$$
\hat{f}_{ki}= e_i + \hat{g}_{ki} \in \H,
$$
(even without any assumption that the $\{\hat{g}_{ki}\}_{i\in E_j}$ are
Bessel), then $\hat{f}_{ki} = f_{ki}$, for all $i\in E_j$.  This
follows from the fact that $P_{E_j}$ is invertible on the range of
$P^{*}_{E_j}$.
\end{remark}

We will now discuss why we believe that Theorem \ref{ABC} gives
a viable approach to constructing a counterexample to KS.  Basically,
we want to construct a sequence of vectors $f_i$ in
${\ell}_2(\N)$ each having at least one big coordinate but so
that whenever we partition $\N$ into a finite number of sets,
one of these sets has sufficient density to guarantee that the vectors
$f_{ki}$ cannot be Bessel.  As we have seen, these vectors are unique.
To get the vectors $f_{ki}$ we have to ``row reduce'' the
$\{f_i\}_{i\in E_j}$ accross the coefficients of $E_j$.  If the
$\{f_i\}$ are chosen appropriately, we believe that this row
reduction process will leave us with $g_{ki}$ which are no
longer Bessel.

This may seem esoteric, but all of this was built
on existing deep constructions in the Banach space
approximation property due to
Szankowski \cite{Sz3,Sz1,Sz2} (see \cite{C5,LT}).
A look at \cite{Sz3} shows that
Szankowski constructs vectors with 6 ones in each vector 
(and this is their only support) in such
a way that when these vectors sit in ${\ell}_p$,
$p\not= 2$ in a careful way, they span a sublattice failing the
approximation property.  Of course, Hilbert spaces have the
approximation property.  But our above propositions show that
the Kadison-Singer Problem is asking for a specialized class
of operators to give the required approximation.  That is,
Kadison-Singer is a {\it restricted} approximation property for
${\ell}_2$.  What we need to do is add a bounded set of
vectors onto the Szankowski construction so that the set is
Bessel, but when we do the required row reduction to get the
vectors in Theorem \ref{ABC}, we end up with a non-Bessel
sequence $\{g_{ki}\}_{i\in E_j}$ for one of the $j=1,2,\ldots ,r$.

\end{document}